\documentclass[a4paper,10pt]{article}
\usepackage[utf8x]{inputenc}
\usepackage{tracefnt,amsmath,tabu,array}
\usepackage{amssymb,graphicx,setspace,amsfonts,amsbsy}
\usepackage{pifont,latexsym,ifthen,amsthm,rotating,calc,textcase,booktabs,cancel,slashed}

\newtheorem{theorem}{Theorem}[section]
\newtheorem{lemma}[theorem]{Lemma}
\newtheorem{corollary}[theorem]{Corollary}
\newtheorem{proposition}[theorem]{Proposition}

\newtheorem{remark}[theorem]{Remark}
\newcommand{\filledbox}{\leavevmode
  \hbox to.77778em{%
  \hfil\vbox to.675em{\hrule width.6em height.6em}\hfil}}

\newcommand{\Rm}{{\mathbb R}}

\newcommand{\eps}{\varepsilon}

\begin{document}
\tabulinesep=1.0mm
\title{The radiation theory of radial solutions to 3D energy critical wave equations}

\author{Ruipeng Shen\\
Centre for Applied Mathematics\\
Tianjin University\\
Tianjin, China
}

\maketitle

\begin{abstract}
 In this work we consider a wide range of energy critical wave equation in 3-dimensional space with radial data. We are interested in exterior scattering phenomenon, in which the asymptotic behaviour of a solutions $u$ to the non-linear wave equation is similar to that of a linear free wave $v_L$ in an exterior region $\{x: |x|>R+|t|\}$, i.e. 
 \[
  \lim_{t\rightarrow \pm \infty} \int_{|x|>R+|t|} (|\nabla(u-v_L)|^2 + |u_t-\partial_t v_L|^2) dx = 0.
 \]
 We classify all such solutions for a given linear free wave $v_L$ in this work. We also give some applications of our theory on the global behaviours of radial solutions to this kind of equations. In particular we show the scattering of all finite-energy radial solutions to the defocusing energy critical wave equations. 
\end{abstract}

\section{Introduction}

\subsection{Background}

In this work we consider the 3-dimensional nonlinear wave equation 
\[
 \left\{\begin{array}{ll} \partial_t^2 u - \Delta u = F(x,t,u), & (x,t) \in \Rm^ 3\times \Rm; \\ (u,u_t)|_{t=0} = (u_0,u_1)\in \dot{H}^1\times L^2(\Rm^3). & \end{array} \right.\qquad (CP1) 
\]
in the radial case with an energy critical nonlinear term $F(x,t,u)$. The details of our assumption on $F$ will be given later. The most frequently used nonlinear terms of this kind are the power-type ones, i.e. $F(x,t,u) = \pm |u|^4 u$. The $\pm$ signs correspond to focusing ($+$) and defocusing ($-$) case, respectively. These equations with power-type nonlinearity satisfy the following rescaling invariance property: if $u$ is a solution, then $\lambda^{-1/2} u(x/\lambda, t/\lambda)$ solves the same equation for any positive constant $\lambda\in \Rm$. This rescaling invariance property plays an important role in the study of these equations. It was proved in the last few decades of 20th century that any solution $u$ to the defocusing equation can be defined for all time and scatter in both two time directions, i.e. there exist two free waves $v_L^\pm$ so that 
\[
 \lim_{t\rightarrow \pm \infty} \int_{\Rm^3} |\nabla_{t,x}(u-v_L^\pm)(x,t)|^2 dx = 0.
\]
Here $\nabla_{t,x} = (\partial_t, \nabla_x)$. Please see, \cite{enscatter1, enscatter2, ss2, struwe}, for instance. The focusing case is much more subtle. Solutions may scatter, blow up in finite time, stay unchanged for all time, or exhibit more complicated behaviours. Please refer to \cite{se, kenig}, for example. We are particularly interested in the asymptotic behaviours of solutions to these equations, especially the exterior scattering. Exterior scattering means that a solution $u$ to a nonlinear wave equation approaches a linear free wave, i.e. a solution $v_L$ to the homogeneous linear wave equation in the exterior region $\{x: |x|>R+|t|\}$, more precisely
\[
 \lim_{t\rightarrow \pm \infty} \int_{|x|>|t|+R} |\nabla_{t,x} (u-v_L)(x,t)|^2 dx = 0. 
\]
The exterior scattering is a very common phenomenon. By finite speed of propagation and standard center cut-off techniques, one may show that any global solution to (CP1) always scatters in the exterior region $\{x: |x|>|t|+R\}$, as long as the radius $R$ is sufficiently large. We believe that by discussing the exterior scattering behaviour systemically, we will be able to gain information about the global behaviour of solutions in the whole space. A lot of results in this work are previously known (see \cite{kenig, enscatter2}) for the wave equation with power-type nonlinearity $\pm |u|^4 u$. Please note that many already known results depends on the rescaling invariance of solutions to these power-type equations. Our argument in this work, however, does not depend on the rescaling invariance, thus applies to wave equations with more general nonlinearities.  

\subsection{Basic conceptions}

Before we may discuss the main topics of this work and give our main results, we first give our assumptions on the non-linear term $F$ in (CP1), then introduce a few notations and basic conceptions. 

\paragraph{Assumptions} We assume that the nonlinear term always satisfies (A1), (A2) in this work. In the last few sections we also assume (A3) and (A4). 
\begin{itemize}
 \item[(A1)] The function $F(x,t,u)$ is a radial function of $x \in \Rm^3$; 
 \item[(A2)] The function $F(x,t,u)$ satisfies the following inequalities. Here $\gamma$ is a positive constant.
 \begin{align*}
  &|F(x,t,u)| \leq \gamma |u|^5;& &|F(x,t,u)-F(x,t,v)| \leq  \gamma |u-v|(|u|^4+|v|^4);&
 \end{align*}
 \item[(A3)] The function $F(x,t,u)$ is independent of time $t$; 
 \item[(A4)] The function $F(x,t,u)$ satisfies the inequality
 \begin{align*}
  |F(x,t,u+v)-F(x,t,u)-F(x,t,v)| \leq \gamma |uv|(|u|^3+|v|^3). 
 \end{align*}
\end{itemize} 

\paragraph{Notations} We first define the following regions in $\Rm^3 \times \Rm$ for radii $R \geq r \geq 0$:
\begin{align*}
 & \Omega_R = \{(x,t): |x|>|t|+R\};& &\Omega_{r,R} = \{(x,t): |t|+r<|x|<|t|+R\};& &\tilde{\Omega}_R = \Omega_{R,2R}.& 
\end{align*}
We also use the notations $\chi_R, \chi_{r,R}, \tilde{\chi}_R$ for the corresponding characteristic functions of these regions. If $u$ is defined in $\Rm^3 \times I$, we define the space-time norm
\[
 \|u\|_{Y(I)} = \|u\|_{L^5 L^{10} (I\times \Rm^3)} = \left(\int_I \left(\int_{\Rm^3} |u(x,t)|^{10} dx\right)^{1/2} dt\right)^{1/5}.
\]
In this work $A \lesssim B$ means that there exists a positive constant $c$ so that $A \leq cB$. We may add subscript(s) to the notation $\lesssim$ to imply that the constant $c$ depends on the subscript(s) but nothing else. In particular we use the notation $\lesssim_1$ to emphasize that the constant $c$ is an absolute constant. The notations $\gtrsim$ and $\simeq$ are similar. We let $\mathbf{S}_L$ be the linear wave propagation operator, i.e. $\mathbf{S}_L (u_0,u_1)$ represents the linear free wave with initial data $(u_0,u_1)$. We also use the notation $\vec{u}$ for the data $(u,u_t)$ of a solution $u$ to a wave equation. Finally we use the following notation for the energy of radial data in the exterior region
\[
 \|(u_0,u_1)\|_{\mathcal{H}_R} = \|(u_0,u_1)\|_{\dot{H}^1\times L^2(\{x:|x|>R\})} = \left(\int_{|x|>R} (|\nabla u_0(x)|^2 + |u_1(x)|^2)dx\right)^{1/2}.
\]

\paragraph{The conception of exterior solutions} Let $u$ be a function defined in the exterior region 
\[
 \{(x,t): |x|>|t|+R,\, t\in (-T_1,T_2)\}.
\]
Here $T_1, T_2 > 0$ are either real numbers or $\infty$. We call $u$ an exterior solution to (CP1) in the region with initial data $(u_0,u_1)\in \dot{H}^1 \times L^2(\Rm^3)$, if and only if $\|\chi_R u\|_{Y(J)} < +\infty$ for any bounded closed interval $J \subset (-T_1,T_2)$, so that 
 \[
  u = \mathbf{S}_L (u_0,u_1) + \int_0^t \frac{\sin (t-t')\sqrt{-\Delta}}{\sqrt{-\Delta}} [\chi_R(\cdot,t') F(t', \cdot, u)] dt', \qquad |x|>R+|t|, \; t\in (-T_1,T_2).
 \]
 Please note that we add $\chi_R$ to emphasize that $u$ and $F(x,t,u)$ are only defined in the exterior region. More precisely we understand the product in the following way
 \begin{align*}
  &\chi_R u = \left\{\begin{array}{ll} u(x,t), & x\in \Omega_R; \\ 0, & x\notin \Omega_R. \end{array} \right.&
  &\chi_R F(x,t,u) = \left\{\begin{array}{ll} F(x,t,u(x,t)), & x\in \Omega_R; \\ 0, & x\notin \Omega_R. \end{array} \right.&
 \end{align*}
 Here we define the initial data in the whole space, but the solution only depends on the value of initial data in the exterior region $\{x: |x|>R\}$, by finite speed of propagation. Sometimes we say $u$ is a exterior solution without giving its initial data, this means that there exist initial data $(u_0,u_1)\in \dot{H}^1\times L^2(\Rm^3)$ so that $u$ becomes an exterior solution to (CP1) with these initial data. 
 
 \paragraph{Restriction and global extension} If $u$ is an exterior solution to (CP1) defined in the exterior region $\Omega_R$ for $t \in I$. Then finite speed of propagation shows that its restriction in the exterior region $\Omega_r$ with $r>R$ is also an exterior solution to (CP1) with the same initial data. On the other hand, the function defined for $(x,t)\in \Rm^3 \times I$ by the formula
 \[ 
   \mathbf{S}_L (u_0,u_1) + \int_0^t \frac{\sin (t-t')\sqrt{-\Delta}}{\sqrt{-\Delta}} [\chi_R(\cdot,t') F(t', \cdot, u)] dt'
 \]
coincides with $u$ in the exterior region $\Omega_R$ and solves the modified wave equation below defined in the whole space $\Rm^3$
 \[
  \partial_t^2 u - \Delta u = \chi_R F(x,t,u), \qquad (x,t) \in \Rm^3 \times I.
 \]
We call it the global extension of $u$ and still use the notation $u$ to represent it. Please note that the global extension does depend on the choice of initial data in the interior region $\{x: |x|<R\}$. The Strichartz estimates then show that the global extension satisfies ($J \subset I$ is a bounded closed interval)
 \begin{align*}
  &\|u\|_{Y(J)} < +\infty;& &(u, u_t) \in C(I; \dot{H}^1\times L^2(\Rm^3)).&
 \end{align*}
 
 \paragraph{Local theory} A combination of the Strichartz estimates and a fixed-point argument immediately gives the following results for local theory of exterior solutions. The argument is similar to those given in \cite{loc1, ls} and somewhat standard. We omit the details except for the last result below (scattering criterion, see Lemma \ref{equivalent to exterior scattering})
 \begin{itemize}
  \item Well-posedness: Given any initial data $(u_0,u_1) \in \dot{H}^1 \times L^2(\Rm^3)$ and $R\geq 0$, there exists a unique exterior solution to (CP1) in the region $\Omega_R$ with a maximal lifespan $(-T_-,T_+)$.
  \item Small data scattering: There exists a constant $\delta = \delta(\gamma) >0$, so that if the initial data satisfy $\|\chi_R \mathbf{S}_L (u_0,u_1)\|_{Y(\Rm)} < \delta$, then the corresponding exterior solution $u$ to (CP1) is defined for all $t\in\Rm$, so that $\|\chi_R u\|_{Y(\Rm)} \leq 2 \|\chi_R \mathbf{S}_L (u_0,u_1)\|_{Y(\Rm)}$.
  \item Finite time blow-up criterion: If the exterior solution $u$ blow up in finite time, i.e. $T_+ < +\infty$, then $\|\chi_R u\|_{Y([0,T_+))} = +\infty$. The situation in negative time direction is similar. 
  \item Scattering criterion: The exterior solution $u$ in $\Omega_R$ scatters in the positive time direction, i.e. $T_+ =  +\infty$ and there exists a linear free wave $u_L^+$ so that 
  \[
   \lim_{t\rightarrow +\infty} \|\vec{u} - \vec{u}_L^+\|_{\dot{H}^1 \times L^2(\{x: |x|>|t|+R\})} = 0, 
  \]
  if and only if $\|\chi_R u\|_{Y([0,T_+))} < +\infty$. 
 \end{itemize}
 
\paragraph{Asymptoticly equivalent solutions} Let $u,v \in C(\Rm; \dot{H}^1\times L^2(\Rm^3))$. We say $u$ and $v$ are $R$-weakly asymptotically equivalent if and only if they satisfy
\[
 \lim_{t\rightarrow \pm \infty} \int_{|x|>|t|+R} |\nabla_{t,x} (u-v)|^2 dx = 0. 
\]
In particular, if $R=0$, we say that $u$ and $v$ are asymptotically equivalent. Since the limit above solely depends on the values of $u,v$ in the exterior region $\Omega_R$, the definition above also applies to functions $u,v$ defined in any region containing $\Omega_R$, as long as they coincides with $C(\Rm; \dot{H}^1\times L^2(\Rm^3))$ functions in the exterior region $\Omega_R$. We say that $u$ and $v$ are weakly asymptotically equivalent if they are $R$-weakly asymptotically equivalent  for some radius $R\geq 0$. We are mainly interested in the case when a solution to a nonlinear wave equation is weakly asymptotically equivalent to a linear free wave. This phenomenon is usually called exterior scattering. A very special case is a solution $u$ that is weakly asymptotically equivalent to the zero solution, i.e. a solution satisfying
\[
 \lim_{t\rightarrow \pm \infty} \int_{|x|>|t|+R} |\nabla_{t,x} u(x,t)|^2 dx = 0.
\]
These solutions are usually called (weakly) non-radiative solutions, which have been extensively studied and play an important role in the channel of energy method. The channel of energy have become an important tool in the study of global behaviours of solutions to non-linear wave equation in the past decade. Please see \cite{tkm1} and \cite{channel5d}, for example. The applications of the channel of energy include soliton resolution of focusing energy critical wave equations (see \cite{se}, \cite{oddhigh}) and conditional scattering theory (see \cite{dkm2}, \cite{shen2}).  

\subsection{Goal and main results}

In this work we give a systematic theory concerning the properties and structure of all weakly asymptotically equivalent solutions to (CP1) of a given free wave, in the 3-dimensional case and the radial setting. We show that all the exterior solutions to (CP1) which are weakly asymptotically equivalent to a given radial free wave with a finite energy form a family of a single parameter. We also give the asymptotic behaviour of these solutions in details. Furthermore, if the nonlinear term satisfies (A3) and (A4), the global behaviours of bounded solutions to (CP1) heavily depend on those of the non-radiative solutions. Our main results include: 

\begin{theorem}[One-parameter family] \label{main 1}
Assume that $F(x,t,u)$ satisfies (A1) and (A2). Let $w_L$ be a finite-energy radial free wave. Then there exists a one-parameter family $\{(u^\alpha, R_\alpha)\}_{\alpha \in \Rm}$ so that each pair $(u^\alpha, R_\alpha)$ satisfies either of the following
\begin{itemize}
 \item[(a)] The radial function $u^\alpha$ is an exterior solution to (CP1) in $\Omega_0$ and is asymptotically equivalent to $w_L$; In this case we choose $R_\alpha = 0^-$;
 \item[(b)] The radial function $u^\alpha$ is defined in $\Omega_{R_\alpha}$ with $R_\alpha \geq 0$ and $\|\chi_{R_\alpha} u^\alpha\|_{Y(\Rm)} = +\infty$ so that for any $R > R_\alpha$, $u^\alpha$ is an exterior solution to (CP1) in $\Omega_R$ and is $R$-weakly asymptotically equivalent to $w_L$.
\end{itemize}
In addition, if $u$ is a radial exterior solution to (CP1) defined in $\Omega_R$ so that $u$ is $R$-weakly asymptotically equivalent to $w_L$, then there exists a unique real number $\alpha$, so that\footnote{In the case (a), we understand $0>R_\alpha = 0^-$.} $R>R_\alpha$ and $u(x,t) = u^\alpha(x,t)$ for $(x,t) \in \Omega_R$. We call the number $\alpha$ the characteristic number of $u$. The characteristic number can also be characterized by the asymptotic behaviour of the solutions. More precisely, given $\alpha, \beta \in \Rm$, we have 
\[
 \lim_{r\rightarrow +\infty} r^{1/2} \sup_{t\in \Rm} \left\|\vec{u^\alpha} (\cdot,t) - \vec{u^\beta} (\cdot, t) - ((\alpha-\beta)|x|^{-1}, 0)\right\|_{\dot{H}^1\times L^2(\{x: |x|>|t|+r\})} = 0. 
\]
\end{theorem}

\begin{remark}
 This one-parameter family is unique up to a translation of the parameter. Namely, if $\{(u^\alpha, R_\alpha)\}_{\alpha \in \Rm}$ is such a family, then all possible families of solutions are given by the translations $\{(u^{\alpha+\beta}, R_{\alpha+\beta})\}_{\alpha \in \Rm}$. Here $\beta\in \Rm$ is a constant. Given two exterior solutions $u$ and $v$ which are weakly asymptotically equivalent to the same free wave, the difference of their characteristic numbers is independent of the choice of one-parameter family. We call this difference the characteristic number of $u$ with respect to $v$. Some free waves $w_L$ (especially those small ones, see Subsection \ref{sec: primary}) may admit natural choice of $u^0$, thus give a standard one-parameter family. For example, if $w_L = 0$, then we may choose a standard non-radiative solution $u^0 = 0$. Throughout this work, we always use this standard one-parameter family whenever we are discussing non-radiative solutions. 
\end{remark}

\begin{theorem} \label{main non-radiative}
 Assume that $F$ satisfies (A1)-(A3). Then all the radial weakly non-radiative solutions to (CP1) are independent of time. each pair $(u^\alpha, R_\alpha)$ in the one-parameter family of non-radiative solutions defined above satisfies either of the following
\begin{itemize}
 \item The solution $u^\alpha \in \dot{H}^1(\Rm^3)$ solves (CP1) in $\Rm^3 \times \Rm$ and $R_\alpha = 0^-$. 
 \item The solution $u^\alpha$ is $r$-weakly non-radiative for all $r>R_\alpha$ but we also have 
 \begin{align*}
  &\|u^\alpha\|_{\dot{H}^1(\{x: |x|>R_\alpha\})} = +\infty;& &\|\chi_{R_\alpha} u^\alpha\|_{Y(\Rm)} = +\infty.&
 \end{align*}
\end{itemize}
In addition, we always have $u^0 = 0$. If $\alpha \neq 0$, then $\|u^\alpha\|_{\dot{H}^1(\{x: |x|>R\})} \simeq_1 |\alpha| R^{-1/2}$ for $R\gtrsim_\gamma |\alpha|^2$. 
\end{theorem}

\begin{remark}
 In the case of focusing wave equation $\partial_t^2 u - \Delta u = +|u|^4 u$, all the nontrivial radial non-radiative solutions are explicitly given by the formula 
 \[
  u^\alpha = \frac{1}{\alpha} \left(\frac{1}{3} + \frac{|x|^2}{\alpha^4}\right)^{-1/2}, \qquad \alpha \neq 0.
 \]
 The nonzero radial non-radiative solutions are unique up to rescaling and a plus/minus sign. They come with the same $\dot{H}^1(\Rm^3)$ norm and are usually called ground states. 
\end{remark}

\noindent Next result shows that the global behaviours of general radial solutions heavily depends on those of the non-radiative solutions. 
\begin{theorem} \label{main global}
 Assume that $F$ satisfies (A1)-(A4) and that $\{(R_\alpha, u^\alpha)\}$ is the one-parameter family of non-radiative solutions defined above. Let $A, c_1$ be positive constants. If there exists a radius $r_\alpha > \max\{R_\alpha, c_1 |\alpha|^2\}$ for each $\alpha \neq 0$, so that the non-radiative solution $u^\alpha$ satisfies $\|u^\alpha\|_{\dot{H}^1(\{x: |x|>r_\alpha\})} = A$, then any radial solution to (CP1) with a maximal lifespan $(-T_-,T_+)$ and
 \[
  \limsup_{t\rightarrow T_+} \|(u(\cdot,t), u_t(\cdot,t))\|_{\dot{H}^1 \times L^2(\Rm^3)} < A, 
 \]
 must satisfy $T_+ = +\infty$ and scatter in the whole space in the positive time direction. 
\end{theorem}

\begin{remark}
 In the case of focusing wave equation $\partial_t^2 u - \Delta u = |u|^4 u$, if $A$ is smaller than the $\dot{H}^1(\Rm^3)$ norm $A_0$ of the ground states, then we may find a constant $c_1$, so that the identity $\|u^\alpha\|_{\dot{H}^1(\{x: |x|>c_1 |\alpha|^2\})} = A$ holds for all $\alpha\neq 0$ by the rescaling invariance. It immediately follows Theorem \ref{main global} that if a solution $u$ to the focusing wave equation satisfies 
\[
 \sup_{t\in [0,T_+)} \left\|(u(\cdot,t), u_t(\cdot,t))\right\|_{\dot{H}^1\times L^2(\Rm^3)} < A_0, 
\]
then $u$ must be defined for all time $t>0$ and scatters in the positive time direction. 
\end{remark}

\noindent Finally we give an application 
\begin{corollary}[Scattering in the defocusing case] \label{main corollary defocusing}
We assume that the nonlinear term $F(x,u)$ satisfies (A1)-(A4). In addition, the nonlinear term is defocusing, i.e. 
\[
 u F(x,u) \leq 0, \qquad \forall (x,u) \in \Rm^3 \times \Rm. 
\]
Then any radial solution to (CP1) with initial data $(u_0,u_1) \in \dot{H}^1\times L^2(\Rm^3)$ must be globally defined for all $t$ and scatter in both time directions. 
\end{corollary}

\subsection{Ideas and main tools} 

The major tool of this work is the radiation fields of linear free waves (see Theorem \ref{radiation} below). We compare the radiation profiles of the initial data of two asymptotically equivalent solutions and show that their difference must decay at a certain rate. The decay of non-radiative linear free waves plays an essential role in this argument, as it does in the spacial case of non-radiative solutions (see \cite{nonradialCE}). The general case turns out to be more difficult, because the decay of general solutions to the nonlinear wave equation is much slower than that of the non-radiative solutions. In general, the $Y$ norm of non-radiative solutions in the exterior region $\Omega_R$ decays polynomially $\|\chi_R u\|_{Y(\Rm)} \lesssim R^{-\kappa}$ as $R$ is sufficiently large. In contrast, the decay of a general exterior solution to (CP1) can only be characterized by an $l^2$ boundedness
 \[
  \sum_{k = N}^\infty  \|\tilde{\chi}_{2^k} u\|_{Y(\Rm)}^2 < +\infty, \qquad \forall N\gg 1. 
 \]
Here $\tilde{\chi}_{2^k}$ is the characteristic function of the channel region $\{(x,t): |t| +2^k < |x|<|t|+2^{k+1}\}$.  Another key observation is the following decomposition in a suitable exterior region
\[
 v^\alpha = v^0 + u^\alpha + \hbox{small error terms}.
\]
Here $\{v^\alpha\}$ is the primary one-parameter family associated with a small linear free wave $v_L$ and $\{u^\alpha\}$ the one-parameter family of non-radiative solutions. This decomposition plays an essential role in our discussion of global behaviours of solutions to (CP1). 

\subsection{Structure of this work}
This work is organized as follows: In Section 2 we give a few preliminary results. We then introduce the decay estimate of radial free waves and non-linear solutions in Section 3. These two sections are both preparation work of our main theory. The main results of this work are proved in the last 4 sections. Section 4 and 5 discuss the one-parameter family of any radial linear free wave and its special case of the non-radiative solutions, respectively. Finally we discuss the global behaviours of solutions to (CP1) as an application in Section 6 and 7. 

\section{Preliminary Results}

\subsection{Strichartz estimates} 
We recall the Strichartz estimates: if $u$ solves the linear inhomogeneous wave equation 
\[
 \partial_t^2 u - \Delta u = F(x,t), \qquad (x,t) \in \Rm^3 \times \Rm,
\]
with initial data $(u_0,u_1)\in \dot{H}^1 \times L^2(\Rm^3)$, then the following inequality holds for any time interval $I$ containing zero
\begin{align}
 \sup_{t\in I}\|(u,u_t)\|_{\dot{H}^1\times L^2(\Rm^3)} + \|u\|_{L^5 L^{10} (I \times \Rm^3)} \leq C_1 \left[\|(u_0,u_1)\|_{\dot{H}^1 \times L^2(\Rm^3)} + \|F\|_{L^1 L^2(I \times \Rm^3)}\right]. \label{Strichartz estimates}
\end{align}
Here the constant $C_1$ does not depends on $u$ or the time interval $I$. More general case of Strichartz estimates and their proof can be found in Ginibre-Velo \cite{strichartz}. 

\subsection{Radiation fields of free waves}

The radiation fields have a history of more than 50 years. Please see, Friedlander \cite{radiation1, radiation2} for instance. Generally speaking, radiation fields discuss the asymptotic behaviour of linear free waves. The radiation field is our main study tool of his work. The following version of statement comes from Duyckaerts-Kenig-Merle \cite{dkm3}.

\begin{theorem}[Radiation field] \label{radiation}
Assume that $d\geq 3$ and let $u$ be a solution to the free wave equation $\partial_t^2 u - \Delta u = 0$ with initial data $(u_0,u_1) \in \dot{H}^1 \times L^2(\Rm^d)$. Then ($u_r$ is the derivative in the radial direction)
\[
 \lim_{t\rightarrow \pm \infty} \int_{\Rm^d} \left(|\nabla u(x,t)|^2 - |u_r(x,t)|^2 + \frac{|u(x,t)|^2}{|x|^2}\right) dx = 0
\]
 and there exist two functions $G_\pm \in L^2(\Rm \times \mathbb{S}^{d-1})$ so that
\begin{align*}
 \lim_{t\rightarrow \pm\infty} \int_0^\infty \int_{\mathbb{S}^{d-1}} \left|r^{\frac{d-1}{2}} \partial_t u(r\theta, t) - G_\pm (r\mp t, \theta)\right|^2 d\theta dr &= 0;\\
 \lim_{t\rightarrow \pm\infty} \int_0^\infty \int_{\mathbb{S}^{d-1}} \left|r^{\frac{d-1}{2}} \partial_r u(r\theta, t) \pm G_\pm (r\mp t, \theta)\right|^2 d\theta dr & = 0.
\end{align*}
In addition, the maps $(u_0,u_1) \rightarrow \sqrt{2} G_\pm$ are bijective isometries from $\dot{H}^1 \times L^2(\Rm^d)$ to $L^2 (\Rm \times \mathbb{S}^{d-1})$. 
\end{theorem}

\noindent We call $G_\pm$ the radiation profiles of this linear free wave $u$, or equivalently, of the corresponding initial data $(u_0,u_1)$. The map between radiation profiles $G_\pm$ is an isometry from $L^2(\Rm \times \mathbb{S}^{d-1})$ to itself. The formula between $G_\pm$ is relatively simpler in odd dimensions than even dimensions. In this work we only need to use the 3-dimensional case (please see \cite{newradiation, shenradiation} for all dimensions, for example)
\[
 G_+(s,\theta) = -G_-(-s, -\theta). 
\]
It is clear that the free wave is radial if and only if its radiation profiles are independent of the angle $\theta$. The formula of a free wave in term of its radiation profile can also be given explicitly, see \cite{shenradiation}, for example. In this work we focus on the 3D radial case: 
\[
 u(r,t) = \frac{1}{r} \int_{t-r}^{t+r} G_-(s) ds.  
\]
A basic calculation gives the initial data in term of the radiation profile 
\begin{align} \label{initial data by radiation profile} 
 &u_0(r,t) =  \frac{1}{r} \int_{-r}^{r} G_-(s) ds; & & u_1(r,t) = \frac{G_-(r)+G_-(-r)}{r}.&
\end{align}
Throughout his paper, when we mention radiation profile of free waves, we mean the radiation profile $G_-$ in the negative time direction, unless otherwise specified.

\subsection{Nonlinear radiation profiles}

\begin{lemma} [Radiation fields of inhomogeneous equation] \label{scatter profile of nonlinear solution}
 Let $u$ be a radial solution to the linear wave equation
 \[
  \left\{\begin{array}{ll} \partial_t^2 u - \Delta u = F(t,x); & (x,t)\in \Rm^3 \times \Rm; \\
  (u,u_t)|_{t=0} = (u_0,u_1). & \end{array} \right.
 \]
 If $F\in L^1 L^2 (\Rm \times \Rm^3)$ is a radial function, then there exists $G^\pm \in L^2(\Rm)$ so that 
 \begin{align*}
  \lim_{t\rightarrow +\infty} \int_0^\infty \left(\left|G^+(r-t) - r u_t (r, t)\right|^2 + \left|G^+(r-t) + r u_r (r, t)\right|^2\right) dr & = 0; \\
  \lim_{t\rightarrow -\infty} \int_0^\infty \left(\left|G^-(r+t) - r u_t(r,t)\right|^2 +  \left|G^-(r+t) - r u_r(r,t)\right|^2\right)dr & = 0.
 \end{align*}
 In addition, the following estimates hold for $G^\pm$ given above and the radiation profile $G_0^\pm$ of the initial data $(u_0,u_1)$ in both two time directions ($R \geq 0$)
 \begin{align*}
 2\sqrt{2\pi} \|G^- - G_0^-\|_{L^2([R,+\infty))} & \leq \|\chi_R F\|_{L^1 L^2((-\infty,0]\times \Rm^3)};& \\
  2\sqrt{2\pi}  \|G^+ - G_0^+\|_{L^2([R,+\infty))} & \leq \|\chi_R F\|_{L^1 L^2([0,+\infty)\times \Rm^3)}.&
 \end{align*}
\end{lemma}
\begin{proof}
Without loss of generality we may assume $(u_0,u_1) = 0$. Otherwise we split the solution $u$ into two parts: the linear propagation of the initial data and the contribution of inhomogeneous term. We first apply the Strichartz estimates and obtain 
\begin{align*}
&  \lim_{t_1, t_2 \rightarrow +\infty} \left\|\mathbf{S}_L(-t_1) \begin{pmatrix} u(\cdot,t_1)\\ u_t(\cdot,t_1)\end{pmatrix} - \mathbf{S}_L(-t_2) \begin{pmatrix} u(\cdot,t_2)\\ u_t(\cdot,t_2)\end{pmatrix}\right\|_{\dot{H}^1 \times L^2(\Rm^3)} \\
= & \lim_{t_1, t_2 \rightarrow +\infty} \left\|\mathbf{S}_L(t_2-t_1) \begin{pmatrix} u(\cdot,t_1)\\ u_t(\cdot,t_1)\end{pmatrix} - \begin{pmatrix} u(\cdot,t_2)\\ u_t(\cdot,t_2)\end{pmatrix}\right\|_{\dot{H}^1 \times L^2(\Rm^3)}\\
\lesssim_1 &\lim_{t_1, t_2 \rightarrow +\infty}  \|F\|_{L^1 L^2([t_1,t_2]\times \Rm^3)} = 0.
\end{align*}
By the completeness of the space $\dot{H}^1 \times L^2$, there exists $(u_0^+, u_1^+)\in \dot{H}^1 \times L^2$ so that 
\[
 \lim_{t\rightarrow +\infty} \left\|\mathbf{S}_L(-t) \begin{pmatrix} u(\cdot,t)\\ u_t(\cdot,t)\end{pmatrix} - \begin{pmatrix} u_0^+\\ u_1^+\end{pmatrix}\right\|_{\dot{H}^1 \times L^2(\Rm^3)} = 0.
\] 
Thus we have 
\begin{equation} \label{forward scattering}
\lim_{t\rightarrow +\infty} \left\| \begin{pmatrix} u(\cdot,t)\\ u_t(\cdot,t)\end{pmatrix} - \mathbf{S}_L(t) \begin{pmatrix} u_0^+\\ u_1^+\end{pmatrix}\right\|_{\dot{H}^1 \times L^2(\Rm^3)} = 0.
\end{equation}
Let $G^+$ be the radiation profile of the free wave $u_L^+ = \mathbf{S}_L (u_0^+, u_1^+)$ in the positive time direction. By the property of radiation fields we have 
\[
  \lim_{t\rightarrow +\infty} \int_0^\infty \left( \left|G^+(r-t) - r (\partial_t u_L^+) (r, t)\right|^2 + \left|G^+(r-t) + r (\partial_r u_L^+) (r, t)\right|^2 \right) dr  = 0.
\]
We may combine this with \eqref{forward scattering} to conclude 
\[
 \lim_{t\rightarrow +\infty} \int_0^\infty \left(\left|G^+(r-t) - r u_t (r, t)\right|^2 + \left|G^+(r-t) + r u_r (r, t)\right|^2\right) dr = 0.
\]
In addition, if $R>0$ is a constant, then the limit above implies that 
\begin{align*}
 \|G^+\|_{L^2([R,+\infty))} = \|G^+(r-t)\|_{L_r^2([t+R, +\infty))} & = \lim_{t\rightarrow +\infty} \|r u_t(r,t)\|_{L_r^2([t+R, +\infty))}\\
 & = \lim_{t\rightarrow +\infty} \|ru_r(r,t)\|_{L_r^2([t+R, +\infty))}.
\end{align*}
Finally we combine the Strichartz estimates with finite speed of propagation to conclude 
\begin{align*}
  2\sqrt{2\pi}\|G^+\|_{L^2([R,+\infty))} & \leq \lim_{t\rightarrow +\infty} \|\nabla_{t,x} u\|_{L^2(\{x: |x|>t+R\})} \\
  & \leq \lim_{t\rightarrow +\infty} \|\chi_R F\|_{L^1 L^2([0,t)\times \Rm^3)} = \|\chi_R F\|_{L^1 L^2([0,+\infty)\times \Rm^3)}. 
\end{align*}
The proof in the negative time direction is similar.
\end{proof}

\begin{corollary} \label{comparison of profile}
 Let radial solutions $u$ and $v$ solve the linear wave equation $\partial_t^2 u - \Delta u = F$ and $\partial_t^2 v - \Delta v = \tilde{F}$, respectively. Here $F, \tilde{F} \in L^1 L^2(\Rm \times \Rm^3)$. Their initial data come with radiation profiles $G_\pm, \tilde{G}_\pm \in L^2(\Rm)$, respectively. In addition, we assume 
 \[
  \lim_{t\rightarrow \pm \infty} \int_{|x|>|t|+R} |\nabla_{t,x} u - \nabla_{t,x} v|^2 dx = 0. 
 \]
 Then we have
 \[ 
  \|G_\pm -\tilde{G}_\pm \|_{L^2([R,+\infty))} \lesssim_1 \|\chi_R (F-\tilde{F})\|_{L^1 L^2(\Rm \times \Rm^3)}, \qquad \forall R\geq 0.
 \] 
\end{corollary}
\begin{proof}
 We consider the solution $w = u-v$ to the equation $\partial_t^2 w - \Delta w = F-\tilde{F}$ and apply Lemma \ref{scatter profile of nonlinear solution}. Our assumption of $u$ and $v$ implies that the radiation profile $G^+$ obtained satisfies $G^+(s) = 0$ for $s > R$. According to Lemma \ref{scatter profile of nonlinear solution} we also have 
 \[
  \|G_+ -\tilde{G}_+ \|_{L^2([R,+\infty))} \leq \|(G_+ -\tilde{G}_+)-G^+\|_{L^2([R,+\infty))} \lesssim_1 \|\chi_R (F-\tilde{F})\|_{L^1 L^2(\Rm \times \Rm^3)}.
 \]
 The negative time direction is similar.
 \end{proof}
 
 \begin{remark} \label{scatter profile of nonlinear exterior solution}
  Lemma \ref{scatter profile of nonlinear solution} also applies to exterior solutions of the linear wave equation 
 \[
  \partial_t^2 u - \Delta u = F(x,t), \qquad (x,t) \in \Omega_R;
 \]
 with initial data $(u_0,u_1)\in \dot{H}^1 \times L^2(\Rm^3)$. More precisely, we assume $u$ and $F$ are both defined in the exterior region $\Omega_R$ so that $\|\chi_R u\|_{Y(\Rm)}<+\infty$, $\|\chi_R F\|_{L^1 L^2(\Rm \times \Rm^3)} < +\infty$ and 
 \[
  u = \mathbf{S}_L (u_0,u_1) + \int_0^t \frac{\sin (t-t')\sqrt{-\Delta}}{\sqrt{-\Delta}} [\chi_R(\cdot,t') F(\cdot,t')] dt', \qquad (x,t) \in \Omega_R.
 \]
 Our conclusions are also concerning the behaviour of solutions in the exterior region only. For example, in the positive time direction we have $G^+\in L^2([R,+\infty))$ with
 \[
   \lim_{t\rightarrow +\infty} \int_{R+t}^\infty  \left(\left|G^+(r-t) - r u_t (r, t)\right|^2 + \left|G^+(r-t) + r u_r (r, t)\right|^2\right) dr = 0;
 \]
 and 
 \[
   2\sqrt{2\pi}  \|G^+ - G_0^+\|_{L^2([R,+\infty))} \leq \|\chi_R F\|_{L^1 L^2([0,+\infty)\times \Rm^3)}.
 \]
 In order to prove this exterior version it suffices to consider the global extension
 \[
  u = \mathbf{S}_L (u_0,u_1) + \int_0^t \frac{\sin (t-t')\sqrt{-\Delta}}{\sqrt{-\Delta}} [\chi_R(\cdot,t') F(\cdot,t')] dt', \qquad (x,t) \in \Rm^3 \times \Rm;
 \]
which solves the equation $\partial_t^2 u - \Delta u = \chi_R F$ in the whole space-time $\Rm^3 \times \Rm$, and apply Lemma \ref{scatter profile of nonlinear solution} on it.  Similarly Corollary \ref{comparison of profile} also applies to exterior solutions.
 \end{remark}

\begin{lemma} \label{equivalent to exterior scattering}
 Let $u$ be a radial exterior solution to (CP1) in the exterior region $\Omega_R$ with a maximal lifespan $(-T_-, T_+)$. Then the following three statements are equivalent to each other: 
 \begin{itemize} 
  \item[(i)] $u$ scatters exteriorly in the positive time direction, i.e. $T_+ = \infty$ and there exists a finite-energy free wave $v_L^+$ so that 
\[
 \lim_{t\rightarrow +\infty} \|\vec{u}(x,t)-\vec{v}_L^+(x,t)\|_{\mathcal{H}_{t+R}} = 0,
\]
 \item[(ii)] $\displaystyle \sup_{t\in [0,T_+)} \|(u(\cdot,t), u_t(\cdot,t))\|_{\dot{H}^1 \times L^2(\{x: |x|>|t|+R\})} < +\infty.$
 \item[(iii)] $\|\chi_R u\|_{Y([0,T_+))} < +\infty$. 
\end{itemize} 
\end{lemma}
\begin{proof}
 We consider the global extension of $u$ and still use the same notation $u$. It suffices to show that 
\[
 (iii) \Rightarrow (i) \Rightarrow (ii) \Rightarrow (iii).
\]
We recall that $(u,u_t) \in C([0,T_+); \dot{H}^1 \times L^2(\Rm^3))$, therefore (i) implies (ii). Next we assume (ii) and prove (iii). We first apply a center cut-off technique. By small data scattering and finite speed propagation, there always exists a radius $R_1\gg R$, so that $\|\chi_{R_1} u\|_{Y([0,T_+))} < +\infty$. Thus it suffice to show that $\|\chi_{R,R_1} u\|_{Y([0,T_+))} < +\infty$. It immediately follows the following point-wise upper bound (here we use the radial assumption and (ii))
 \[
  |u(r,t)| \lesssim r^{-1/2} \|u(\cdot,t)\|_{\dot{H}^1(\{x: |x|>R+|t|\})} \lesssim r^{-1/2}, \qquad \forall r > R+|t|, \; t\in [0,T_+).
 \]
 and the assumption $\|\chi_R u\|_{Y([0,T])} < +\infty$ for any $T \in (0, T_+)$. Finally we show $(iii) \Rightarrow (i)$. If $\|\chi_R u\|_{Y([0,T_+))} < +\infty$, then the finite-time blow-up criterion gives that $T_+ = +\infty$. we then observe the fact 
 \[
  \|\chi_R F(x,t,u)\|_{L^1 L^2([0, +\infty))} \lesssim_\gamma \|\chi_R u\|_{Y([0,+\infty))}^5 < +\infty. 
 \]
 and follow the same argument as given in the proof of Lemma \ref{scatter profile of nonlinear solution} to obtain the scattering.
\end{proof}

\subsection{Small data theory}

\begin{lemma}\label{nonlinear comparison by linear}
 Let $u$, $v$ be exterior solutions to (CP1) in the region $\Omega_R$ with initial data $(u_0,u_1)$ and $(v_0,v_1)$, respectively, so that $\|\chi_R u\|_{Y(I)}$, $\|\chi_R v\|_{Y(I)}$ are both sufficiently small, then we have 
 \[
  \|\chi_R (u-v)\|_{Y(I)} \leq 2\|\chi_R \mathbf{S}_L (u_0-v_0, u_1-v_1)\|_{Y(I)}. 
 \]
 Here $I$ is an arbitrary time interval containing zero. 
\end{lemma} 
\begin{proof}
The inequality follows the Strichartz estimates (let $w_L = \mathbf{S}_L (u_0-v_0, u_1-v_1)$)
\begin{align*}
 \|\chi_R (u-v)\|_{Y(I)} &\leq \|\chi_R w_L\|_{Y(I)} + C_1 \|\chi_R (F(x,t,u) - F(x,t,v))\|_{L^1 L^2(I \times \Rm^3)}\\
 & \leq \|\chi_R w_L\|_{Y(I)} + C_1\gamma (\|\chi_R u\|_{Y(I)}^4 + \|\chi_R v\|_{Y(I)}^4) \|\chi_R (u-v)\|_{Y(I)}
\end{align*}
Here $C_1$ is an absolute constant. The conclusion clearly holds as long as $\|\chi_R u\|_{Y(I)}, \|\chi_R v\|_{Y(I)} \leq \delta(\gamma) \doteq (4C_1\gamma)^{-1/4}$. 
\end{proof}

\begin{lemma}\label{nonlinear comparison by linear 2}
 Assume that $u$ is an exterior solution to (CP1) in the region $\Omega_R$ with initial data $(u_0,u_1)$. Let $(v_0,v_1) \in \dot{H}^1 \times L^2$ be initial data and $w_L = \mathbf{S}_L (u_0-v_0, u_1-v_1)$. If both $\|\chi_R u\|_{Y(\Rm)}$ and $\|\chi_R w_L\|_{Y(\Rm)}$ are both sufficiently small, then the exterior solution $v$ to (CP1) in the exterior region $\Omega_R$ with initial data $(v_0,v_1)$ must be well-defined for all $t$, so that 
  \[
  \|\chi_R (u-v)\|_{Y(\Rm)} \leq 2\|\chi_R w_L\|_{Y(\Rm)}. 
 \]
\end{lemma}
\begin{proof}
 We assume that $\|\chi_R u\|_{Y(\Rm)}, \|\chi_R w_L\|_{Y(\Rm)} \leq \delta(\gamma) /4$. Here $\delta(\gamma)$ is the constant in Lemma \ref{nonlinear comparison by linear}. It suffices to show that $v$ is globally defined for all $t$ with $\|\chi_R v\|_{Y(\Rm)} \leq \delta(\gamma)$ and apply Lemma \ref{nonlinear comparison by linear}. If $v$ were not globally defined or $\|\chi_R v\|_{Y(\Rm)} > \delta(\gamma)$, then we could find a bounded closed time interval $I$ contained in its maximal lifespan, so that $0\in I$ and $\|\chi_R v\|_{Y(I)} = \delta(\gamma)$. We then apply Lemma \ref{nonlinear comparison by linear} and obtain $\|\chi_R (u-v)\|_{Y(I)} \leq \delta(\gamma)/2$, thus $\|\chi_R v\|_{Y(I)} \leq 3\delta(\gamma)/4$. This is a contradiction thus finishes the proof. 
\end{proof}

\begin{lemma} \label{asymptoticly equivalent to square small}
 There exists a constant $\delta = \delta(\gamma)>0$, so that if $v_L = \mathbf{S}_L (v_0,v_1)$ is a finite-energy radial free wave with $\|\chi_R v_L\|_{Y(\Rm)} < \delta$, then
 there exists at least one exterior solution $u$ to (CP1) with initial data $(u_0,u_1) \in \dot{H}^1 \times L^2(\Rm^3)$ satisfying:
 \begin{itemize}
  \item The solution $u$ is defined in the region $\Omega_R$;
  \item The solution $u$ is $R$-weakly asymptoticly equivalent to $v_L$; 
  \item The initial data of $u$ and $v_L$ satisfies $\|(u_0,u_1)-(v_0,v_1)\|_{\dot{H}^1 \times L^2(\Rm^3)} \lesssim_\gamma \|\chi_R v_L\|_{Y(\Rm)}^5 $.
 \end{itemize}
\end{lemma}
\begin{proof}
We consider a complete distance space of radiation profiles
 \[
  X = \left\{G \in L^2(\Rm): \|G\|_{L^2(\Rm)} \leq 2\|\chi_R v_L\|_{Y(\Rm)}; G(s) = 0, \forall s\in [-R,R]\right\}, 
 \]
 whose distance is the $L^2$ norm of the difference. Given any $G \in X$, we let $(u_0,u_1) \in \dot{H}^1 \times L^2$ be the radial initial data so that the radiation profile of initial data $(u_0-v_0, u_1-v_1)$ is exactly $G$. Therefore the linear free wave $u_L = \mathbf{S}_L (u_0,u_1)$ satisfies
 \[
  \|\chi_R u_L\|_{Y(\Rm)} \lesssim_1 \|\chi_R v_L\|_{Y(\Rm)} + \|G\|_{L^2} \lesssim_1 \|\chi_R v_L\|_{Y(\Rm)} \ll 1.
 \]
 Thus the exterior solution to (CP1) with initial data $(u_0,u_1)$ is defined for all $t$ so that 
 \[
  \|\chi_R u\|_{Y(\Rm)} \leq 2\|\chi_R u_L\|_{Y(\Rm)} \lesssim_1 \|\chi_R v_L\|_{Y(\Rm)}. 
 \]
Next we consider the global extension of $u$ (also called $u$) and let $G^\pm$ be the nonlinear radiation profile of it, as given in Lemma \ref{scatter profile of nonlinear solution}. We have 
\begin{align}
 \|G^+ + G(-s)  + G_0(-s)\|_{L^2([R,+\infty))} & \lesssim_1 \|\chi_R F(x,t,u)\|_{L^1 L^2(\Rm \times \Rm^3)}\lesssim_\gamma \|\chi_R u\|_{Y(\Rm)}^5 \nonumber \\& \lesssim_\gamma \|\chi_R v_L\|_{Y(\Rm)}^5. \label{upper bound G 1 int}
\end{align}
Here $G_0$ is the radiation profile of $v_L$. Similarly we have 
\begin{equation} \label{upper bound G 2 int}
 \|G^- - G(s)  - G_0\|_{L^2([R,+\infty))} \lesssim_\gamma \|\chi_R v_L\|_{Y(\Rm)}^5. 
\end{equation} 
Therefore when $\|\chi_R v_L\|_{Y(\Rm)} < \delta(\gamma)$ is sufficiently small, the map 
\[
 (\mathbf{T} G) (s) = \left\{\begin{array}{ll}G(s) -G^-(s) + G_0 (s), & s > R; \\ 0, & s\in [-R,R]; \\ G(s) + G^+(-s) + G_0 (s), & s< -R.\end{array} \right. 
\]
is a map from $X$ to itself. Next We show that $\mathbf{T}$ is actually a contraction map. Given another $\tilde{G} \in X$, we define $\tilde{u}_L$, $\tilde{u}_L$, $\tilde{G}^\pm$ accordingly. We have 
\[
 \|\chi_R (u - \tilde{u})\|_{Y(\Rm)} \leq 2\|\chi_R (u_L - \tilde{u}_L)\|_{Y(\Rm)} \lesssim_1 d(G, \tilde{G}). 
\]
Therefore we have 
\[
 \|\chi_R (F(x,t,u) - F(x,t,\tilde{u}))\|_{L^1 L^2(\Rm \times \Rm^3)} \lesssim_\gamma \|\chi_R v_L\|_{Y(\Rm)}^4 d(G, \tilde{G}). 
\]
Next we observe that $u-\tilde{u}$ solves $\partial_t^2 w - \Delta w = \chi_R (F(x,t,u) - F(x,t,\tilde{u}))$, apply Lemma \ref{scatter profile of nonlinear solution}, utilize the estimate above and obtain
\begin{align*}
 \|G^+ -\tilde{G}^+ + G(-s) - \tilde{G}(-s)\|_{L^2([R,+\infty))} & \lesssim_\gamma \|\chi_R v_L\|_{Y(\Rm)}^4 d(G, \tilde{G});\\
 \|G^- -\tilde{G}^- - G + \tilde{G}\|_{L^2([R,+\infty))} & \lesssim_\gamma \|\chi_R v_L\|_{Y(\Rm)}^4 d(G, \tilde{G}).
\end{align*}
In summary we have 
\[
 d(\mathbf{T} G, \mathbf{T}\tilde{G}) \lesssim_\gamma  \|\chi_R v_L\|_{Y(\Rm)}^4 d(G, \tilde{G}). 
\]
This verifies that $\mathbf{T}$ is a contraction map when $\|\chi_R v_L\|_{Y(\Rm)} < \delta(\gamma)$ is sufficiently small. Let $G$ be the unique fixed point and $u$ be the corresponding solution. By $\mathbf{T} G = G$ we have 
\begin{align} \label{indeed asymptotically equivalent}
 &G^- (s) = G_0(s),& &G^+(s) = - G_0(-s),& &s>R.& 
\end{align}
Thus the solution $u$ is indeed $R$-weakly asymptoticly equivalent solution of $v_L$. Finally we combine the two identities \eqref{indeed asymptotically equivalent} with \eqref{upper bound G 1 int}, \eqref{upper bound G 2 int} and obtain
\[
 \|(u_0-v_0,u_1-v_1)\|_{\dot{H}^1 \times L^2} \lesssim_1 \|G\|_{L^2} \lesssim_\gamma \|\chi_R v_L\|_{Y(\Rm)}^5. 
\] 
\end{proof}

\begin{corollary} \label{nonlinear asymptotic with small norm profiles}
 There exist a constant $\delta = \delta (\gamma) > 0$, so that if $v_L = \mathbf{S}_L (v_0,v_1)$ is a radial free wave satisfying $\|(v_0,v_1)\|_{\dot{H}^1 \times L^2(\Rm^3)} < \delta$, then there exists a global-in-time solution $u$ to (CP1) satisfying
 \begin{itemize}
  \item The solution $u$ scatters in both two time directions;
  \item The solution $u$ is asymptoticly equivalent to $v_L$;
  \item The initial data $(u_0,u_1)$ satisfy $\|(u_0, u_1)\|_{\dot{H}^1 \times L^2(\Rm^3)} \leq 2\|(v_0,v_1)\|_{\dot{H}^1 \times L^2(\Rm^3)}$. 
 \end{itemize}
\end{corollary}

\subsection{Extension of scattering exterior solutions} 
 
\begin{lemma} 
Assume $R>0$. Let $u$ be an exterior solution to (CP1) defined in the region $\Omega_R$ with initial data $(u_0,u_1)\in \dot{H}^1 \times L^2(\Rm^3)$ so that $\|\chi_R u\|_{Y(\Rm)} < +\infty$. Then there exists a radius $r\in (0,R)$, so that we may extend the domain of $u$ to $\Omega_r$, so that $u$ is still an exterior solution to (CP1) in $\Omega_r$ with the same initial data and $\|\chi_r u\|_{Y(\Rm)} < +\infty$. 
\end{lemma}
\begin{proof}
 We slightly abuse the notation to let $u$ represent the global extension of the exterior solution $u$. If suffices to show that when $r<R$ is sufficiently close to $R$, the unique exterior solution $v$ to (CP1) with initial data $(u_0,u_1)$ in the region $\Omega_r$ can be defined for all time $t$ and satisfies $\|\chi_r v\|_{Y(\Rm)} < +\infty$. We may apply Strichartz estimates and obtain 
\[
  \|u\|_{Y(\Rm)} \lesssim_1 \|(u_0,u_1)\|_{\dot{H}^1 \times L^2(\Rm^3)} + \|\chi_R u\|_{Y(\Rm)}^5 < +\infty.
\]
We choose $r \in (0,R)$ so that $\|\chi_{r,R} u\|_{Y(\Rm)} < \delta(\gamma) \ll 1$ is sufficiently small. Let $v$ be the exterior solution to (CP1) with initial data $(u_0,v_0)$ in the region $\Omega_r$ (as well as its global extension) and $J = [-T_1,T_2]$ be a time interval contained in its maximal lifespan. By finite speed of propagation we have
\[
 u(x,t) = v(x,t), \qquad t\in J, \; |x|>|t|+R.
\]
We observe that $w = v-u$ solves the equation $\partial_t^2 w - \Delta w = \chi_r F(x,t,v) - \chi_R F(x,t,u)$ with zero initial data. The Strichartz estimates then give
\begin{align*}
 \|\chi_{r,R}(v-u)\|_{Y(J)} &\lesssim_1 \|\chi_r F(x,t,v) - \chi_R F(x,t,u)\|_{L^1 L^2 (J\times \Rm^3)} = \|\chi_{r,R} F(x,t,v)\|_{L^1 L^2 (J\times \Rm^3)} \\
 & \lesssim_\gamma \|\chi_{r,R} v\|_{Y(J)}^5. 
\end{align*} 
Therefore we have 
\[
 \|\chi_{r,R} v\|_{Y(J)} \leq \|\chi_{r,R} u\|_{Y(\Rm)} + C_\gamma \|\chi_{r,R} v\|_{Y(J)}^5.
\]
Since the norm $\|\chi_{r,R} u\|_{Y(\Rm)}$ is a small number, a continuity argument shows that the inequality $\|\chi_{r,R} v\|_{Y(J)} \leq 2 \|\chi_{r,R} u\|_{Y(\Rm)}$ holds for any time interval $J$ contained in the maximal lifespan of $v$. Thus $\|\chi_r v\|_{Y(J)}$ is uniformly bounded for any time interval $J$ contained in the maximal lifespan of $v$. This implies that $v$ is defined for all time $t$ and $\|\chi_r v\|_{Y(\Rm)} < +\infty$. 
\end{proof}

\subsection{Technical Lemma}

\begin{lemma} \label{sequence lemma}
 Let $\{a_k\}_{k\geq 0}$ and $\{y_k\}_{k\geq 0}$ be both nonnegative sequences so that $a_k \in l_k^1$ and the inequality 
 \[
  y_k \leq \sum_{j=0}^{k-1} \gamma^{k-j} a_j y_j
 \]
 holds for sufficiently large $k$ and a constant $\gamma \in (0,1)$. Then we have $y_k \lesssim \gamma^k$ for sufficiently large $k$. 
\end{lemma}
\begin{proof}
 Without loss of generality we may assume that the inequality holds for all $k \geq 1$. We use the new notation $z_k = \gamma^{-k} y_k$ and rewrite the inequality in the form of 
 \[
  z_k \leq \sum_{j=0}^{k-1} a_j z_j, \qquad k\geq 1. 
 \]
 This gives
 \[
  z_k \leq a_{k-1} z_{k-1} + \sum_{j=0}^{k-2} a_j z_j \leq (a_{k-1} +1) \sum_{j=0}^{k-2} a_j z_j \leq \cdots \leq \left(\prod_{j=1}^{k-1} (1+a_j)\right)a_0 z_0.
 \]
 Finally we conclude that $z_k$ is uniformly bounded thus finish the proof, by observing the basic fact that 
 \[
  a_k \in l^1 \quad \Rightarrow \quad \prod_{j=1}^{+\infty} (1+a_j) < +\infty. 
 \]
\end{proof}

\section{Decay of radial solutions}

In this section we discuss the asymptotic decay of any global solutions to (CP1). We first recall a few notations. Let $\chi_R$ and $\tilde{\chi}_{R}$ be characteristic functions of the regions $\{(x,t): |x|>|t|+R\}$ and $\{(x,t): |t|+R<|x|<|t|+2R\}$, respectively. 

\subsection{Decay of free waves}

\begin{lemma} \label{channel by localized radiation}
 Let $u_L$ be a radial free wave with a compactly-supported radiation profile $G \in L^2([-2R,-R]\cup [R,2R])$. Then we have 
\[
 \|\tilde{\chi}_{R_1} u_L\|_{Y(\Rm)} \lesssim_1 (R_1/R)^{1/10} \|G\|_{L^2}, \qquad \forall R_1\in (0,R].
\]
\end{lemma}
\begin{proof}
 We recall the explicit formula of the free wave $u_L$ in term of its radiation profile 
 \[
  u_L (r,t) = \frac{1}{r} \int_{t-r}^{t+r} G(s) ds. 
 \]
 This immediately gives a pointwise estimate 
 \[
  |u_L(r,t)| \lesssim_1 \frac{R^{1/2}}{r} \|G\|_{L^2}. 
 \]
 A straight-forward calculation then gives the upper bound
 \[
  \|u_L (\cdot,t)\|_{L^{10} (\{x: |t|+R_1 < |x| < |t|+2R_1\})} \lesssim_1 \frac{R_1^{1/10} R^{1/2}}{(t+R_1)^{4/5}} \|G\|_{L^2}. 
 \]
 This immediately gives the upper bound of $\|\tilde{\chi}_{R_1} u_L\|_{Y(\Rm)}$ if $R_1 \geq R/4$. On the other hand, if $R_1 < R/4$, the explicit formula above implies that $\tilde{\chi}_{R_1} u_L = 0$ for all $x\in \Rm^3$ when $|t| < R/4$. Thus in this case we utilize the $L^{10}$ estimate above and obtain 
 \[
  \|\tilde{\chi}_{R_1} u_L\|_{Y(\Rm)} \lesssim_1 \left(\int_{|t|>R/4} \frac{R_1^{1/2} R^{5/2}}{(t+R_1)^{4}} dt \right)^{1/5} \|G\|_{L^2} \lesssim_1 (R_1/R)^{1/10} \|G\|_{L^2}. 
 \]
\end{proof}

\begin{remark} \label{channel by small support data}
Let $u_L$ be a radial free wave with a compactly-supported radiation profile $G \in L^2([-2R,-R]\cup [R,2R])$. In addition, we assume 
$\left|\{s\in \Rm: G(s) \neq 0\}\right| < R_0$, then we have 
\[
 \|\tilde{\chi}_{R_1} u_L\|_{Y(\Rm)} \lesssim_1 (R_0/R)^{1/2} (R_1/R)^{1/10} \|G\|_{L^2}, \qquad \forall R_1\in (0,R].
\]
The proof is similar to that of Lemma \ref{channel by localized radiation}, but we have a better point-wise estimate $|u_L(r,t)| \lesssim_1 R_0^{1/2} r^{-1} \|G\|_{L^2}$. 
\end{remark}

\begin{lemma} \label{wide channel by localized radiation}
 Let $u_L$ be a radial free wave with a radiation profile $G \in L^2(\Rm)$ so that 
 \[
  \left|\{s\in \Rm: G(s) \neq 0\}\right| \leq R.
 \]
 Then we have
 \[
  \|\chi_{R_1} u_L\|_{Y(\Rm)} \lesssim_1 (R/R_1)^{1/2} \|G\|_{L^2}. 
 \]
\end{lemma}
\begin{proof}
 By the same argument as in the proof of the previous lemma, we have the pointwise estimate $|u_L(r,t)| \lesssim_1 R^{1/2} r^{-1} \|G\|_{L^2}$. A straight-forward calculation then finishes the proof. 
\end{proof}

\begin{proposition} \label{l2 linear}
 Let $u_L$ be a radial free wave with initial data $(u_0,u_1)$. Then we have 
 \[
  \sum_{k=-\infty}^\infty \|\tilde{\chi}_{2^k} u_L\|_{Y(\Rm)}^2 \lesssim_1 \|(u_0,u_1)\|_{\dot{H}^1 \times L^2(\Rm^3)}^2
 \]
\end{proposition}
\begin{proof}
 Let the radiation profile of $u_L$ be $G \in L^2(\Rm)$. We make a dyadic decomposition 
\begin{align*}
 &G = \sum_{k=-\infty}^\infty G_k;& &\hbox{Supp} \, G_k \subset [-2^{k+1}, -2^k]\cup [2^k, 2^{k+1}].&
\end{align*}
Then we way apply Lemma \ref{channel by localized radiation}, Lemma \ref{wide channel by localized radiation} and obtain 
\[
 \|\tilde{\chi}_{2^k} u_L\|_{Y(\Rm)} \lesssim_1 \sum_{j=-\infty}^\infty 2^{-\frac{|j-k|}{10}} \|G_j\|_{L^2}.
\]
The right hand side is the convolution of the $l^2$ sequence $\|G_j\|_{L^2}$ with an $l^1$ sequence $c_j \doteq 2^{-|j|/10}$. Therefore we may apply Young's inequality and obtain
\[
 \sum_{k=-\infty}^\infty \|\tilde{\chi}_{2^k} u_L\|_{Y(\Rm)}^2 \lesssim_1 \|c_j\|_{l^1}^2 \sum_{k=-\infty}^\infty \|G_k\|_{L^2}^2\lesssim_1 \|(u_0,u_1)\|_{\dot{H}^1\times L^2}^2. 
\]
\end{proof}

\subsection{Decay of nonlinear solution}

\begin{lemma}\label{channel by nonlinear preliminary} 
 Assume that $R\geq 4R_1 > 0$ and $t_0\geq 0$ . Let $u_1 \in L^2(\Rm^3)$ be radial initial data compactly supported in the region $\{x: t_0+R<|x|<t_0+2R\}$. Then the free wave $u_L$ with initial data $(0,u_1)$ satisfies 
 \[
  \|u_L\|_{L^5(\Rm; L^{10}(\{x: |t|+t_0+R_1 < |x|< |t|+t_0+2R_1\}))} \lesssim_1 (R_1/R)^{1/10} \|u_1\|_{L^2(\Rm^3)}. 
 \]
\end{lemma}
\begin{proof}
 We first recall the explicit formula 
 \[
  u_L(r,t) = \frac{1}{2r} \int_{r-t}^{r+t} s u_1(s) ds.
 \]
 This gives a pointwise estimate 
 \begin{align*}
  |u_L(r,t)| & \lesssim_1 \frac{1}{r} \int_{t_0+R<|s|<t_0+2R} |s u_1(s)| ds\\
  & \lesssim_1 \frac{R^{1/2}}{r} \left(\int_{t_0+R<|s|<t_0+2R} |s u_1(s)|^2 ds\right)^{1/2} \lesssim_1 \frac{R^{1/2}}{r}\|u_1\|_{L^2(\Rm^3)}
 \end{align*}
 In addition, finite speed of propagation gives 
 \[
  u_L(x,t) = 0, \qquad |x|<|t|+t_0+2R_1, |t|<R/4.
 \]
 Thus a straight-forward calculation as given in the proof of Lemma \ref{channel by localized radiation} finishes the proof. 
\end{proof}

\noindent Next we recall the solution of $\partial_t^2 u - \Delta u = F$ with zero initial data is given by the Duhamel's formula 
\[
 u(\cdot,t) = \int_0^t \mathbf{S}_L (t-\tau) (0, F(\cdot,\tau)) d\tau. 
\]
and obtain
\begin{corollary} \label{channel by nonlinear}
 Assume $R\geq R_1>0$. Let $u$ be the solution to the linear wave equation $\partial_t^2 u - \Delta u = F$ with zero initial data. Here $F(x,t) \in L^1 L^2 (\Rm \times \Rm^3)$ is a radial function supported in $\tilde{\Omega}_R$. Then we have 
 \[
  \|\tilde{\chi}_{R_1} u\|_{Y(\Rm)} \lesssim (R_1/R)^{1/10} \|F\|_{L^1 L^2(\Rm \times \Rm^3)}. 
 \]
\end{corollary}

\begin{lemma} \label{l2 small nonlinear} 
Let $(u_0,u_1)\in \dot{H}^1 \times L^2$ be radial initial data. There exists a constant $\delta = \delta (\gamma)$, so that if the corresponding free wave $u_L$ satisfies
\[ 
 \sum_{k=-\infty}^\infty \|\tilde{\chi}_{2^k} u_L\|_{Y(\Rm)}^2 < \delta^2. 
\]
Then the exterior solution $u$ to (CP1) in the region $\Omega_0$ is defined for all time $t$. In addition, we have 
\[
 \sum_{k=-\infty}^\infty  \|\tilde{\chi}_{2^k} u\|_{Y(\Rm)}^2 \leq 4 \sum_{k=-\infty}^\infty \|\tilde{\chi}_{2^k} u_L\|_{Y(\Rm)}^2. 
\]
\end{lemma}
\begin{proof}
 First of all, we have 
\[
 \|\chi_0 u_L\| \leq \left(\sum_{k=-\infty}^\infty \|\tilde{\chi}_{2^k} u_L\|_{Y(\Rm)}^5\right)^{1/5}\leq \delta \ll 1.
\]
Thus $u$ can be obtained by a standard fixed-point argument. More precisely we define the sequence 
 \begin{align*}
  &u^{0} = 0;& &u^{n+1} = u_L + \int_0^t \mathbf{S}_L (t-\tau) (0, \chi_0 F(\cdot,\tau, u^n(\cdot, \tau))) d\tau.& 
 \end{align*}
 The solution $u$ is the limit of $u^n$ in the space $L_t^5 L^{10} (\{x: |x|>|t|\})$. For convenience we introduce the following notations 
 \begin{align*}
  &b_k = \|\tilde{\chi}_{2^k} u_L\|_{Y(\Rm)};& &a_{n,k} = \|\tilde{\chi}_{2^k} u^n\|_{Y(\Rm)}.&
 \end{align*}
We then combine finite speed of propagation and Lemma \ref{channel by nonlinear} to deduce an recursion formula
 \[
  a_{n+1,k} \leq b_k + c \sum_{j\geq k} 2^{-(j-k)/10}\|\tilde{\chi}_{2^j} F(x,t,u^n)\|_{L^1 L^2(\Rm \times \Rm^3)} \leq b_k + c\gamma \sum_{j \geq k} 2^{-(j-k)/10} a_{n,j}^5. 
 \]
 Thus we have 
 \[
  \|a_{n+1, k}\|_{l_k^2} \leq \|b_k\|_{l_k^2} + c\gamma \|a_{n,k}\|_{l_k^\infty}^4 \left\| \sum_{j \geq k} 2^{-(j-k)/10} a_{n,j}\right\|_{l_k^2} \leq \|b_k\|_{l_k^2} + c_1 \gamma \|a_{n,k}\|_{l_k^2}^5. 
 \]
 Here we apply Young's inequality on the convolution of an $l^1$ function and $\{a_{n,k}\}_k$. The letters $c$ and $c_1$ above represent absolute constants. Therefore if $\|b_k\|_{l_k^2} < \delta(\gamma)$ is sufficiently small, an induction immediately gives that $\|a_{n,k}\|_{l_k^2} \leq 2 \|b_k\|_{l_k^2}$ for all nonnegative integers $n$. We then pass to limit and obtain 
 \[
  \sum_{k=-\infty}^\infty  \|\tilde{\chi}_{2^k} u\|_{Y(\Rm)}^2 \leq 4\|b_k\|_{l_k^2}^2. 
 \]
\end{proof}

\noindent Combining this lemma with Proposition \ref{l2 linear}, we obtain 

\begin{proposition} \label{L2 small nonlinear} 
 Let $u$ be a solution to (CP1) with radial initial data $(u_0,u_1)$. There exists a constant $\delta = \delta (\gamma)$, so that if $\|(u_0,u_1)\|_{\dot{H}^1\times L^2(\Rm^3)} < \delta$, then $u$ is a global solution satisfying the inequality 
 \[
  \sum_{k=-\infty}^\infty  \|\tilde{\chi}_{2^k} u\|_{Y(\Rm)}^2 \lesssim \|(u_0,u_1)\|_{\dot{H}^1 \times L^2(\Rm^3)}^2. 
 \]
\end{proposition}

\begin{corollary} \label{l2 exterior nonlinear}
 Let $u$ be an exterior solution of (CP1) defined in the region $\Omega_R$. Then we have 
 \[
  \sum_{k = N}^\infty  \|\tilde{\chi}_{2^k} u\|_{Y(\Rm)}^2 < +\infty, \qquad \forall N\gg 1. 
 \]
\end{corollary}
\begin{proof}
 Let $\phi: \Rm^3 \rightarrow [0,1]$ be a radial, smooth center cut-off function so that $\phi(x)=1$ if $|x|\geq 1$ and $\phi(x) = 0$ if $|x|\leq 1/2$. Given initial data $(u_0,u_1)$ of the solution $u$, there exists a sufficiently large radius $R_1 > R$, so that the cut-off version of initial data $\phi(x/R_1) (u_0,u_1)$ satisfy
 \[
  \|\phi(x/R_1) (u_0,u_1)\|_{\dot{H}^1\times L^2(\Rm^3)} < \delta.
 \] 
 We than apply Proposition \ref{L2 small nonlinear} and conclude that the solution $v$ to (CP1) with initial data $\phi(x/R)(u_0,u_1)$ satisfies the inequality 
 \[
  \sum_{k=-\infty}^\infty  \|\tilde{\chi}_{2^k} v\|_{Y(\Rm)}^2 < +\infty.
 \]
 Finally finite speed of propagation gives that fact $u = v$ in the region $\Omega_{2^N}$ as long as $2^N > R_1$. This immediately finishes the proof. 
\end{proof}

\section{One-parameter Family}

In this section, we prove Theorem \ref{main 1}. Given any radial finite-energy free wave $v_L$, there exists at least one exterior solution $u^0$ so that $u^0$ is weakly asymptotically equivalent to $v_L$, thanks to Lemma \ref{asymptoticly equivalent to square small}. We then figure out all other asymptotically equivalent solutions $u^\alpha$. The proof consists of the following four steps. We will do each step in an individual subsection. 
\begin{itemize}
 \item[(1)] We first show that any weakly equivalent solution of $v_L$ comes with a characteristic number (with respect to the chosen $u^0$);
 \item[(2)] We show that given any real number $\alpha$, we may find a weakly equivalent solution $u^\alpha$ with characteristic number $\alpha$;
 \item[(3)] We extend the domain of $u^\alpha$ to the maximum, so that it satisfies either of the two conditions given in Theorem \ref{main 1}; 
 \item[(4)] We prove that any weakly equivalent exterior solution of $v_L$ is covered by the one-parameter family.
\end{itemize}
The fifth and sixth subsections of this section discuss additional properties of asymptotic solutions and will be used later. 

\subsection{The characteristic number}

\begin{lemma} \label{L1 lemma}
 Let $u$ and $v$ be two radial exterior solutions to (CP1) in the region $\Omega_R$ so that they are $R$-weakly asymptotically equivalent to each other. Let $G$ and $\tilde{G}$ be the radiation profiles of their initial data, respectively. Then we have 
 \begin{itemize}
  \item[(i)] The inequality $\|\chi_r (u-v)\|_{Y(\Rm)} \lesssim r^{-1/2}$ holds for all sufficiently large $r$;
  \item[(ii)]  We have $G-\tilde{G} \in L^1(\Rm)$. 
  \item[(iii)] Let $\alpha = \int_{\Rm} [G(s)-\tilde{G}(s)] ds$. The following limit holds
  \[
 \lim_{r\rightarrow +\infty} r^{1/2} \sup_{t\in \Rm}\|\vec{u}(x,t) - \vec{v}(x,t)-(\alpha|x|^{-1},0)\|_{\mathcal{H}_{|t|+r}} = 0.
\]
 \end{itemize}
\end{lemma}
\begin{proof} 
 First of all, Corollary \ref{l2 exterior nonlinear} (as well as its proof) guarantees that there exists a number $N \gg 1$, so that 
 \begin{align*}
  &\sum_{k=N}^\infty \|\tilde{\chi}_{2^k} u\|_{Y(\Rm)}^2 + \|\tilde{\chi}_{2^k} v\|_{Y(\Rm)}^2 < +\infty;& &\|\chi_{2^N} u\|_{Y(\Rm)} + \|\chi_{2^N} v\|_{Y(\Rm)} \ll 1.&
 \end{align*}
 Without loss of generally we assume that $R=1$ and $N=0$, the proof of other cases are exactly the same. For convenience we introduce the notations 
 \begin{align*}
  &b_k = \|\tilde{\chi}_{2^k} u\|_{Y(\Rm)} + \|\tilde{\chi}_{2^k} v\|_{Y(\Rm)};& &y_k = \|\chi_{2^k} (u-v)\|_{Y(\Rm)}.&
 \end{align*}
 We also slightly abuse the notation, so that $u$ and $v$ also represent their corresponding global extensions. We apply Corollary \ref{comparison of profile} and obtain ($k \geq 0$)
\begin{align}
 \|G-\tilde{G}\|_{L^2(\{s: |s|>2^k\})} & \lesssim_1 \|\chi_{2^k} (F(x,t,u)-F(x,t,v))\|_{L^1 L^2 (\Rm \times \Rm^3)} \nonumber\\
 & \lesssim_1 \sum_{j \geq k} \|\tilde{\chi}_{2^j} (F(x,t,u)-F(x,t,v))\|_{L^1 L^2 (\Rm \times \Rm^3)} \nonumber \\
 & \lesssim_\gamma \sum_{j\geq k} b_j^4 y_j. \label{difference of G}
\end{align}
This gives the following estimates concerning the initial data $(u_0,u_1)$ and $(v_0,v_1)$ by Lemma \ref{wide channel by localized radiation} and Strichartz estimates
\begin{align*}
 \|\chi_{2^k} \mathbf{S}_L (u_0-v_0,u_1-v_1)\|_{Y(\Rm)} & \lesssim_1 2^{-\frac{k}{2}} \|G-\tilde{G}\|_{L^2([-1,1])} + \sum_{j=0}^{k-1} 2^{-\frac{k-j}{2}} \|G-\tilde{G}\|_{L^2(\{s: 2^j<|s|<2^{j+1}\})}\\
 & \qquad + \|G-\tilde{G}\|_{L^2(\{s: |s|>2^{k}\})}\\
 & \lesssim_\gamma 2^{-\frac{k}{2}} \|G-\tilde{G}\|_{L^2([-1,1])} + \sum_{j=0}^{k-1} \sum_{i\geq j} 2^{-\frac{k-j}{2}} b_i^4 y_i + \sum_{j\geq k} b_j^4 y_j\\
 & \lesssim_\gamma 2^{-\frac{k}{2}} \|G-\tilde{G}\|_{L^2([-1,1])} + \sum_{j=0}^{k-1} 2^{-\frac{k-j}{2}} b_j^4 y_j + \sum_{j\geq k} b_j^4 y_j.
\end{align*}
We then apply Lemma \ref{nonlinear comparison by linear} to conclude 
\[
 y_k \lesssim_\gamma 2^{-\frac{k}{2}} \|G-\tilde{G}\|_{L^2([-1,1])} + \sum_{j=0}^{k-1} 2^{-\frac{k-j}{2}} b_j^4 y_j + \sum_{j\geq k} b_j^4 y_j. 
\]
We observe that $y_k$ is a decreasing sequence of $k$ and $\sum_j b_j^4 < +\infty$. Thus the final term in the right hand side above can be absorbed by the left hand side if $k$ is sufficiently large. This gives
\[
 y_k \lesssim_\gamma 2^{-\frac{k}{2}} \|G-\tilde{G}\|_{L^2([-1,1])} + \sum_{j=0}^{k-1} 2^{-\frac{k-j}{2}} b_j^4 y_j, \qquad \forall k\gg 1. 
\]
Lemma \ref{sequence lemma} then guarantees that $y_k \lesssim 2^{-k/2}$ for sufficiently large $k$. It immediately follows that $y_k \lesssim 2^{-k/2}$ for all $k\geq 0$. This proves part (i). Next we recall \eqref{difference of G} and obtain
\begin{equation} \label{L2 tail estimate}
 \|G-\tilde{G}\|_{L^2\{s: |s|>2^k\}} \lesssim \sum_{j\geq k} b_j^4 y_j \lesssim \sum_{j\geq k}  2^{-j/2} b_j^4, \qquad \forall k \geq 0. 
\end{equation}
This immediately gives 
\[
 \int_{2^k < |s| < 2^{k+1}} |G(s)-\tilde{G}(s)| ds \lesssim \sum_{j\geq k}  2^{-\frac{j-k}{2}} b_j^4, \qquad k\geq 0.
\]
Thus for any integer $N \geq 0$, we have
\begin{equation} \label{L1 tail estimate}
 \int_{|s|>2^N} |G(s)-\tilde{G}(s)| ds \lesssim \sum_{k\geq N} \sum_{j\geq k}  2^{-\frac{j-k}{2}} b_j^4 \lesssim \sum_{j\geq N} b_j^4 < +\infty. 
\end{equation} 
This shows that $G-\tilde{G} \in L^1(\Rm)$. In order to prove (iii), we use the explicit formula of initial data in term of the radiation profiles and write
\begin{align*}
 &u_0(r) - v_0(r) = \frac{1}{r} \int_{-r}^r [G(s) - \tilde{G}(s)] ds;& &u_1(r) - v_1 (r) = \frac{G(s)-\tilde{G}(r)}{r} - \frac{G(-r)-\tilde{G}(-r)}{r}.&
\end{align*}
Therefore if $r\gg 1$ is a large radius, then
\begin{align*}
 \left|\frac{d}{dr}\left(u_0(r) - v_0(r) - \frac{\alpha}{r} \right)\right| & = \left|\frac{d}{dr} \left( \frac{1}{r} \int_{|s|>r}  [G(s)-\tilde{G}(s)] ds \right)\right| \\
 & \leq \frac{1}{r^2} \int_{|s|>r} |G(s)-\tilde{G}(s)| ds + \frac{|G(r)-\tilde{G}(r)|}{r} + \frac{|G(-r)-\tilde{G}(-r)|}{r}\\
 & \lesssim \frac{1}{r^2} \sum_{2^j > r/2} b_j^4 + \frac{|G(r)-\tilde{G}(r)|}{r} + \frac{|G(-r)-\tilde{G}(-r)|}{r}. 
\end{align*}
Thus 
\[
 \|u_0 - v_0 - \alpha |x|^{-1} \|_{\dot{H}^1(\{x: |x|>r\})} \lesssim r^{-1/2} \sum_{2^j >r/2} b_j^4 = o(r^{-1/2}).
\]
Similar we have 
\[
 \|u_1-v_1\|_{L^2(\{x:|x|>r\})} \lesssim r^{-1/2} \sum_{2^j >r/2} b_j^4 = o(r^{-1/2}). 
\]
Here we use \eqref{L2 tail estimate} and \eqref{L1 tail estimate}. Finally we combine Strichartz estimates, finite speed of propagation to obtain
\begin{align*}
 \sup_{t\in \Rm}\|\vec{u}(x,t) - \vec{v}(x,t)-(\alpha|x|^{-1},0)\|_{\mathcal{H}_{|t|+r}} &\lesssim_1 \|(u_0,u_1)-(v_0,v_1)-(\alpha|x|^{-1},0)\|_{\mathcal{H}_r} \\
 & \qquad + \|\chi_r(F(x,t,u)-F(x,t,v))\|_{L^1 L^2(\Rm \times \Rm^3)}
\end{align*}
We finish the proof by observing that the right hand side is still $o(r^{-1/2})$ by our estimate of initial data and the estimates $\|\chi_r(u-v)\|_{Y(\Rm)} \lesssim r^{-1/2}$.
\end{proof}

\paragraph{The characteristic number} Let $u$, $v$ be two exterior solutions to (CP1) as in Lemma \ref{L1 lemma}. We define the value of the integral 
\[
 \int_{\Rm} [G(s) -\tilde{G}(s)] ds
\]
to be the characteristic number of $u$ with respect to $v$. This definition applies to more general solutions $u$ and $v$, even if $u$, $v$ are defined in different exterior regions, or they are not asymptotically equivalent in the whole overlap region of their domains, as long as they are $R$-weakly asymptotically equivalent for a sufficiently large number $R>0$. We simply apply Lemma \ref{L1 lemma} on the restriction of them in $\Omega_R$. Please note that the choice of large radius $R$ does not affect the characteristic number. This is because although $G$ and $\tilde{G}$ do depend on the choice of initial data in the interior region $\{x: |x|\leq R\}$, which are irrelevant to the exterior solution defined in $\Omega_R$, the value of the integral above, however, can be unique determined by the behaviour of solutions near infinity, as shown in part (iii) of Lemma \ref{L1 lemma}. Here we use the fact 
\[
 \||x|^{-1}\|_{\dot{H}^1(\{x: |x|>r\})} \simeq r^{-1/2}. 
\]
We recall that we have fixed an exterior solution $u^0$ to (CP1), so that $u^0$ is weakly asymptotically equivalent to $v_L$. If another radial exterior solution $u$ to (CP1) is weakly asymptoticly equivalent to $v_L$, then $u$ and $u^0$ are also weakly asymptotically equivalent to each other, we call the characteristic number of $u$ with respect to $u^0$ to be the characteristic number of $u$ for simplicity. 

\subsection{Existence}

In this subsection we show that given any radial exterior solution $v$ to (CP1) and a real number $\alpha$, there exists an exterior solution $u$ to (CP1), which is weakly asymptoticly equivalent to $v$ and whose characteristic number (with respect to $v$) is exactly $\alpha$. Please note that the domain of $u$ may be smaller than that of $v$. If we choose $v=u^0$, then this gives asymptotically equivalent solutions with every possible characteristic number. 

\begin{lemma} \label{existence of solutions far} 
 Given a radial exterior solution $v$ to (CP1) defined for all time $t\in \Rm$ and a real number $\alpha$, there exists a radial exterior solution $u$ to (CP1), so that $u$ is weakly asymptotically equivalent to $v$ and the characteristic number of $u$ with respect to $v$ is exactly $\alpha$. 
\end{lemma}
\begin{proof}
We first choose a radius $R>0$ so that $v$ is defined in $\Omega_R$ and $\|\chi_R v\|_{Y(\Rm)} < +\infty$. For convenience we let $G_0$ be the radiation profile of initial data $(v_0,v_1)$ and introduce the notations 
\[
 b_k = \|\tilde{\chi}_{2^k} v\|_{Y(\Rm)}, \qquad \forall k \geq \log_2 R,\, k\in \mathbb{Z}.
\]
According to Corollary \ref{l2 exterior nonlinear}, we have $\{b_k\} \in l^2$. We consider the complete distance space 
\[
 X = \left\{G \in L^2(\{s: |s|>2^N\}): \|G\|_{L^2(\{s: |s|>2^k\})} \leq c\left(|\alpha|^5 2^{-5k/2} + |\alpha|\sum_{j\geq k} b_j^4 2^{-j/2}\right) , \forall k \geq N\right\},
\]
whose distance is defined by 
\[
 d(G_1,G_2) = \sup_{k\geq N} \left(|\alpha|^5 2^{-5k/2} + |\alpha|\sum_{j\geq k} b_j^4 2^{-j/2}\right)^{-1} \|G_1-G_2\|_{L^2(\{s: |s|>2^k\})}. 
\]
Here $c = c(\gamma) \gg 1$ is a number to be determined later and $N = N(\gamma, c, \{b_k\})\geq \log_2 R$  is a positive integer satisfying
\begin{align*}
 &c\sum_{j \geq N} b_j^4 \leq 1/100;& &2^N \geq 100 c^{1/2} |\alpha|^2.&
\end{align*}
The choice of $N$ guarantees that we have 
\[
 \|\chi_{2^N} v\|_{Y(\Rm)} \lesssim \left(\sum_{j\geq N} b_j^5\right)^{1/5} \ll 1. 
\]
We also define a map $\mathbf{T}$ from $X$ to itself. Given $G\in X$, we first extend its domain to all $s\in \Rm$ so that 
\[
 G(s) = \frac{\alpha}{2^{N+1}} - \frac{1}{2^{N+1}}\int_{|s'|>2^N} G(s') ds', \qquad |s|\leq 2^N;.
\]
and let $(u_0,u_1)\in \dot{H}^1 \times L^2(\Rm^3)$ be initial data so that the corresponding radiation profile is exactly $G_0 + G$. Let $u$ be the corresponding solution to the non-linear wave equation 
\[ 
 \partial_t^2 u - \Delta u = \chi_{2^N} F(t,x,u)
\]
In the argument below we slightly abuse the notation to let $v$ be the global extension of the restriction of $v$ in the region $\Omega_{2^N}$. This will not affect the values of $v$ in the exterior region $\Omega_{2^N}$. We use the notation $w_L = \mathbf{S}_L (u_0-v_0,u_1-v_1)$ and consider the upper bound of $\|\chi_{2^k} w_L\|_{Y(\Rm)}$. We observe that the radiation profile of $w_L$ is exactly $G$ thus we have 
\begin{align*}
 w_L (r,t) = \frac{1}{r} \int_{t-r}^{t+r} G(s) ds.
\end{align*}
If $r > |t| + 2^k$ for an integer $k\geq N$, then we have 
\begin{align*}
 |w_L(r,t)| \leq \frac{|\alpha|}{r} + \frac{1}{r} \int_{|s|>2^k} |G(s)| ds \lesssim_1 \left(|\alpha| + c|\alpha| \sum_{j \geq k} b_j^4 + c|\alpha|^5 2^{-2k}\right) \frac{1}{r}. 
\end{align*}
A straight-forward calculation shows (here we recall the choice of $N$)
\begin{equation} \label{upper bound of wL tem}
 \|\chi_{2^k} w_L\|_{Y(\Rm)} \lesssim_1 |\alpha| 2^{-k/2}+ c|\alpha| 2^{-k/2} \sum_{j \geq k} b_j^4 + c|\alpha|^5 2^{-5k/2} \lesssim_1 |\alpha| 2^{-k/2} \ll 1. 
\end{equation}
We then apply Lemma \ref{nonlinear comparison by linear 2} and obtain 
\begin{equation}
 \|\chi_{2^k} (u-v)\|_{Y(\Rm)} \lesssim_1 |\alpha| 2^{-k/2}. \label{upper bound on Y u}
\end{equation}
Let $G^\pm$, $G_0^\pm$ be the radiation profiles of the nonlinear solution $u$ and $v$, respectively. We apply Lemma \ref{scatter profile of nonlinear solution} on $u-v$ and obtain ($k \geq N$)
\begin{align}
 \|G^+ - G_0^+ + G(-s)\|_{L^2([2^k,+\infty))} & \lesssim_1 \|\chi_{2^k} F(x,t,u) - \chi_{2^k} F(x,t,v)\|_{L^1 L^2}\nonumber \\
 & \lesssim_1 \sum_{j \geq k} \|\tilde{\chi}_{2^j} F(x,t,u) - \tilde{\chi}_{2^j} F(x,t,v)\|_{L^1 L^2}\nonumber \\
 & \lesssim_\gamma \sum_{j \geq k} (\|\tilde{\chi}_{2^j} v\|_{Y(\Rm)}^4 + \|\tilde{\chi}_{2^j} (u-v)\|_{Y(\Rm)}^4) \|\tilde{\chi}_{2^j} (u-v)\|_{Y(\Rm)}\nonumber \\
 & \lesssim_\gamma |\alpha|\sum_{j \geq k} b_j^4 2^{-j/2} + |\alpha|^5 2^{-5k/2}. \label{T estimate 1}
\end{align}
Similarly in the negative time direction we have ($k \geq N$)
\begin{equation}
\|G^- - G_0^- - G(s)\|_{L^2([2^k,+\infty))}  \lesssim_\gamma |\alpha|\sum_{j \geq k} b_j^4 2^{-j/2} + |\alpha|^5 2^{-5k/2}. \label{T estimate 2}
\end{equation} 
We then define a map $\mathbf{T}$ on $X$
\[
 (\mathbf{T} G) (s) = \left\{\begin{array}{ll} G(s) - G^-(s) + G_0^-(s), & s>2^N; \\ G(s) + G^+(-s) - G_0^+(-s), & s<-2^N. \end{array} \right.
\]
If $c = c(\gamma) \gg 1$ is sufficiently large, our assumption on $N$ and the estimates \eqref{T estimate 1}, \eqref{T estimate 2} guarantee that $\mathbf{T} X \subseteq X$. Next we claim that $\mathbf{T}$ is a contraction map from $X$ to itself. Thus there exists a fixed-point $G \in X$. This implies that the nonlinear radiation profiles $G^\pm$ of the corresponding solution $u$ satisfy $G^\pm (s) = G_0^\pm (s)$ for all $s>2^N$, thus verifies that $u$ is weakly asymptoticly equivalent to $v$. The way we define the values of $G(s)$ for $|s|<2^N$ guarantees that the characteristic number of $u$ is exactly $\alpha$. The remaining work is to show that $\mathbf{T}$ is a contraction map. Given $G, \tilde{G} \in X$, we use the notations $(u_0,u_1), (\tilde{u}_0, \tilde{u}_1)$, $u, \tilde{u}$ and $G^\pm, \tilde{G}^\pm$ for the corresponding initial data, solutions and nonlinear radiation profiles. Again we first consider the linear free wave $\bar{w}_L = \mathbf{S}_L (u_0-\tilde{u}_0, u_1-\tilde{u}_1)$, whose radiation profile is $G-\tilde{G}$. if $k\geq N$ and $r>|t|+2^k$, then we have
\[
 |\bar{w}_L(r,t)| = \left|\frac{1}{r} \int_{t-r}^{t+r} (G(s) -\tilde{G}(s)) ds \right| \leq \frac{1}{r} \int_{|s|>2^k} |G(s) - \tilde{G}(s)| ds \lesssim_1 \frac{d}{r}\left(|\alpha| \sum_{j \geq k} b_j^4 + |\alpha|^5 2^{-2k}\right).
\]
Here $d = d(G, \tilde{G})$. A basic calculation shows that
\[
 \|\chi_{2^k} \bar{w}_L\|_{Y(\Rm)} \lesssim_1 \left(|\alpha| 2^{-k/2} \sum_{j \geq k} b_j^4 + |\alpha|^5 2^{-5k/2}\right) d. 
\]
In addition, we may combine $\|\tilde{\chi}_{2^k}v\|_{Y(\Rm)} = b_k$ with \eqref{upper bound on Y u} and write 
\[ 
 \|\tilde{\chi}_{2^k} u\|_{Y(\Rm)}, \|\tilde{\chi}_{2^k} \tilde{u}\|_{Y(\Rm)} \lesssim_1 |\alpha| 2^{-k/2} + b_k \ll 1. 
\]
We then apply Lemma \ref{nonlinear comparison by linear} and obtain 
\[
 \|\chi_{2^k} (u-\tilde{u})\|_{Y(\Rm)} \lesssim_1 \left(|\alpha| 2^{-k/2} \sum_{j \geq k} b_j^4 + |\alpha|^5 2^{-5k/2}\right) d, \qquad k\geq N. 
\]
We apply Lemma \ref{scatter profile of nonlinear solution} on $u-\tilde{u}$ and obtain 
\begin{align*}
 \|G^+ - \tilde{G}^+ + G(-s) - \tilde{G}(-s)\|_{L^2([2^k,+\infty))} & \lesssim_1 \|\chi_{2^k} F(x,t,u) - \chi_{2^k} F(x,t,\tilde{u})\|_{L^1 L^2} \\
 & \lesssim_1 \sum_{j \geq k} \|\tilde{\chi}_{2^j} F(x,t,u) - \tilde{\chi}_{2^j} F(x,t,\tilde{u})\|_{L^1 L^2} \\
 & \lesssim_\gamma \sum_{j \geq k} (\|\tilde{\chi}_{2^j} u\|_{Y(\Rm)}^4 + \|\tilde{\chi}_{2^j} \tilde{u}\|_{Y(\Rm)}^4) \|\tilde{\chi}_{2^j} (u-\tilde{u})\|_{Y(\Rm)}. 
\end{align*}
We insert the upper bounds of $\|\tilde{\chi}_{2^k} u\|_{Y(\Rm)}$, $\|\tilde{\chi}_{2^k} \tilde{u}\|_{Y(\Rm)}$, $\|\chi_{2^k} (u-\tilde{u})\|_{Y(\Rm)}$ given above, and obtain 
\begin{align*}
 \|G^+ & - \tilde{G}^+ + G(-s) - \tilde{G}(-s)\|_{L^2([2^k,+\infty))}\\
 & \lesssim_\gamma d \sum_{j\geq k} \left(\left(|\alpha|^4 2^{-2j} + b_j^4 \right) |\alpha| 2^{-j/2} \left(|\alpha|^4 2^{-2j} + \sum_{i \geq j} b_i^4\right)\right)\\
 & \lesssim_\gamma d \left(|\alpha|^4 2^{-2k} + \sum_{i\geq k} b_i^4\right)\left(|\alpha|\sum_{j \geq k} 2^{-j/2} b_j^4 + |\alpha|^5 2^{-5k/2}\right).
\end{align*} 
A similar argument in the negative time direction gives 
\[
 \|G^- - \tilde{G}^- - G+ \tilde{G}\|_{L^2([2^k,+\infty))} \lesssim_\gamma d \left(|\alpha|^4 2^{-2k} + \sum_{i\geq k} b_i^4\right)\left(|\alpha|\sum_{j \geq k} 2^{-j/2} b_j^4 + |\alpha|^5 2^{-5k/2}\right).
\]
Finally we recall the definition of $\mathbf{T}$ and conclude 
\begin{align*}
 \|\mathbf{T} G - \mathbf{T} \tilde{G}\|_{L^2(\{s: |s|>2^k\})} & \leq \|G - G^- - \tilde{G} +\tilde{G}^-\|_{L^2([2^k, +\infty))} \\
 & \qquad +\|G+G^+(-s) - \tilde{G} -\tilde{G}^+(-s)\|_{L^2((-\infty, -2^k])}\\
 & \lesssim_\gamma d \left(|\alpha|^4 2^{-2k} + \sum_{i\geq k} b_i^4\right)\left(|\alpha|\sum_{j \geq k} 2^{-j/2} b_j^4 + |\alpha|^5 2^{-5k/2}\right). 
\end{align*}
Thus 
\begin{align*}
 d(\mathbf{T} G, \mathbf{T} \tilde{G}) & =  \sup_{k\geq N} \left(|\alpha|\sum_{j \geq k} 2^{-j/2} b_j^4 + |\alpha|^5 2^{-5k/2}\right)^{-1}  \|\mathbf{T} G - \mathbf{T} \tilde{G}\|_{L^2(\{s: |s|>2^k\})} \\
 & \lesssim_\gamma \sup_{k\geq N}  \left(|\alpha|^4 2^{-2k} + \sum_{i\geq k} b_i^4\right) d\\
 & \lesssim_\gamma \left(|\alpha|^4 2^{-2N} + \sum_{i\geq N} b_i^4\right) d. 
\end{align*}
Our assumption on $N$ then guarantees that the number in the inequality above satisfies
\[
 |\alpha|^4 2^{-2N} + \sum_{i\geq N} b_i^4 < \frac{1}{99} c^{-1}.
\]
As a result, $\mathbf{T}$ is a contraction map as long as $c=c(\gamma) \gg 1$ is sufficiently large. This finishes the proof. 
\end{proof}

\paragraph{More details} Part (iii) of Lemma \ref{L1 lemma} implies that the solutions $u$ and $v$ in the previous lemma satisty
 \[
   \|u(\cdot,0)-v(\cdot,0)\|_{\dot{H}^1(\{x: |x|>R\})} \simeq_1 |\alpha| R^{-1/2}, \qquad \forall R\gg 1. 
 \]
 We are interested in how large the radius $R$ need to be in this estimate, if $v$ is a sufficiently small solution. 

\begin{corollary} \label{lower bound far outside}
 If an exterior solution $v$ to (CP1) is defined in $\Omega_0$ so that ($c(\gamma)$ is the constant $c$ in the proof of Lemma \ref{existence of solutions far})
  \[
  \sum_{k=-\infty}^\infty \|\tilde{\chi}_{2^k} v\|_{Y(\Rm)}^4 < c(\gamma)/100, 
 \]
then the weakly asymptotically equivalent solution $u$ given by Lemma \ref{existence of solutions far} satisfies 
 \[
  \|u(\cdot,0)-v(\cdot,0)\|_{\dot{H}^1(\{x: |x|>R\})} \simeq_1 |\alpha| R^{-1/2}, \qquad \forall R \gtrsim_\gamma |\alpha|^2. 
 \]
\end{corollary}
\begin{proof}
This result follows a careful review of Lemma \ref{existence of solutions far}'s proof. We will use the same notation as in the proof of Lemma \ref{existence of solutions far}. Our assumption on $v$ guarantees that we may choose $2^N \simeq_\gamma |\alpha|^2$. Let $G \in X$ be the fixed-point of the contraction map $\mathbf{T}$. We have 
 \[
  u(r,0)-v(r,0) = \frac{1}{r} \int_{-r}^r G(s) ds. 
 \]
 Therefore 
 \begin{align*}
  u_r(r,0) - v_r(r,0) & = -\frac{1}{r^2} \int_{-r}^r G(s) ds + \frac{G(r) + G(-r)}{r}\\
  & = -\frac{\alpha}{r^2} + \frac{1}{r^2} \int_{|s|>r} G(s) ds +  \frac{G(r) + G(-r)}{r}. 
 \end{align*}
 A careful calculation shows that our definition of the space $X$ gives
 \begin{align*}
  &\int_{|s|>2^N} |G(s)| ds \leq \frac{|\alpha|}{5};& &\left\|\frac{G(r)+G(-r)}{r}\right\|_{L^2([R, +\infty); r^2 dr)} \leq \frac{|\alpha| R^{-1/2}}{10}, \quad R\geq 2^N.&
 \end{align*} 
 Therefore we have for $R \geq 2^N \simeq_\gamma |\alpha|^2$,
 \begin{align*}
  \|u_r(\cdot,0)-v_r(\cdot,0)\|_{L^2([R,+\infty); r^2 dr)} & \geq \left\|\frac{4|\alpha|}{5r^2}\right\|_{L^2([R,+\infty); r^2 dr)} - \left\|\frac{G(r)+G(-r)}{r}\right\|_{L^2([R, +\infty); r^2 dr)}\\
  & \geq \frac{7}{10} |\alpha| R^{-1/2}. 
 \end{align*}
 Similarly we have for $R \geq 2^N \simeq_\gamma |\alpha|^2$,
 \begin{align*}
  \|u_r(\cdot,0)-v_r(\cdot,0)\|_{L^2([R,+\infty); r^2 dr)} & \leq  \left\|\frac{6|\alpha|}{5r^2}\right\|_{L^2([R,+\infty); r^2 dr)} + \left\|\frac{G(r)+G(-r)}{r}\right\|_{L^2([R, +\infty); r^2 dr)}\\
  & \leq \frac{13}{10} |\alpha| R^{-1/2}. 
 \end{align*}
A combination of lower and upper bounds above finishes our proof.
\end{proof}

\subsection{The maximal domain} 

We recall that given any real number $\alpha$, there exists a radial exterior solution $u^\alpha$, so that $u^\alpha$ is weakly asymptotically equivalent to $v_L$ with characteristic number $\alpha$. In this subsection we extend the domain of $u^\alpha$ to the maximum, so that it is still weakly asymptotically equivalent to $v_L$ in the larger domain $\Omega_{R_\alpha}$, at least in a weak sense. This completes the construction of one-parameter family. The main result of this subsection is 

\begin{proposition} \label{extension behaviour}
 Let $v_L$ be a finite-energy radial free wave. If $u$ is an exterior solution to (CP1) defined in $\Omega_R$ so that $u$ is $R$-weakly asymptotically equivalent to $v_L$, then we may extend the domain of $u$ to a maximum so that exactly one of the following holds:
 \begin{itemize}
  \item $u$ becomes an exterior solution defined in $\Omega_0$ and is asymptotically equivalent to $v_L$;
  \item $u$ can be defined in $\Omega_{R_0}$ with $\|\chi_{R_0} u\|_{Y(\Rm)} = +\infty$; In addition, given any $r>R_0$, $u$ is an exterior solution in $\Omega_r$ and is $r$-weakly asymptotically equivalent to $v_L$.
 \end{itemize}
 In the first case we call $u$ a scattering asymptotically equivalent solution; while in the second case we call it a blow-up asymptotically equivalent solution. 
\end{proposition}

\noindent Before we prove this proposition, we first do some preparation work. 

\begin{lemma} \label{coincidence inside}
 Let $u$ and $v$ be two exterior solutions to (CP1) in $\Omega_R$ so that they are both $R$-weakly asymptotically equivalent to the same linear free wave $v_L$. If the identity $u(x,t) = v(x,t)$ holds in an exterior region $\Omega_{R_0}$ for some $R_0>R$, then the same identity also holds in the whole exterior region $\Omega_R$.  
\end{lemma}

\begin{proof}
First of all, Lemma \ref{equivalent to exterior scattering} implies that $\|\chi_R u\|_{Y(\Rm)}, \|\chi_R v\|_{Y(\Rm)} < +\infty$. Let $\mathcal{E}$ be the set of all real number $r\in [R, R_0]$ so that the identity $u(x,t) = v(x,t)$ holds in the exterior region $\Omega_{r}$. Clearly $R_0 \in \mathcal{E}$ thus $\mathcal{E}$ is nonempty. Let $R_1 = \inf \mathcal{E} \geq R$. We claim that $R_1 \in \mathcal{E}$. If this were false, then $R_1$ would be the limit of a decreasing sequence $r_k \in \mathcal{E}$. By definition of $\mathcal{E}$, we have 
 \[
  u(x,t) = v(x,t), \qquad \forall (x,t) \in \Omega_{r_k}. 
 \]
 Because $\Omega_{R_1} = \cup \Omega_{r_k}$, thus we have $R_1 \in \mathcal{E}$. Next we show $R_1=R$ by a proof of contradiction. If $R_1> R$, we would find a radius $R_2 \in (R, R_1)$ so that $R_2 \in \mathcal{E}$ in the argument below. This is a contradiction. Our argument starts by considering the radiation profiles $G$, $\tilde{G}$ of the initial data of $u$ and $v$, respectively. Because $u$ and $v$ coincides in $\Omega_{R_1}$, their initial data $(u_0,u_1)$ and $(v_0,v_1)$ also coincide in the exterior region $\{x: |x|>R_1\}$. We then recall the finite speed of wave propagation and the definition of radiation profiles (as well as the relationship between radiation profiles in two time directions) to conclude that $G(s)=\tilde{G}(s)$ for $|s|>R_1$. We then combine this identity with the explicit formula of initial data in term of radiation profiles to write
 \[
  u_0(r) - v_0(r) = \frac{1}{r} \int_{-R_1}^{R_1} [G(s) -\tilde{G}(s)] ds, \qquad r>R_1.
 \]
 Therefore we have
\[
 \int_\Rm [G(s) -\tilde{G}(s)] ds = 0.
\]
Next we consider $R_2 = R_1 -\delta > \max\{R,R_1/2\}$ and let $d = \|G-\tilde{G}\|_{L^2(\{s: R_2<|s|<R_1\})}$. We consider the free wave $w_L = \mathbf{S}_L(u_0-v_0,u_1-v_1)$. The finite speed of propagation and the coincidence of initial data in the exterior region gives 
\[
 w_L(r,t) = 0, \qquad \forall r>|t|+R_1. 
\]
In addition, if $|t|+R_2 < r < |t|+R_1$, then we use all the information of $G-\tilde{G}$ given above and obtain 
\[
 |w_L(r,t)| = \left|\frac{1}{r} \int_{t-r}^{t+r} [G(s) - \tilde{G}(s)] ds \right| \leq \frac{1}{r} \int_{R_2 < |s| <R_1} |G(s)-\tilde{G}(s)| ds \lesssim_1 \frac{\delta^{1/2} d}{r}.
\]
A straightforward calculation shows
\[
 \|\chi_{R_2} w_L\|_{Y(\Rm)} = \|\chi_{R_2, R_1}  w_L\|_{Y(\Rm)}\lesssim_1 (\delta/R_1)^{3/5} d. 
\]
We next apply Strichartz estimates, recall $u=v$ in $\Omega_{R_1}$ and obtain
\begin{align*}
 \|\chi_{R_2}(u-v)\|_{Y(\Rm)} & \leq \|\chi_{R_2} w_L\|_{Y(\Rm)} + C_1 \|\chi_{R_2} (F(x,t,u)-F(x,t,v))\|_{L^1 L^2(\Rm\times \Rm^3)} \\
 & \leq \|\chi_{R_2} w_L\|_{Y(\Rm)} + C_1\gamma (\|\chi_{R_2,R_1} u\|_{Y(\Rm)}^4 + \|\chi_{R_2,R_1} v\|_{Y(\Rm)}^4) \|\chi_{R_2}(u-v)\|_{Y(\Rm)}. 
\end{align*}
Here $C_1$ is an absolute constant. When $\delta$ is sufficiently small, we have 
\[
 \|\chi_{R_2,R_1} u\|_{Y(\Rm)}, \|\chi_{R_2,R_1} v\|_{Y(\Rm)} \ll 1. 
\]
Therefore in this case we have 
\begin{equation} \label{estimate difference small delta}
 \|\chi_{R_2}(u-v)\|_{Y(\Rm)} \leq 2\|\chi_{R_2} w_L\|_{Y(\Rm)} \lesssim_1 (\delta/R_1)^{3/5} d. 
\end{equation}
We then apply Corollary \ref{comparison of profile} and obtain 
\begin{align*}
 d & \lesssim_1 \|\chi_{R_2}(F(x,t,u)-F(x,t,v))\|_{L^1 L^2(\Rm \times \Rm^3)}\\
 & \lesssim_\gamma (\|\chi_{R_2,R_1} u\|_{Y(\Rm)}^4 + \|\chi_{R_2,R_1} v\|_{Y(\Rm)}^4) \|\chi_{R_2}(u-v)\|_{Y(\Rm)}\\
 & \lesssim_\gamma (\|\chi_{R_2,R_1} u\|_{Y(\Rm)}^4 + \|\chi_{R_2,R_1} v\|_{Y(\Rm)}^4) (\delta/R_1)^{3/5} d. 
\end{align*}
This immediately gives $d=0$ as long as $\delta$ is sufficiently small. We immediately obtain the identity $u(x,t) = v(x,t)$ for $(x,t) \in \Omega_{R_1-\delta}$ by \eqref{estimate difference small delta}.
\end{proof}

\paragraph{Domain extension} Given a finite-energy radial free wave $v_L$ and an $R$-weakly asymptotically equivalent solution $u$ to it, we let $\mathcal{D}$ be set of all real numbers $r\in [0,R]$ so that there exists a radial exterior solution $u^{(r)}$ to (CP1), which is $r$-weakly asymptotically equivalent to $v_L$ and whose restriction on $\Omega_R$ is $u$. We then let $R_0 = \inf \mathcal{D}$ and define a function in $\Omega_{R_0}$ by
\[
 u(x,t) = u^{(r)}(x,t), \qquad \forall (x,t) \in \Omega_{r}
\]
Given $r_1 \leq r_2$, since the solutions $u^{(r_1)}$ and $u^{(r_2)}$ are both $r_2$-weakly asymptotically equivalent to $v_L$ and coincide in $\Omega_R$, we must have $u^{(r_1)}(x,t) = u^{(r_2)}(x,t)$ for $(x,t)\in \Omega_{r_2}$ by Lemma \ref{coincidence inside}. Thus the function $u$ above is well-defined although the value of it at a given point may be defined for multiple times. Next we show that this extension $u$ satisfies the conditions in Proposition \ref{extension behaviour}. It is clear that $u$ is $r$-weakly asymptotically equivalent to $v_L$ for all $r > R_0$. In order to finish the proof we still need to show that if $\|\chi_{R_0}\|_{Y(\Rm)} < +\infty$, then we must have that $R_0 = 0$ and the extension $u$ is an asymptotically equivalent exterior solution to (CP1) of $v_L$. This immediately follows a combination of the two lemmata below. 
\begin{lemma}
 If the function $u$ defined above satisfies $\|\chi_{R_0} u\|_{Y(\Rm)} < +\infty$, then $u$ must be an exterior solution to (CP1) in $\Omega_{R_0}$ and $R_0$-weakly asymptotically equivalent to $v_L$. 
\end{lemma}
\begin{proof}
We choose decreasing radii $r_1 > r_2 > r_3 > \cdots$ so that $r_k \rightarrow R_0$. Let us first define initial data $(u_0,u_1)$ in the exterior region $\{x: |x|>R_0\}$. 
\[
 (u_0(r), u_1(r)) = (u^{(r_k)} (r,0), u_t^{(r_k)} (r,0)), \qquad r>r_k. 
\]
We claim that 
\begin{equation} \label{finite energy outside passing}
 \int_{|x|>R_0} (|\nabla u_0(x)|^2 + |u_1(x)|^2) dx < +\infty. 
\end{equation}
We start by considering the radiation profile $G_k^\pm $ of the initial data of $u^{(r_k)}$. An application of Corollary \ref{comparison of profile} shows 
\[
 \|G_k^\pm - G_v^\pm\|_{L^2([r_k, +\infty))} \lesssim_1 \|\chi_{r_k} F(x,t,u)\|_{L^1 L^2} \lesssim_\gamma \|\chi_{R_0} u\|_{Y(\Rm)}^5. 
\]
Here $G_v^\pm$ are the radiation profiles of $v_L$. Therefore we have the uniform boundedness
\[
 \|G_k^\pm\|_{L^2([r_k,+\infty))} \lesssim_\gamma \|G_v^+\|_{L^2(\Rm)} + \|\chi_{R_0} u\|_{Y(\Rm)}^5 < +\infty.
\]
We then utilize the explicit formula of initial data in term of radiation profiles and obtain 
\begin{align}
 \left\|(u^{(r_k)} (x,0), u_t^{(r_k)} (x,0)) - (a_k |x|^{-1},0)\right\|_{\mathcal{H}_{r_k}} & \lesssim_1 \|G_k^-\|_{L^2(\{s: |s|>r_k\})} \nonumber\\
 &\lesssim_\gamma \|G_v^+\|_{L^2(\Rm)} + \|\chi_{R_0} u\|_{Y(\Rm)}^5. \label{outer estimate urk}
\end{align}
Here the constants $a_k$'s are defined by 
\[
 a_k = \int_{-r_k}^{r_k} G_k^- (s) ds. 
\]
Because all $u^{(r_k)}$ coincide in the exterior region $\Omega_R$. The inequality above gives
\[
 \|(a_k - a_j)(|x|^{-1}, 0)\|_{\mathcal{H}_R} \lesssim_\gamma \|G_v^+\|_{L^2(\Rm)} + \|\chi_{R_0} u\|_{Y(\Rm)}^5 < +\infty, \qquad \forall j,k \geq 1.
\]
Therefore $a_k$ are also uniformly bounded. By passing to a subsequence if necessary we may assume $a_k \rightarrow a$. We then fix $r>R_0$, let $k\rightarrow +\infty$ in \eqref{outer estimate urk}, and obtain 
\begin{equation} \label{limit behaviour of u}
  \left\|(u_0(x), u_1 (x)) - (a |x|^{-1},0)\right\|_{\mathcal{H}_{r}} \lesssim_\gamma \|G_v^+\|_{L^2(\Rm)} + \|\chi_{R_0} u\|_{Y(\Rm)}^5.
\end{equation}
If $R_0>0$, we may make $r\rightarrow R_0$ and finish the proof of \eqref{finite energy outside passing}. The same argument still works in the case $R_0=0$ if we can show that $a=0$. In fact, we may combine \eqref{limit behaviour of u}, Strichartz estimates and finite speed of propagation to write
\[
 \left\|\chi_r(u - a|x|^{-1})\right\|_{Y(\Rm)} \lesssim_\gamma \|G_v^+\|_{L^2(\Rm)} + \|\chi_{R_0} u\|_{Y(\Rm)}^5, \qquad \forall r>R_0.
\]
By considering the limit $r\rightarrow 0^+$ we may verify that $a=0$ if $R_0 = 0$. This finally proves \eqref{finite energy outside passing}. Next if $R_0 > 0$ we define $(u_0,u_1)$ in the interior region $\{x: |x|<R_0\}$ by 
\[
 (u_0(x), u_1(x)) = (u_0(R_0), 0), \qquad |x|<R_0.
\]
We verify that $u$ is indeed an exterior solution to (CP1) with initial data $(u_0,u_1)$. In fact, if $\tilde{u}$ is the exterior solution of (CP1) in the region $\Omega_{R_0}$ with initial data $(u_0,u_1)$, then finite speed of propagation of wave equation immediately gives 
\[ 
 u(x,t) = \tilde{u}(x,t), \qquad (x,t) \in \Omega_{R_0}, \; t\in (-T_-,T_+).
\]
Here $(-T_-,T_+)$ is the maximal lifespan of $\tilde{u}$. We then obtain $T_-=T_+=+\infty$ by finite time blow-up criterion, since $\|\chi_{R_0} \tilde{u}\|_{Y((-T_-,T_+))} = \|\chi_{R_0} u\|_{Y((-T_-,T_+))} < +\infty$. As a result, $u$ is an exterior solution to (CP1) with initial data $(u_0,u_1)$. Finally our assumption $\|\chi_{R_0} u\|_{Y(\Rm)} < +\infty$ implies that $u$ must be $R_0$-weakly asymptotically equivalent to a free wave, by Lemma \ref{equivalent to exterior scattering}. The way we construct $u$ guarantees that $u$ must be $R_0$-weakly asymptotically equivalent to $v_L$. 
\end{proof}

\begin{lemma} \label{blow up of asymptotic equivalent}
 Assume $R > 0$. Let $u$ be an exterior solutions to (CP1) defined in the regions $\Omega_{R}$ and $v_L$ be a free wave with a finite energy, so that they are $R$-weakly asymptoticly equivalent to each other.  Then there exists a radius $r \in (0,R)$ so that we may extend the domain of $u$ to $\Omega_r$ so that $u$ is an exterior solution to (CP1), which is $r$-weakly asymptoticly equivalent to $v_L$. 
\end{lemma}
\begin{proof}
By Lemma \ref{equivalent to exterior scattering}, we must have $\|\chi_R u\|_{Y(\Rm)} < +\infty$. For convenience we use the same notations $u$ for its corresponding global extensions. It solves the equation
 \begin{align*}
  \partial_t^2 u - \Delta u = \chi_{R} F(x,t,u).
 \end{align*}
Please note that our assumptions guarantee that 
\begin{align*}
 \|u\|_{Y(\Rm)}, \|\chi_{R} F(x,t,u)\|_{L^1 L^2(\Rm \times \Rm^3)} < +\infty. 
\end{align*}
 Let $G_u^\pm, G_v^\pm \in L^2(\Rm)$ be the nonlinear radiation profiles of $u$ and the radiation profiles of $v_L$, respectively. We choose constant $\eps = \eps(\gamma) \ll 1$, and $r\in (R/2, R)$, so that 
 \begin{align}
  &\|\chi_{r, R} u\|_{Y(\Rm)} < \eps;& &\|G_u^\pm - G_v^\pm\|_{L^2([r, R])} < \eps/2.&  \label{assumption delta}
 \end{align}
 We then consider a complete distance space $X = \{G\in L^2(\{s: r<|s|<R\}): \|G\|_{L^2} < 2\eps\}$ whose distance is defined by the $L^2$ norm of the difference. Given $G \in X$, we first extend its domain to $\Rm$ by 
 \begin{align*}
  &G(s) = 0, \quad |s|\geq R;& &G(s) = \frac{-1}{2r} \int_{r<|s'|<R} G(s') ds', \quad |s|<r. &
 \end{align*}
 Let $(w_0,w_1)\in \dot{H}^1 \times L^2(\Rm^3)$ and $w_L$ be initial data and the corresponding linear free wave so that the radiation profile of $(w_0,w_1)-(u_0,u_1)$ is exactly $G$. The way we construct $G$ implies that $(w_0,w_1) = (u_0,u_1)$ in the exterior region $\{x: |x|>R\}$, by the explicit formula of initial data in term of radiation profiles. We consider the solution to 
 \[
  \partial_t^2 w - \Delta w = \chi_{r} F(x,t,w). 
 \]
 We claim that $w$ is well-defined for all time $t$. In fact, if $J = [-T_1,T_2]$ is an interval contained in the maximal lifespan of $w$, then the Strichartz estimate gives ($C_1$ is a constant)
 \[
  \|\chi_{r, R} (w-u)\|_{Y(J)} \leq \|\chi_{r, R} (w_L -u_L)\|_{Y(J)} + C_1 \|\chi_{r} F(x,t,w) - \chi_{R} F(x,t,u)\|_{L^1 L^2(J \times \Rm^3)}.
 \]
 By finite speed of propagation, we have $w = u$ in the region $\Omega_{R}$. Therefore we have 
 \begin{align*}
  \|\chi_{r, R} (w-u)\|_{Y(J)} &\leq \|\chi_{r, R} (w_L -u_L)\|_{Y(J)} + C_1 \|\chi_{r, R} F(x,t,w)\|_{L^1 L^2(J \times \Rm^3)}\\
  & \leq  \|\chi_{r, R} (w_L -u_L)\|_{Y(J)}  + C_1 \gamma \|\chi_{r, R} w\|_{Y(J)}^5. 
 \end{align*}
 Thus 
 \[
  \|\chi_{r, R} w\|_{Y(J)} \leq \|\chi_{r, R} (w_L -u_L)\|_{Y(J)} + \|\chi_{r,R} u\|_{Y(J)} + C_1 \gamma \|\chi_{r, R} w\|_{Y(J)}^5.
 \]
 Here the radiation profile of $w_L - u_L$ is $G$, thus the Strichartz estimates and our assumption \eqref{assumption delta} give
 \[
  \|\chi_{r, R} (w_L -u_L)\|_{Y(J)} +  \|\chi_{r,R} u\|_{Y(J)}  \lesssim_1 \|G\|_{L^2}  +  \|\chi_{r,R} u\|_{Y(\Rm)} \lesssim_1  \eps.
 \]
 Since $\eps = \eps(\gamma) \ll 1$, a continuity argument shows that
 \[
  \|\chi_{r, R} w\|_{Y(J)} \leq 2\|\chi_{r, R} (w_L -u_L)\|_{Y(J)} +2 \|\chi_{r,R} u\|_{Y(J)} \lesssim_1 \eps. 
 \]
 Therefore $\|\chi_r w\|_{Y(J)}$ is uniformly bounded for any interval time interval $J$ in the maximal span of $w$. This implies that $w$ is defined for all $t$. In addition, we have 
 \begin{equation} \label{upper bound of w R delta}
  \|\chi_{r,R} w\|_{Y(\Rm)} \lesssim_1 \eps. 
 \end{equation} 
 Let $G_w^\pm$ be the nonlinear radiation profiles of $w$. Since $w-u$ solves the equation 
 \[
  \partial_t^2 (w-u) - \Delta (w-u) = \chi_{r} F(x,t,w) - \chi_{R} F(x,t,u),
 \]
 by Lemma \ref{scatter profile of nonlinear solution} we have 
 \begin{align*}
  \|G_w^+ - G_u^+ + G(-s)\|_{L^2([r, +\infty))} & \lesssim_1 \|\chi_{r} F(x,t,w) - \chi_{R} F(x,t,u)\|_{L^1 L^2(\Rm \times \Rm^3)} \\
  & \lesssim_\gamma \|\chi_{r, R} w\|_{Y(\Rm)}^5 \lesssim_\gamma \eps^5. 
 \end{align*}
 Similarly we have 
 \begin{align*}
  \|G_w^- - G_u^- - G\|_{L^2([r, +\infty))} \lesssim_\gamma \eps^5. 
 \end{align*}
 Next we define a map $\mathbf{T}$ on  $X$  
 \[
  (\mathbf{T} G)(s) = \left\{\begin{array}{ll} G(s) - G_w^- (s) + G_v^-(s), & s \in (r,R); \\ G(s) + G_w^+(-s) - G_v^+(-s), & s\in (-R, -r).\end{array}\right. 
 \]
 This map is actually a map from $X$ to itself if $\eps=\eps(\gamma)$ is sufficiently small. In fact, we have ($C_\gamma$ is a constant determined by $\gamma$)
 \begin{align*}
  \|G-G_w^- + G_v^-\|_{L^2((r,R))} & \leq \|G-G_w^- + G_u^-\|_{L^2((r,R))} + \|G_u^- - G_v^-\|_{L^2((r,R))}\\
  & \leq  C_\gamma \eps^5 + \eps/2.
 \end{align*}
 The upper bound of $\|G(s) + G_w^+(-s) - G_v^+(-s)\|_{L^2((-R, -r))}$ is the same. Therefore (please note $\eps = \eps(\gamma)\ll 1$)
 \[
  \|\mathbf{T} G\|_{L^2(\{s: r < |s|<R\})} \leq 2 C_\gamma \eps^5 + \eps < 2\eps. 
 \]
 We claim that $\mathbf{T}$ is a contraction map on $X$. This immediately gives a fixed-point $G \in X$. The corresponding exterior solution $w$ then satisfies $G_w^\pm (s)= G_v^\pm(s)$ for $s\in (r, R)$, thus is indeed weakly asymptoticly equivalent to $v$ in $\Omega_{r}$. The remaining task is to show $\mathbf{T}$ is a contraction map on $X$. We let $\tilde{G}$ be another element in $X$ and define $\tilde{w}$, $\tilde{w}_L$, $G_{\tilde{w}}^\pm$ accordingly. We utilize the Strichartz estimate and \eqref{upper bound of w R delta}, follow the same argument as above, then obtain
\[
 \|\chi_{r, R} (w-\tilde{w})\|_{Y(\Rm)} \leq 2 \|\chi_{r, R} (w_L -\tilde{w}_L)\|_{Y(\Rm)} \lesssim_1 d(G,\tilde{G}). 
\]
We then consider the nonlinear radiation profiles
\begin{align*}
 \|G_w^+ - G_{\tilde{w}}^+ + G(-s) - \tilde{G}(-s)\|_{L^2([r, +\infty))} &\lesssim_1\|\chi_{r} (F(x,t,w)-F(x,t,\tilde{w}))\|_{L^1 L^2(\Rm \times \Rm^3)}\\
 & \lesssim_1 \|\chi_{r,R} (F(x,t,w)-F(x,t,\tilde{w}))\|_{L^1 L^2(\Rm \times \Rm^3)}\\
 & \lesssim_\gamma \eps^4 d(G,\tilde{G}).
\end{align*}
Similarly we have 
\[
 \|G_w^- - G_{\tilde{w}}^- - G + \tilde{G}\|_{L^2([r, +\infty))} \lesssim_\gamma \eps^4 d(G,\tilde{G}).
\]
Therefore we have $d(\mathbf{T} G, \mathbf{T} \tilde{G}) \lesssim_\gamma \eps^4 d(G, \tilde{G})$. Our choice $\eps = \eps(\gamma) \ll 1$ then implies that $\mathbf{T}$ is a contraction map. 
\end{proof}

\subsection{Uniqueness} 

In this subsection we show that every weakly asymptotically equivalent exterior solution of $v_L$ is covered by the one-parameter family defined earlier. In fact, given a standard radial weakly asymptotically equivalent solution to a free wave, the asymptotically equivalent solution to it with a given characteristic number is essentially unique. We have 

\begin{lemma}
 Assume that $u$ and $v$ are radial exterior solutions to (CP1) in $\Omega_R$ so that they are both $R$-weakly asymptoticly equivalent to the same finite-energy free wave. If the characteristic numbers of $u$ respect to $v$ is zero, then we have 
\[
 u(x,t) = v(x,t), \qquad \forall |x|>|t|+R. 
\]
\end{lemma}
\begin{proof}
We use the same notation as in the proof of Lemma \ref{L1 lemma}. It suffices to show that the identity $u(x,t) = v(x,t)$ holds in $\Omega_{R_1}$ for a sufficiently large radius $R_1>R$, thanks to Lemma \ref{coincidence inside}. We recall the radiation profiles $G$ and $\tilde{G}$ satisfy ($C$ is a constant)
\[ 
 \|G-\tilde{G}\|_{L^2\{s: |s|>2^k\}} \leq  C \sum_{j\geq k}  2^{-j/2} b_j^4, \qquad \forall k \gg 1. 
\]
We then use the explicit formula of $w_L = \mathbf{S}_L (u_0-v_0,u_1-v_1)$ to give an upper bound when $r > 2^k + |t|$:
\begin{align*}
 |w_L (r,t)| = \left|\frac{1}{r} \int_{t-r}^{t+r} (G(s) - \tilde{G}(s)) ds\right| \leq \frac{1}{r} \sum_{j\geq k} \int_{2^j < |s| < 2^{j+1}} |G(s)-\tilde{G}(s)| ds.
\end{align*}
Please note that here our assumption on the characteristic number implies that $\int_\Rm (G-\tilde{G}) ds = 0$. Inserting the $L^2$ norm given above, we obtain 
\[
 |w_L(r,t)| \lesssim_1 \frac{C}{r} \sum_{j\geq k} b_j^4, \qquad r>|t|+2^k
\]
Thus when $k$ is sufficiently large, we have 
\[
 \|\chi_{2^k} w_L\|_{Y(\Rm)} \lesssim_1 C \cdot 2^{-k/2} \sum_{j\geq k} b_j^4. 
\]
We choose $K \gg 1$ so that the inequalities above hold for $k \geq K$ and 
\[
 \sum_{j \geq K} b_j^4 \ll 1 \qquad \Rightarrow \qquad \|\chi_{2^K} u\|_{Y(\Rm)},  \|\chi_{2^K} v\|_{Y(\Rm)} \ll 1.
\]
Thus if $k \geq K$ is sufficiently large, then we may apply Lemma \ref{nonlinear comparison by linear} and obtain
\[
 \|\chi_{2^k} (u-v)\|_{Y(\Rm)} \lesssim_1 \|\chi_{2^k} w_L\|_{Y(\Rm)} \lesssim_1 C \cdot 2^{-k/2} \sum_{j\geq k} b_j^4\lesssim_1 C \cdot 2^{-k/2} \sum_{j\geq K} b_j^4.
\]
We then apply Corollary \ref{comparison of profile} and give a new upper bound if $k \geq K$: 
\begin{align*}
 \|G-\tilde{G}\|_{L^2(\{s: |s|>2^k\})} & \lesssim_1 \|\chi_{2^k} (F(x,t,u)-F(x,t,v))\|_{L^1 L^2 (\Rm \times \Rm^3)} \\
 & \lesssim_1 \sum_{j \geq k} \|\tilde{\chi}_{2^j} (F(x,t,u)-F(x,t,v))\|_{L^1 L^2 (\Rm \times \Rm^3)} \\
 & \lesssim_\gamma  \sum_{j \geq k}  (\|\tilde{\chi}_{2^j} u\|_{Y(\Rm)}^4 + \|\tilde{\chi}_{2^j} v\|_{Y(\Rm)}^4) \|\tilde{\chi}_{2^j} (u-v)\|_{Y(\Rm)}\\
 & \lesssim_\gamma C \left(\sum_{j\geq K} b_j^4\right) \sum_{j \geq k} 2^{-j/2} b_j^4. 
\end{align*}
We may iterate this argument and obtain that there exists a constant $c_\gamma$ so that 
\[
 \|G-\tilde{G}\|_{L^2(\{s: |s|>2^k\})} \leq \left(c_\gamma \sum_{j\geq K} b_j^4\right)^N C \sum_{j \geq k} 2^{-j/2} b_j^4, \qquad \forall N \geq 1, \, k\geq K.
\]
We may let $N \rightarrow +\infty$ and conclude that $\|G-\tilde{G}\|_{L^2\{s: |s|>2^K\}} = 0$ as long as $K$ is sufficiently large. The explicit formula of initial data in term of radiation profile then gives the coincidence of initial data in the exterior region $\{x: |x|>2^K\}$, which in turn gives $u(x,t)= v(x,t)$ in $\Omega_{2^K}$ by finite speed of propagation. 
\end{proof}

\begin{corollary} \label{coincidence with same characteristic} 
Assume that a radial exterior solution $v$ is weakly asymptotically equivalent to a finite-energy free wave $v_L$. If radial exterior solutions $u$, $\tilde{u}$ are both $R$-weakly asymptotically equivalent to $v_L$, so that the characteristic numbers of $u$ and $\tilde{u}$ respect to $v$ are the same, then we have 
\[
 u(x,t) = \tilde{u}(x,t), \qquad \forall |x|>|t|+R. 
\]
\end{corollary}
\begin{proof}
A basic observation is that the characteristic number of $u$ with respect to $\tilde{u}$ is zero. We then apply the lemma above.  
\end{proof}

\noindent Finally we show that all radial exterior solutions $u$ that is weakly asymptotically equivalent to $v_L$ are covered by the one-parameter family $\{(u^\alpha, R_\alpha)\}$. Let $u$ be a radial exterior solutions that is $R$-weakly asymptotically equivalent to $v_L$. We let $\alpha$ be its characteristic number (with respect to $u^0$). If $R > R_\alpha$, we may apply Corollary \ref{coincidence with same characteristic} in $\Omega_R$ to conclude that 
\[
 u(x,t) = u^\alpha (x,t), \qquad \forall (x,t)\in \Omega_R. 
\]
The remaining work is to show that $R \leq R_\alpha$ is impossible. In fact, if $R\leq R_\alpha$ held, then $u^\alpha$ would have to be a blow-up asymptotically equivalent solution. We might apply Corollary \ref{coincidence with same characteristic} in $\Omega_{r_k}$ with a decreasing sequence of radii $r_k \rightarrow R_\alpha$ and obtain $u(x,t) = u^\alpha (x,t)$ for all $(x,t) \in \Omega_{R_\alpha}$. Therefore $\|\chi_{R_\alpha} u\|_{Y(\Rm)} = \|\chi_{R_\alpha} u^\alpha\|_{Y(\Rm)} = +\infty$. This is a contradiction, because our assumption guarantees $\|\chi_{R} u\|_{Y(\Rm)} < +\infty$ by Lemma \ref{equivalent to exterior scattering}. 

\subsection{Primary Asymptotically Equivalent Solutions} \label{sec: primary}

If $v_L = \mathbf{S}_L(v_0,v_1)$ is a radial linear free wave so that 
\[
 \delta^2 \doteq \sum_{k=-\infty}^\infty \|\tilde{\chi}_{2^k} v_L\|_{Y(\Rm)}^2 < \delta(\gamma)^2
\]
is sufficiently small, then Lemma \ref{asymptoticly equivalent to square small} immediately gives an asymptotically equivalent solution $u$ to $v_L$ defined in $\Omega_0$ so that 
\begin{equation} \label{smallness in family}
 \|(u(\cdot,0),u_t(\cdot,0))-(v_0,v_1)\|_{\dot{H}^1\times L^2(\Rm^3)} \lesssim_\gamma \|\chi_0 v_L\|_{Y(\Rm)}^5 \lesssim_\gamma \delta^5. 
\end{equation}
We then recall Proposition \ref{l2 linear} and Lemma \ref{l2 small nonlinear} to conclude that 
\[
 \sum_{k=-\infty}^\infty \|\tilde{\chi}_{2^k} u\|_{Y(\Rm)}^2 \lesssim_\gamma \delta^2,
\]
as long as $\delta$ is sufficiently small. Let $\{(u^\alpha, R_\alpha)\}$ be the one-parameter family of weakly asymptotically equivalent solutions to $v_L$ with $u^0 = u$. For any given $\alpha \neq 0$, we may apply corollary \ref{lower bound far outside} and conclude that there exists $R > R_\alpha$, so that 
\[
 \|(u^\alpha(\cdot,0), u_t^\alpha(\cdot.0))- (u(\cdot,0),u_t(\cdot,0))\|_{\dot{H}^1\times L^2(\{x: |x|>R\})} \gtrsim_\gamma 1. 
\]
As a result, $u$ is the only solution satisfying the inequality
\[
 \|(u(\cdot,0),u_t(\cdot,0))-(v_0,v_1)\|_{\dot{H}^1\times L^2(\Rm^3)} \ll 1
\]
in the one-parameter family of weakly asymptotically equivalent solutions, as long as $\delta$ is sufficiently small. We call $u$ the primary asymptotically equivalent solution of $v_L$ and the one-parameter family $\{u^\alpha, R_\alpha\}_{\alpha \in \Rm}$ defined above the primary one-parameter family associated to $v_L$.
 
\subsection{Time Translations} 

Before we conclude this section, we consider the characteristic numbers of time-translated solutions. This will be very useful when we consider non-radiative solutions in the next section. Let $u$ and $v$ be two exterior solution of (CP1) defined in $\Omega_R$ so that they are $R$-weakly asymptoticly equivalent to each other. Given any time $t_0 \in \Rm$, it is clear that the time-translated versions $u(x,t+t_0)$ and $v(x,t+t_0)$ are both exterior solutions of the modified equation 
\[
 \partial_t^2 w - \Delta w = F(x,t+t_0, w)
\]
in the region $\Omega_{R+|t_0|}$. In addition, the time-translated versions of solutions are still weakly asymptoticly equivalent to each other. A natural question is whether the characteristic number $\alpha'$ of $u(x,t+t_0)$ respect to $v(x,t+t_0)$ remains the same number as the original characteristic number $\alpha$ of $u$ respect to $v$. The answer is affirmative, this can be verified by considering the asymptotic behaviour of the difference of the data of these two solutions. We apply Lemma \ref{L1 lemma} on both the original solutions and time-translated solutions, thus
\begin{align*}
 \|(u(\cdot,t_0)-v(\cdot,t_0),u_t (\cdot,t_0)-v_t(\cdot,t_0))-(\alpha |x|^{-1},0)\|_{\dot{H}^1 \times L^2(\{x: |x|>r+|t_0|\})} & = o(r^{-1/2});\\
 \|(u(\cdot,t_0)-v(\cdot,t_0),u_t (\cdot,t_0)-v_t(\cdot,t_0))-(\alpha' |x|^{-1},0)\|_{\dot{H}^1 \times L^2(\{x: |x|>r\})} & = o(r^{-1/2}).
\end{align*}
This implies that $\alpha' = \alpha$ since we have $\|(|x|^{-1},0)\|_{\dot{H}^1(\{x: |x|>r\})} \simeq r^{-1/2}$.  

\section{Non-radiative Solutions}

Let us consider the special case, when scattering profile of the solutions are exactly zero. In the remaining part of this work, we always assume $F(t,x,u)$ is independent of the time $t$. We have 
\begin{lemma} \label{independence of time}
 Assume that $F(x,t,u)$ satisfies (A1)-(A3). Let $u$ be an exterior solution to (CP1) in the region $\Omega_R$ so that 
 \[
  \lim_{t\rightarrow \pm \infty} \int_{|x|>R+|t|} |\nabla_{t,x} u(x,t)|^2 dx = 0. 
 \]
 Then $u(x,t)$ is independent of $t$ in the region $\Omega_R$. 
\end{lemma}
\begin{proof}
 Clearly the zero solution is another solution to (CP1), which is weakly asymptoticly equivalent to zero. The results in the previous section show that the exterior solutions $u(x,t+t_0)$ and $u(x,t)$ defined in the region $\Omega_{R+|t_0|}$ shares the same characteristic number with respect to the zero solution. By uniqueness we have 
 \[
  u(x,t+t_0) = u(x,t), \qquad |x|>|t|+R+|t_0|.
 \]
 This finishes the proof. 
\end{proof}

\paragraph{One-parameter family} Now we consider the one-parameter family $\{(u^\alpha, R_\alpha)\}_{\alpha \in \Rm}$ of non-radiative solutions. We always choose $u^0$ to be the zero solution in the case of non-radiative solutions. We first recall the result of Corollary \ref{lower bound far outside}:
\begin{equation*} 
 \|u^\alpha\|_{\dot{H}^1(\{x:|x|>R\})} \simeq_1 |\alpha| R^{-1/2}, \qquad R\gtrsim_\gamma |\alpha|^2. 
\end{equation*}
Next we consider additional properties of blow-up/scattering non-radiative solutions besides the properties given in Theorem \ref{main 1}, and prove Theorem \ref{main non-radiative}. We first consider blow-up non-radiative solutions, i.e. those $u^\alpha$'s with $R_\alpha\geq 0$ and $\|\chi_{R_\alpha} u^\alpha\|_{Y(\Rm)} = +\infty$. We start by

\begin{lemma} \label{Y by Hdot1}
 Let $v \in \dot{H}^1(\Rm^3)$ be radial. If we view $v$ as a function of $(x,t) \in \Rm^3 \times \Rm$ (independent of $t$), then we have the inequlaity
 \[
  \|\chi_0 v\|_{Y(\Rm)} \lesssim_1 \|v\|_{\dot{H}^1(\Rm^3)}. 
 \]
\end{lemma}
\begin{proof}
 We let $a_k = \|\nabla v\|_{L^2(\{x: 2^k < |x|<2^{k+1}\})}$ and $\Phi_k = \{(x,t): 2^k < |x| < 2^{k+1}, |t|<|x|\}$. We have 
 \[
  \sup_{(x,t) \in \Phi_k} |v(x)| \leq \sum_{j \geq k} \int_{2^j}^{2^{j+1}} |v_r| dr 
  \leq \sum_{j \geq k} \left(\int_{2^j}^{2^{j+1}} \frac{dr}{r^2}\right)^{1/2} \left(\int_{2^j}^{2^{j+1}} r^2 |v_r|^2 dr\right)^{1/2} \lesssim_1 \sum_{j \geq k} 2^{-j/2} a_j. 
 \]
 Thus we have 
 \[
  \|\chi_{\Phi_k} v\|_{Y(\Rm)} \lesssim_1 \sum_{j\geq k} 2^{-\frac{j-k}{2}} a_j.
 \]
 We then finish the proof by Young's inequality 
 \[
  \|\chi_0 v\|_{Y(\Rm)} \leq \left(\sum_{k=-\infty}^\infty \|\chi_{\Phi_k} v\|_{Y(\Rm)}^5 \right)^{1/5} \leq \left(\sum_{k=-\infty}^\infty \|\chi_{\Phi_k} v\|_{Y(\Rm)}^2\right)^{1/2} \lesssim_1 \|\{a_k\}\|_{l^2} = \|v\|_{\dot{H}^1(\Rm^3)}. 
 \]
\end{proof}

\begin{corollary} \label{non-radiative Y by H}
 Let $u^\alpha$ be a member of the one-parameter family of non-radiative solutions to (CP1). Then we have 
 \[
  \|\chi_R u^\alpha\|_{Y(\Rm)} \lesssim_1 \|u^\alpha\|_{\dot{H}^1(\{x: |x|>R\})}, \qquad \forall R \geq R_\alpha.
 \]
 Therefore all blow-up non-radiative solutions satisfy $\|u^\alpha\|_{\dot{H}^1(\{x: |x|>R_\alpha\})} = +\infty$. 
\end{corollary}
\begin{proof}
 It suffice to consider the radial $\dot{H}^1(\Rm^3)$ function $v$ defined by 
 \[
  v(x) = \left\{\begin{array}{ll} u^\alpha(x), & |x|>R; \\ u^\alpha(R), & |x|\leq R; \end{array}\right.
 \]
 and apply Lemma \ref{Y by Hdot1}. 
\end{proof}

\paragraph{ODE method} All the radial non-radiative solutions can be obtained by a method of ordinary differential equation. This will be very useful in the discussion concerning further properties of non-radiative solutions. We consider the one-variable function $w(r) = r u(r)$. If $u$ solves the elliptic equation $-\Delta u = F(x,u)$, then $w$ solves the elliptic equation 
\begin{equation} \label{one-dimensional elliptic w}
 - w_{rr} = F(r, u(r)) = F(r, w(r)/r).
\end{equation}
In addition, the asymptotic behaviour of $u$ gives that $w(r) \rightarrow \alpha$ as $r\rightarrow +\infty$. We may solve this ordinary differential equation by a fixed-point argument near infinity. We consider the distance space 
\begin{align*}
 &X = \left\{w \in C([R, +\infty)): |w(r)|\leq 2|\alpha|, \, \forall r\geq R\right\};& &d(w_1,w_2) = \sup_{r\geq R} |w_1(r)-w_2(r)|.&
\end{align*}
A straight-forward calculation shows that if $R$ is sufficiently large, then the map (a similar argument can be found in Shen \cite{shen2})
\[
 (\mathbf{T} w)(r) = \alpha - \int_r^\infty \int_s^\infty F(s, w(s)/s) ds dr 
\]
is a contraction map from $X$ to $X$. The unique fixed-point $w$ is a $C^{1,1}$ solution to the elliptic equation \eqref{one-dimensional elliptic w} for $r\geq R$. We then solve this equation backward from $r= R$ and obtain a solution $w$ with a maximal domain $(R_0, +\infty)$. Please note that if $R_0>0$, then we must have $\|w_r\|_{L^2((R_0,+\infty))} = +\infty$ by the basic theory of ordinary differential equations. We then consider the function $u(x)=w(|x|)/|x|$. A basic calculation shows that $u$ is an exterior solution to (CP1) in $\Omega_r$ for all $r > R_0$, which is also non-radiative with characteristic number $\alpha$. Thus we have $u = u^\alpha$ as long as $|x|>\max\{R_0, R_\alpha\} + |t|$. Next we observe the facts 
\begin{align*}
 &\|u\|_{\dot{H}^1(\{x: |x|>r\})} < +\infty, \, \forall r > R_0;& &\|u\|_{\dot{H}^1(\{x: |x|>R_0\})}\gtrsim_1 \|w_r\|_{L^2((R_0,+\infty))}.&
\end{align*}
We claim that $R_0 = R_\alpha$ because 
\begin{itemize}
 \item If $R_0 > R_\alpha$, then we would obtain $R_0>0$ and 
 \[
  \|u^\alpha\|_{\dot{H}^1(\{x: |x|>R_0\})} = \|u\|_{\dot{H}^1(\{x: |x|>R_0\})}\gtrsim_1 \|w_r\|_{L^2((R_0,+\infty))} = +\infty.
 \]
 This is a contradiction.  
 \item If $R_0 < R_\alpha$, then we would obtain $R_\alpha > 0$ and 
 \[
  \|\chi_{R_\alpha} u^\alpha\|_{Y(\Rm)} \lesssim_1 \|u^\alpha\|_{\dot{H}^1(\{x: |x|>R_\alpha\})} = \|u\|_{\dot{H}^1(\{x: |x|>R_\alpha\})} < +\infty.
 \]
 This contradicts with the basic property of blow-up asymptotically equivalent solutions. 
\end{itemize} 
Therefore the solution $u^\alpha$ are completely given by a $C_{loc}^{1,1}$ solution $w$ of the elliptic equation \eqref{one-dimensional elliptic w}.

\paragraph{Scattering non-radiative solutions} The only remaining work to prove Theorem \ref{main non-radiative} is to show that all scattering non-radiative solutions, i.e. those $u^\alpha$'s satisfying $\|\chi_0 u^\alpha\|_{Y(\Rm)} < +\infty$, are actually (or extend naturally to) solutions to (CP1) defined in the whole space-time $\Rm^3 \times \Rm$. We recall that the solution $u^\alpha$ is given by a solution $w$ to the elliptic equation \eqref{one-dimensional elliptic w}.  Our radial assumption and the fact $u^\alpha \in \dot{H}^1(\Rm^3)$ imply that $r^{1/2} u^\alpha (r) \rightarrow 0$ as $r\rightarrow 0^+$. In addition, we have 
\[
 \int_0^\infty r^2 |u_r^\alpha (r)|^2 dr \simeq \|u^\alpha\|_{\dot{H}^1(\Rm^3)}^2 < +\infty \qquad \Rightarrow \qquad \liminf_{r\rightarrow 0^+} r^{3/2} u_r^\alpha (r) = 0. 
\]
By $w_r = u^\alpha + ru_r^\alpha$ we also have
\[
 \liminf_{r\rightarrow 0^+} r^{1/2} w_r (r) = 0.
\]
Thus given any positive constant $c_0 \ll 1$, we may find a radius $r_0  \ll 1$, so that 
\begin{align*}
 &|u^\alpha (r)| \leq 2c_0 r^{-1/2}, \; \forall r\in (0,r_0];& &|w_r (r_0)| \leq c_0 r_0^{-1/2}.&
\end{align*}
The first inequality above gives $w_{rr} \leq 32\gamma c_0^5 r^{-3/2}$ for all $r\in (0,r_0]$. Therefore we have 
\begin{align*}
 |w_r(r)| & \leq |w_r(r_0)| + \int_r^{r_0} |w_{rr}(r')| dr' \leq \left(64\gamma c_0^5 + 2^{-1/2} c_0 \right) r^{-1/2}, & & \forall r\in (0,r_0/2];\\
 |u^\alpha (r)| & \leq r^{-1} \int_0^r |w_r(r')| dr' \leq 2 \left(64\gamma c_0^5 + 2^{-1/2} c_0 \right) r^{-1/2}, & & \forall r\in (0,r_0/2].
\end{align*}
We may iterate this argument and obtain
\begin{align*}
 &|u^\alpha (r)| \leq 2c_k r^{-1/2}, \; \forall r\in (0,2^{-k} r_0];& &|w_r (r)| \leq c_k r^{-1/2}, \; r\in (0,2^{-k} r_0].&
\end{align*}
Here $c_k$ is a sequence defined by the induction formula $c_{k+1} = 2^{-1/2} c_k + 64\gamma c_k^5$. As long as $c_0$ is sufficiently small, we have $c_k \lesssim 2^{-\beta k}$ holds for any given $\beta < 1/2$ and sufficiently large $k$. This gives the decay 
\[
 |u^\alpha (r)|, |w_r(r)| \lesssim r^{\beta-1/2}, \forall r \ll 1. 
\]
We may choose $\beta$ slightly smaller than $1/2$ and obtain $w_{rr} \rightarrow 0$ as $r \rightarrow 0^+$. This means that $w_r(r)$ converges to a finite number $\kappa$ as $r\rightarrow 0^+$. We then obtain the following asymptotic behaviours of $w$ and $u^\alpha$ near $r = 0$:
\begin{align*}
 &w(r) = \kappa r + O(r^3);& &w_r(r) = \kappa + O(r^2);& &w_{rr} = O(r).&\\
 &u^\alpha(r) = \kappa + O(r^2);& &u_r^\alpha (r) = O(r); & &u_{rr}^\alpha = O(1).&
\end{align*}
These imply that $u^\alpha \in C^1 \cap L^{10}(\Rm^3)$ with weak second derivatives in $L^2(\Rm^3)$ so that $u^\alpha$ solves the elliptic equation $-\Delta u = F(x,u)$ in the (almost everywhere) point-wise sense. This verifies that $u^\alpha$ is a solution to (CP1) defined in the whole space-time $\Rm^3 \times \Rm$.  

\section{Scattering of Solutions with A Priori Bound}

In this section we assume that $F$ satisfies (A1)-(A4). 

\begin{lemma} \label{partition of u with u alpha}
Assume that $c_1$ and $M$ are positive constants. There exists a positive constant $\delta_0 = \delta_0 (\gamma, M, c_1)$ so that if a non-radiative solution $u^\alpha$ defined in the previous section and an exterior solution $v$ defined in $\Omega_0$ satisfy: 
\begin{itemize}
 \item there exists a radius $R > \max\{R_\alpha, c_1 |\alpha|^2\}$ so that $\|\chi_R u^\alpha\|_{Y(\Rm)} \leq M$;
 \item the exterior solution $v$ satisfies
 \[
  \delta^2 \doteq \sum_{k=-\infty}^\infty \|\tilde{\chi}_{2^k} v\|_{Y(\Rm)}^2 < \delta_0^2;
 \] 
\end{itemize} 
 then there exists an exterior solution $u$ to (CP1) defined in $\Omega_R$, which is $R$-weakly asymptotically equivalent to $v$, with characteristic number $\alpha$ and satisfies 
\begin{align*}
  &\|\chi_R(u-u^\alpha)\|_{Y(\Rm)} \lesssim_{\gamma, M, c_1} \delta; & \|u(\cdot,0)-u^\alpha - v(\cdot,0)\|_{\dot{H}^1(\{x: |x|>R\})} \lesssim_{\gamma, M, c_1} \delta. 
\end{align*}
\end{lemma}
\begin{remark} 
 We may substitute the assumption $\|\chi_R u^\alpha\|_{Y(\Rm)} \leq M$ by $\|u^\alpha\|_{\dot{H}^1(\{x: |x|>R\})} \leq M$ because Corollary \ref{non-radiative Y by H} gives 
 \[
   \|\chi_R u^\alpha\|_{Y(\Rm)} \lesssim_1 \|u^\alpha\|_{\dot{H}^1(\{x: |x|>R\})}, \qquad \forall R \geq R_\alpha.
 \]
\end{remark}
\begin{remark} \label{uniqueness of uvw alpha}
 Our assumption on $v$ implies that $v$ is asymptotically equivalent to some linear free wave $v_L$. Therefore any exterior solution to (CP1), that is $R$-weakly asymptotically equivalent to $v$, with characteristic number $\alpha$ must coincide with the solution $u$ given in the Lemma above in the exterior region $\Omega_R$, thanks to Corollary \ref{coincidence with same characteristic}. In fact we may consider the one-parameter family $\{(v^\alpha, R'_\alpha)\}$ associated with this linear free wave $v_L$ with $v^0 = v$. Clearly $u = v^\alpha$ is the weakly asymptotically equivalent solution we are looking for in Lemma \ref{partition of u with u alpha}. 
\end{remark}
\begin{proof}
Without loss of generality we assume $\alpha\neq 0$, otherwise we may simply choose $u=v$. Let $u = v^\alpha$ be the member of one-parameter family given in Remark \ref{uniqueness of uvw alpha}. We still need to show that $R>R'_\alpha$ and verifies the inequalities. We use the notation $b_k = \|\tilde{\chi}_{2^k} v\|_{Y(\Rm)}$ and let $c = c(\gamma)$, $N = N(\gamma, c, \{b_k\})$ be constants given in Lemma \ref{existence of solutions far}. We recall that $u$ are initially constructed in the proof of Lemma \ref{existence of solutions far} with a domain $\Omega_{2^N}$. By choosing $\delta$ sufficiently small (depending on $\gamma$), we have $N \simeq_\gamma |\alpha|^2$. In addition, the non-radiative solution $u^\alpha$ is also initially constructed in $\Omega_{2^N}$ in the same manner (via a substitution of $v$ by zero, with the same constant $N$).  A review of the proof for Lemma \ref{existence of solutions far} reveals that the solutions $u$ and $u^\alpha$ satisfy
\begin{equation}
 \|\chi_{2^k} (u-v)\|_{Y(\Rm)},\, \|\chi_{2^k} u^\alpha\| \lesssim_1 |\alpha| 2^{-k/2}, \qquad \forall k \geq N. \label{upper bound on Y u new}
\end{equation}
Without loss of generality, we assume $R>R'_\alpha$, otherwise we may substitute $R$ by a radius $R'$ slightly larger than $R_\alpha$, then proceed. Our final conclusion shows that $\|\chi_{R'}(u-v)\|_{Y(\Rm)} \lesssim_{\gamma, c_1, M} \delta$. This immediately gives a contradiction if we make $R' \rightarrow R'_\alpha$ since we have $\|\chi_{R_\alpha} u\|_{Y(\Rm)} = +\infty$. Next we consider the solution $w = u-v-u^\alpha$ defined in $\Omega_R$ solving the equation
\[
 \partial_t^2 w - \Delta w = F(x,u) - F(x,u^\alpha) - F(x, v), \qquad (x,t) \in \Omega_R;
\]
Please note we may choose a pair of initial data $(w_0,w_1) \in \dot{H}^1 \times L^2(\Rm^3)$ for $w$, because 
\begin{itemize} 
 \item Our assumption $R>\max\{R'_\alpha,R_\alpha\}$ guarantees that the initial data of $w$ in the exterior region $\{x: |x|\geq R\}$ comes with a finite energy;
 \item The exterior solution $w$ is independent of $(w_0,w_1)$ in the interior region $\{x: |x|<R\}$ by finite speed of propagation.
\end{itemize}
Let $w_L$ be the finite-energy linear free wave with initial data $(w_0,w_1)$ and $G$ be the radiation profile of $w_L$. For convenience we also use the notations 
\begin{align*} 
  y_k = \|\chi_{2^k} w_L\|_{Y(\Rm)}.
\end{align*}
We may write $\chi_{2^k} w_L = \chi_{2^k} w_L^1 + \chi_{2^k} w_L^2$ for any $k\geq N$. Here $w_L^1, w_L^2$ are linear free waves whose initial data in the exterior region $\{x: |x|>2^N\}$ are $(u(\cdot,0), u_t(\cdot,0))-(v(\cdot,0), v_t(\cdot,0))$ and $(u^\alpha(\cdot,0), u_t^\alpha(\cdot,0))$, respectively. According to \eqref{upper bound of wL tem}, we have $\|\chi_{2^k} w_L^1\|_{Y(\Rm)}, \|\chi_{2^k} w_L^2\|_{Y(\Rm)} \lesssim_1 |\alpha| 2^{-k/2}$. Thus $y_k \lesssim_1 |\alpha| 2^{-k/2}$. The proof of the lemma consists of three steps: 

\paragraph{Step 1} We first show that the inequality $\|\chi_{2^N} w\| \lesssim_{\gamma} \delta$ holds. We start by applying Strichartz estimates and obtain ($k\geq N$)
\begin{align*}
 \|\chi_{2^k} w\|_{Y(\Rm)} & \leq  \|\chi_{2^k} w_L\|_{Y(\Rm)} + C \|\chi_{2^k} (F(x,u) - F(x,u^\alpha) - F(x, v))\|_{L^1 L^2 (\Rm \times \Rm^3)}\\
 & \leq \|\chi_{2^k} w_L\|_{Y(\Rm)} + C\|\chi_{2^k} (F(x,u) - F(x,v+u^\alpha))\|_{L^1 L^2(\Rm \times \Rm^3)} \\
 & \qquad + C\|\chi_{2^k} (F(x, v+u^\alpha)-F(x,v)-F(x,u^\alpha)) \|_{L^1 L^2(\Rm \times \Rm^3)}\\
 & \leq \|\chi_{2^k} w_L\|_{Y(\Rm)} + C\gamma \|\chi_{2^k} w\|_{Y(\Rm)} \left(\|\chi_{2^k} u\|_{Y(\Rm)}^4 + \|\chi_{2^k} (v+u^\alpha) \|_{Y(\Rm)}^4 \right)\\
 & \qquad + C\|\chi_{2^k} (F(x, v+u^\alpha)-F(x,v)-F(x,u^\alpha)) \|_{L^1 L^2}.
\end{align*}
Here $C$ is a constant. Since we have (if necessary, we slightly enlarge the value of $c = c(\gamma)$)
\[
 \|\chi_{2^N} u\|_{Y(\Rm)}, \|\chi_{2^N} u^\alpha\|_{Y(\Rm)}, \|\chi_{2^N} v\|_{Y(\Rm)} \ll 1,
\]
the inequality above implies that
\begin{equation} \label{equation 1 for uvw}
 \|\chi_{2^k} w\|_{Y(\Rm)} \leq 2 \|\chi_{2^k} w_L\|_{Y(\Rm)} + 2C\|\chi_{2^k} (F(x, v+u^\alpha)-F(x,v)-F(x,u^\alpha)) \|_{L^1 L^2}. 
\end{equation}
We then observe that the nonlinear radiation profiles of $w$ vanishes for all $s>R$ in both time directions. This gives ($k \geq N$)
\begin{align*}
 & \|G\|_{L^2(\{s: |s|>2^k\})} \\
 & \lesssim_1 \|\chi_{2^k} (F(x,u) - F(x,u^\alpha) - F(x, v))\|_{L^1 L^2 (\Rm \times \Rm^3)}\\
 & \lesssim_1 \|\chi_{2^k} (F(x,u) - F(x,v+u^\alpha))\|_{L^1 L^2} + \|\chi_{2^k} (F(x, v+u^\alpha)-F(x,v)-F(x,u^\alpha)) \|_{L^1 L^2}\\
 & \lesssim_\gamma \sum_{j\geq k} \|\tilde{\chi}_{2^j} w\|_{Y(\Rm)} \left(\|\tilde{\chi}_{2^j} u\|_{Y(\Rm)}^4 + \|\tilde{\chi}_{2^j}(v+u^\alpha)\|_{Y(\Rm)}^4\right)\\ 
 & \qquad +  \|\chi_{2^k} (F(x, v+u^\alpha)-F(x,v)-F(x,u^\alpha)) \|_{L^1 L^2}\\
 & \lesssim_\gamma \sum_{j\geq k} \|\chi_{2^j} w\|_{Y(\Rm)} \left(b_j^4 + |\alpha|^4 2^{-2j}\right) +  \|\chi_{2^k} (F(x, v+u^\alpha)-F(x,v)-F(x,u^\alpha)) \|_{L^1 L^2}
 \end{align*}
Here we use the definition $b_j = \|\tilde{\chi}_{2^k} v\|_{Y(\Rm)}$ and the upper bounds given in \eqref{upper bound on Y u new}. We then insert \eqref{equation 1 for uvw} and obtain
 \begin{align*}
 \|G\|_{L^2(\{s: |s|>2^k\})} \lesssim_\gamma \sum_{j\geq k} y_j \left(b_j^4 + |\alpha|^4 2^{-2j}\right) +  \|\chi_{2^k} (F(x, v+u^\alpha)-F(x,v)-F(x,u^\alpha)) \|_{L^1 L^2}. 
\end{align*}
The second term in the right hand side can be dominated by 
\begin{align}
  \|\chi_{2^k} & (F(x, v+u^\alpha)-F(x,v)-F(x,u^\alpha)) \|_{L^1 L^2} \nonumber\\
  & \lesssim_\gamma \sum_{j \geq k} \|\tilde{\chi}_{2^j} v\|_{Y(\Rm)} \|\tilde{\chi}_{2^j} u^\alpha\|_{Y(\Rm)} \left(\|\tilde{\chi}_{2^j} v\|_{Y(\Rm)}^3 + \|\tilde{\chi}_{2^j} u^\alpha\|_{Y(\Rm)}^3\right) \nonumber \\
  & \lesssim_\gamma \sum_{j \geq k} \left( |\alpha| 2^{-j/2}b_j^4 + |\alpha|^4 2^{-2j} b_j\right). \label{difference F w u alpha}
\end{align}
Thus we have 
\begin{align*}
 \|G\|_{L^2(\{s: |s|>2^k\})} & \lesssim_\gamma \sum_{j\geq k} y_j \left(b_j^4 + |\alpha|^4 2^{-2j}\right)  + \sum_{j \geq k} \left( |\alpha| 2^{-j/2}b_j^4 + |\alpha|^4 2^{-2j} b_j\right)\\
 & \lesssim_\gamma |\alpha|^4 2^{-2k} y_k  + \sum_{j \geq k} \left( |\alpha| 2^{-j/2}b_j^4 + |\alpha|^4 2^{-2j} b_j\right).
\end{align*}
Here we recall the upper bound $y_k \lesssim_1 |\alpha| 2^{-k}$ and the fact that $\{y_k\}$ is a decreasing sequence.  We then recall the characteristic numbers and obtain $\int_\Rm G(s) ds = 0$. Thus 
\begin{align*}
  \left|\int_{-2^k}^{2^k} G(s) ds\right| & = \left|\int_{|s|>2^k} G(s) ds\right| \leq \sum_{j\geq k} \int_{2^j <|s|<2^{j+1}} |G(s)| ds\\
 & \lesssim_\gamma |\alpha|^4 2^{-3k/2} y_k + \sum_{j \geq k} \left( |\alpha| b_j^4 + |\alpha|^4 2^{-3j/2} b_j\right).
\end{align*}
Therefore we have 
\begin{align*}
 y_k & = \|\chi_{2^k} w_L\|_{Y(\Rm)} \lesssim_1 2^{-k/2} \left|\int_{-2^k}^{2^k} G(s) ds\right| + \|G\|_{L^2(\{s: |s|>2^k\})}\\
 & \lesssim_\gamma |\alpha|^4 2^{-2k} y_k + 2^{-k/2} \sum_{j \geq k} \left( |\alpha| b_j^4 + |\alpha|^4 2^{-3j/2} b_j\right)\\
 & \lesssim_\gamma c^{-1} y_k + 2^{-k/2} \sum_{j \geq k} \left( |\alpha| b_j^4 + |\alpha|^4 2^{-3j/2} b_j\right).
\end{align*}
Please note that our assumption on $N$ implies $2^k \geq 2^N \gtrsim_1 c^{1/2} |\alpha|^2$. Therefore we have (Again, if necessary we slightly enlarge the value of $c$)
\[
 y_k \lesssim_\gamma 2^{-k/2} \sum_{j \geq k}(|\alpha| b_j^4 + |\alpha|^4 2^{-3j/2} b_j), \qquad k\geq N.
\]
Finally we recall \eqref{equation 1 for uvw}, \eqref{difference F w u alpha}, $2^N \simeq_\gamma |\alpha|^2$  and obtain 
\begin{align*}
  \|\chi_{2^N} w\|_{Y(\Rm)} & \lesssim_1 \|\chi_{2^N} w_L\|_{Y(\Rm)} + \|\chi_{2^N} (F(x, v+u^\alpha)-F(x,v)-F(x,u^\alpha)) \|_{L^1 L^2}  \\
  & \lesssim_\gamma y_N + \sum_{j \geq N} \left( |\alpha| 2^{-j/2}b_j^4 + |\alpha|^4 2^{-2j} b_j\right)\\
 & \lesssim_\gamma 2^{-N/2} \sum_{j \geq N}(|\alpha| b_j^4 + |\alpha|^4 2^{-3j/2} b_j) \\
 & \lesssim_\gamma \sum_{j\geq N} b_j^4 + 2^{-N/2} |\alpha|^4 \left(\sum_{j\geq N} 2^{-3j}\right)^{1/2}\left(\sum_{j \geq N} b_j^2\right)^{1/2} \lesssim_\gamma \delta.
\end{align*}
We also recall the upper bound of $\int_{-2^k}^{2^k} G(s) ds$ given above and obtain
\begin{align*}
 \left|\int_{-2^N}^{2^N} G(s) ds\right| & \lesssim_\gamma |\alpha|^4 2^{-3N/2} y_N + \sum_{j \geq N} \left( |\alpha| b_j^4 + |\alpha|^4 2^{-3j/2} b_j\right) \\
 & \lesssim_\gamma \sum_{j \geq N} \left( |\alpha| b_j^4 + |\alpha|^4 2^{-3j/2} b_j\right) \lesssim_\gamma 2^{N/2} \delta. 
\end{align*}
\paragraph{Step 2} We then show that the inequality $\|\chi_{R} w\| \lesssim_{\gamma, c_1, M} \delta$ holds by an induction. We first fix a positive constant $\eta = \eta(\gamma) \ll 1$ so that $C_1 \gamma \eta^4 < 1/256$. Here $C_1\geq 1$ is the constant in the Strichartz estimate \eqref{Strichartz estimates}. We then split the region $\Omega_{R,2^N}$ into $N_1 = N_1(\gamma, c_1, M)$ pieces $\Phi_{k}$ accordingly, with 
\[
 \Phi_k = \Omega_{R_{k+1}, R_k}, \qquad R = R_{N_1+1} < R_{N_1} < \cdots < R_3 < R_2 < R_1 = 2^N,
\]
so that $\|\chi_{R_{k+1}, R_k} u^\alpha\|_{Y(\Rm)} < \eta$ and $R_{k+1}/R_k \geq 1/2$. For convenience, we also define $\Phi_0 = \Omega_{2^N}$. In this step we apply an argument of induction to show that the following inequalities hold for each $k=1,2,\cdots, N_1$
\begin{align} \label{step 2 induction hypothesis} 
 &\|\chi_{R_{k+1}} w\|_{Y(\Rm)} \lesssim_{\gamma, k} \delta& &R_{k+1}^{-1/2} \left|\int_{-R_{k+1}}^{R_{k+1}} G(s) ds\right| \lesssim_{\gamma, k} \delta&
\end{align}
as long as $\delta$ is sufficiently small. The case of $k=0$ has been verified in Step 1. Now we assume the inequalities above hold for $k-1$, then show that they also hold for $k$. First of all, we claim that $\|\chi_{R_{k+1}, R_k} u\|_{Y(\Rm)} \leq 2\eta$, as long as $\delta$ is sufficiently small. Otherwise we might substitute $R_{k+1}$ by $R' \in (R_{k+1}, R_k)$ so that $\|\chi_{R', R_k} u\|_{Y(\Rm)} = 2\eta$ and proceed as normal. Finally we would obtain 
\begin{align*}
 \|\chi_{R', R_k} u\|_{Y(\Rm)} & \leq \|\chi_{R', R_k} u^\alpha\|_{Y(\Rm)} + \|\chi_{R', R_k} w\|_{Y(\Rm)} + \|\chi_{R',R_k} v\|_{Y(\Rm)} \\
 & \leq \eta + C(\gamma, k) \delta + \delta
\end{align*}
This contracts with our assumption $\|\chi_{R', R_k} u\|_{Y(\Rm)} = 2\eta$ when $\delta$ is sufficiently small. We start by applying the Strichartz estimates
\begin{align*}
 \|\chi_{R_{k+1}, R_k} w\|_{Y(\Rm)} & \leq \|\chi_{R_{k+1}, R_k} w_L\|_{Y(\Rm)} + C_1 \|\chi_{R_{k+1}} (F(x,u)-F(x,v)-F(x,u^\alpha))\|_{L^1 L^2}\\
 & \leq \|\chi_{R_{k+1}, R_k} w_L\|_{Y(\Rm)} + J_1 + J_2 + J_3.
\end{align*}
The definitions of $J$'s and their upper bounds are given below (we always choose $\delta < \eta$)
\begin{align*}
 J_1 & = C_1\|\chi_{R_{k+1}, R_k} (F(x,u)-F(x,v+u^\alpha))\|_{L^1 L^2} \\
 & \leq C_1 \gamma \|\chi_{R_{k+1}, R_k} w\|_{Y(\Rm)} \left(\|\chi_{R_{k+1}, R_k} u\|_{Y(\Rm)}^4 + \|\chi_{R_{k+1}, R_k} (v+u^\alpha)\|_{Y(\Rm)}^4 \right)\\
 & \leq 32 C_1\gamma \eta^4 \|\chi_{R_{k+1}, R_k} w\|_{Y(\Rm)}. 
\end{align*}
We then use induction hypothesis 
\begin{align*}
 J_2 & = C_1 \|\chi_{R_k} (F(x,u)-F(x,v+u^\alpha))\|_{L^1 L^2} \\
 & \leq C_1 \gamma \|\chi_{R_k} w\|_{Y(\Rm)} \left(\|\chi_{R_k} u\|_{Y(\Rm)}^4 + \|\chi_{R_k} (v+u^\alpha)\|_{Y(\Rm)}^4 \right)\\
 & \lesssim_{\gamma, k} \delta.
\end{align*}
We also have 
\begin{align*}
 J_3 & = C_1 \|\chi_{R_{k+1}} (F(x,v+u^\alpha)-F(x,v)-F(x,u^\alpha))\|_{L^1 L^2}\\
 & \lesssim_\gamma \|\chi_{R_{k+1}} v\|_{Y(\Rm)} \|\chi_{R_{k+1}} u^\alpha\|_{Y(\Rm)} \left(\|\chi_{R_{k+1}} u^\alpha\|_{Y(\Rm)}^3 + \|\chi_{R_{k+1}} v\|_{Y(\Rm)}^3 \right)\\
 & \lesssim_{\gamma, k} \delta.
\end{align*}
As a result, we recall $32 C_1\gamma \eta^4 < 1/2$ and obtain
\begin{equation} \label{step 2 wL to w}
 \|\chi_{R_{k+1}, R_k} w\|_{Y(\Rm)} \leq 2 \|\chi_{R_{k+1}, R_k} w_L\|_{Y(\Rm)} + C(\gamma, k) \delta. 
\end{equation}
The notation $C(\gamma, k)$ represents a constant depending on $\gamma, k$ but nothing else. We slightly abuse the notation so that it may represent different constants at different places. We then observe that the nonlinear radiation profile of $w$ is zero for $|s|>R$ in both two time direction. Thus we may apply Lemma \ref{scatter profile of nonlinear solution}, utilize the upper bounds of $J_1, J_2, J_3$ and obtain 
\begin{align*}
 2\sqrt{2\pi} \|G\|_{L^2(\{s:|s|>R_{k+1}\})} & \leq \|\chi_{R_{k+1}} (F(x,u)-F(x,v)-F(x,u^\alpha))\|_{L^1 L^2} \\
 & \leq 32 \gamma \eta^4 \|\chi_{R_{k+1}, R_k} w\|_{Y(\Rm)} + C(\gamma, k)\delta \\
 & \leq 64 \gamma \eta^4 \|\chi_{R_{k+1}, R_k} w_L\|_{Y(\Rm)} + C(\gamma, k)\delta\\
 & \leq (4C_1)^{-1} \|\chi_{R_{k+1}, R_k} w_L\|_{Y(\Rm)} + C(\gamma, k)\delta.
\end{align*}
Here we recall our assumption $C_1 \gamma \eta^4 < 1/256$. We may combine this with the induction hypothesis and obtain 
\begin{align}
  R_{k+1}^{-1/2} \left|\int_{-R_{k+1}}^{R_{k+1}} G(s) ds\right| & \leq R_{k+1}^{-1/2} \left|\int_{R_{k+1} < |s|<R_{k}} G(s) ds\right| +R_{k+1}^{-1/2} \left|\int_{-R_{k}}^{R_k} G(s) ds\right| \nonumber\\
 & \leq (8C_1)^{-1}\pi^{-1/2} \|\chi_{R_{k+1}, R_k} w_L\|_{Y(\Rm)} + C(\gamma,k) \delta. \label{step 2 wL to int G}
\end{align}
We consider the linear free wave $w'_L$ whose radiation profile $G'$ is given by 
\[
 G'(s) = \left\{\begin{array}{ll} \frac{1}{2R_{k+1}} \int_{-R_{k+1}}^{R_{k+1}} G(s') ds', & |s|<R_{k+1}; \\ G(s), & |s|>R_{k+1}. \end{array} \right.
\]
It is clear that $\chi_{R_{k+1}} w_L = \chi_{R_{k+1}} w'_L$ holds by the explicit formula of linear free wave in term of radiation profiles. We also have 
\begin{align*}
 \|G'\|_{L^2(\Rm)} & \leq (2R_{k+1})^{-1/2} \left|\int_{-R_{k+1}}^{R_{k+1}} G(s) ds\right| + \|G\|_{L^2(\{s:|s|>R_{k+1}\})}\\
 & \leq (4C_1)^{-1} (2\pi)^{-1/2} \|\chi_{R_{k+1}, R_k} w_L\|_{Y(\Rm)} + C(\gamma,k) \delta. 
\end{align*}
Thus by the Strichartz estimates we have 
\begin{align*}
 \|\chi_{R_{k+1}, R_k} w_L\|_{Y(\Rm)} \leq C_1 2\sqrt{2\pi} \|G'\|_{L^2(\Rm)}
  \leq (1/2) \|\chi_{R_{k+1}, R_k} w_L\|_{Y(\Rm)} + C(\gamma,k) \delta. 
\end{align*}
Therefore we have 
\[
  \|\chi_{R_{k+1}, R_k} w_L\|_{Y(\Rm)} \lesssim_{\gamma, k} \delta. 
\]
Finally we insert this upper bound into \eqref{step 2 wL to w} and \eqref{step 2 wL to int G} to verify \eqref{step 2 induction hypothesis} thus finish the proof of this step. 
\paragraph{Step 3} Finally we collect all the estimates given in Step 1 and 2. Considering the case $k= N_1$ in step 2, we have 
\begin{align*} 
 &\|\chi_{R} w\|_{Y(\Rm)} \lesssim_{\gamma, c_1, M} \delta;& &R^{-1/2} \left|\int_{-R}^{R} G(s) ds\right| \lesssim_{\gamma, c_1, M} \delta;& &\|G\|_{L^2(\{s:|s|>R\})} \lesssim_{\gamma, c_1, M} \delta. &
\end{align*}
This gives
\[
 \|\chi_R(u-u^\alpha)\|_{Y(\Rm)} \leq \|\chi_R w\|_{Y(\Rm)} + \|\chi_R v\|_{Y(\Rm)} \lesssim_{\gamma, c_1, M} \delta;
\]
and 
\[
 \|(w_0, w_1)\|_{\dot{H}^1 \times L^2(\{x: |x|>R\})} \lesssim_1 R^{-1/2}  \left|\int_{-R}^{R} G(s) ds\right| + \|G\|_{L^2(\{s:|s|>R\})} \lesssim_{\gamma, c_1, M} \delta;
\]
thus finishes the proof. 
\end{proof}

\noindent The remaining part of this section is devoted to the proof of Theorem \ref{main global}.

\begin{proof}[Proof of Theorem \ref{main global}]
 The proof consists of two parts. We first show that the solutions satisfying the conditions must be globally defined for all $t > 0$ and then show that these solutions scatter in the positive time direction. 
  
 \paragraph{Part 1} We define 
 \[
  \mathcal{X} = \left\{A_1 \in \Rm^+: \begin{array}{c} \hbox{any solution to (CP1) with a maximal lifespan}\, (-T_+,T_+)\, \hbox{and} \\
  \displaystyle \limsup_{t\rightarrow T_+}  \|(u(\cdot,t), u_t(\cdot,t))\|_{\dot{H}^1 \times L^2} < A_1\, \hbox{must satisfy}\, T_+ = +\infty \end{array} \right\}
 \]
By scattering theory with small initial data, the set $\mathcal{X}$ is nonempty. The way we define $\mathcal{X}$ implies that $\mathcal{X} = (0,A_\mathcal{X}]$ for some $A_\mathcal{X} \in \Rm^+$ or $\mathcal{X} = \Rm^+$. The first part of our proof is to prove that $A \in \mathcal{X}$. If this were false, we would have that $A_\mathcal{X} < A$. By the definition of $\mathcal{X}$, given any $\eps > 0$, there exists a solution $u$ to (CP1) with a maximal lifespan $(-T_-,T_+)$, $T_+ \in \Rm^+$ and 
\[
 \limsup_{t\rightarrow T_+} \|(u(\cdot,t), u_t(\cdot,t))\|_{\dot{H}^1 \times L^2} < \sqrt{A_\mathcal{X}^2 + \eps^2}. 
\]
We next show that this can never happen when $\eps$ is sufficiently small. We can always assume 
\begin{equation} \label{boundedness at T plus}
 \sup_{t\in [0,T_+)} \|(u(\cdot,t), u_t(\cdot,t))\|_{\dot{H}^1 \times L^2} < \sqrt{A_\mathcal{X}^2 + \eps^2},
\end{equation}
by a time translation if necessary. This gives a universal upper bound
\[
 u(r,t) \lesssim_1 r^{-1/2} \|u(\cdot,t)\|_{\dot{H}^1(\Rm^3)} \lesssim_1 r^{-1/2} \sqrt{A_\mathcal{X}^2 + \eps^2}, \qquad (r,t) \in \Rm^+ \times [0,T_+).
\]
This implies that given $t_1 \in (0,T_+)$, if we define $\Phi_{1,t} = \{(x,t): t\in [t_1, T_+), \, |x|>t-t_1\}$, then
\[ 
 \|\chi_{\Phi_{1,t_1}} u\|_{Y([t_1, T_+))} < +\infty. 
\]
Therefore $T_+$ is still contained in the maximal lifespan of the exterior solution to the following equation, whose initial data are given at time $t_1$
\[
 \left\{\begin{array}{ll} \partial_t^2 u - \Delta u = F(x,u), & |x|>|t-t_1|; \\ (u,u_t)|_{t=t_1} = (u(\cdot,t_1), u_t(\cdot,t_1)) & \end{array} \right.
\]
Please note that although the exterior solution here can be defined for some time $t>T_+$, we still use the same notation $u$ since this exterior solution always coincides with the original solution in the overlap region of their domains. Next we define data $(v_0,v_1)$ by 
\[
 (v_0(x), v_1(x)) = (u(x,T_+), u_t(x,T_+)), \qquad |x|>T_+ - t_1.
\] 
The right hand side is the data of the exterior solution $u$ defined above at the time $T_+$. Given different times $t_1, t_2 \in (0,T_+)$, finite speed of propagation shows that the exterior solutions defined above with initial times $t_1, t_2$ share the same data at time $T_+$ in the region $\{x: |x|>T_+ - \min\{t_1,t_2\}\}$. Therefore we may define $(v_0,v_1)$ in the region $\Rm^3 \setminus \{0\}$ without any conflict. Our assumption \eqref{boundedness at T plus} and the continuity of (global extension of) exterior solutions in $\dot{H}^1 \times L^2(\Rm^3)$ then implies that 
\[
 \int_{\Rm^3} (|\nabla v_0(x)|^2 + |v_1(x)|^2) dx < A_\mathcal{X}^2  + \eps^2.  
\]
Therefore we have $(v_0,v_1) \in \dot{H}^1 \times L^2(\Rm^3)$. We then consider the solution $v$ to (CP1) with initial data $(v(\cdot,T_+), v_t(\cdot,T_+)) = (v_0,v_1)$. We claim that for $t$ in the maximal lifespans of both solutions $u$ and $v$, we have 
\[
 v(x,t) = u(x,t), \qquad |x|>|T_+-t|. 
\]
In fact, a combination of the way we define $v_0,v_1$ and the finite speed of propagation implies that the identity holds as long as $|x|>|T_+ - t_1|+|T_+-t|$ for any $t_1\in [0,T_+)$. We then obtain the coincidence of $u$ and $v$ in the region $|x| > |T_+ - t|$ by letting $t_1 \rightarrow T_+$. By continuity of data with respect to the time, we obtain: given any small positive constant $\eps_1 > 0$, there exists $r_1 > 0$ and $\delta_1 > 0$ so that 
\[
 \|\nabla_{x,t} v(\cdot,t) \|_{L^2(\{x: |x|<r_1\})} + \|v(\cdot,t)\|_{L^6(\{x: |x|<r_1\})}< \eps_1, \qquad t\in (T_+-\delta_1, T_+].
\]
Without loss of generality, we may choose $\delta_1 < r_1/100$. This then gives 
\begin{equation} \label{empty barrier}
 \|\nabla_{x,t} u(\cdot,t) \|_{L^2(\{x: |T_+ -t|<|x|<r_1\})} + \|u(\cdot,t)\|_{L^6(\{x: |T_+-t| < |x|<r_1\})} < \eps_1, \qquad t\in (T_+-\delta_1, T_+].
\end{equation}
We fix a smooth cut-off function $\varphi : \Rm \rightarrow [0,1]$ so that $\varphi(r) = 1$ if $r\leq 3/4$ and $\varphi(r) = 0$ if $r\geq  1$. Given $t_1\in (T_+-\delta_1, T_+)$, we consider two new pairs of initial data
\begin{align*}
 (u_0^{t_1}(x), u_1^{t_1}(x)) &= \varphi(|x|/r_1) (u(x,t_1), u_t(x,t_1));\\
 (\tilde{u}_0^{t_1}(x), \tilde{u}_1^{t_1}(x)) & = [1-\varphi(2|x|/r_1)]\varphi(|x|/r_1) (u(x,t_1), u_t(x,t_1));
\end{align*} 
and consider the solutions $u^{t_1}, \tilde{u}^{t_1}$ to (CP1) with these initial data at time $t_1$. Since we have 
\[
 (u_0^{t_1}(x), u_1^{t_1}(x)) = (u(x,t_1), u_t(x,t_1)), \qquad |x|<3r_1/4,
\]
We always have $u(x,t) = u^{t_1}(x,t)$ in the region $\{x: |x|<3r_1/4 - |t-t_1|\}$ as long as $t\in (t_1,T_+)$ is still contained in the maximal lifespan of $u^{t_1}$.   This implies that $u^{t_1}$ must blow up at a time $T_{t_1,+} \leq T_+$ because $u$ blows up at time $T_+$ and at the origin. We also have for $t\in [t_1, T_{t_1,+})$, 
\begin{align}
 \int_{|x|<3r_1/4-|t-t_1|} |\nabla_{t,x} u^{t_1}(x,t)|^2 dx & = \int_{|x|<3r_1/4-|t-t_1|} |\nabla_{t,x} u(x,t)|^2 dx\nonumber \\
 & \leq \int_{|x|<|T_+-t|} |\nabla_{t,x} u(x,t)|^2 dx + \eps_1^2. \label{center part estimate}
\end{align}
Here we use the upper bound \eqref{empty barrier}. In addition, a straight-forward calculation shows that 
\[
 \|(\tilde{u}_0^{t_1}(x), \tilde{u}_1^{t_1}(x))\|_{\dot{H}^1 \times L^2(\Rm^3)} \lesssim_1 \eps_1. 
\]
When $\eps_1$ is sufficiently small, the solutions $\tilde{u}^{t_1}$ are all global scattering solutions so that 
\[
  \sup_{t\in \Rm} \|(\tilde{u}^{t_1}(x,t), \tilde{u}_t^{t_1}(x,t))\|_{\dot{H}^1 \times L^2 (\Rm^3)} \lesssim_1 \eps_1. 
\]
We then observe that $(u_0^{t_1}(x), u_1^{t_1}(x)) = (\tilde{u}_0^{t_1}(x), \tilde{u}_1^{t_1}(x))$ for $|x|>r_1/2$, thus $u^{t_1} (x,t) = \tilde{u}^{t_1}(x,t)$ in the exterior region $|x|>|t-t_1|+r_1/2$ as long as $u^{t_1}$ is still well-defined at time $t$. Therefore we have 
\[
 \int_{|x|>r_1/2 +|t-t_1|} |\nabla_{t,x} u^{t_1}(x,t)|^2 dx \lesssim_1 \eps_1^2, \qquad t\in [t_1, T_{t_1,+}). 
\]
We may combine this with \eqref{center part estimate} and obtain
\[
 \|u^{t_1} (\cdot,t)\|_{\dot{H}^1\times L^2(\Rm^3)}^2 \leq \int_{|x|<|T_+-t|} |\nabla_{t,x} u(x,t)|^2 dx + C_1 \eps_1^2, \qquad t\in [t_1, T_{t_1,+}).
\]
Here $C_1>0$ is an absolute constant. We use the definition of $A_\mathcal{X}$, consider the blow-up behaviour of $u^{t_1}$ and obtain that given any $t_1 \in (T_+-\delta_1,T_+)$, there exists a time $t' \in [t_1, T_{t_1,+})$ so that 
\[
 \int_{|x|<|T_+-t'|} |\nabla_{t,x} u(x,t')|^2 dx \geq A_\mathcal{X}^2 - (C_1+1) \eps_1^2
\]
Therefore we have 
\[
 \limsup_{t\rightarrow T_+} \int_{|x|<|T_+-t|} |\nabla_{t,x} u(x,t)|^2 dx\geq A_\mathcal{X}^2 - (C_1+1) \eps_1^2.
\]
We make $\eps_1\rightarrow 0^+$ and obtain the energy concentration 
\[
 \limsup_{t\rightarrow T_+} \int_{|x|<|T_+-t|} |\nabla_{t,x} u(x,t)|^2 dx\geq A_\mathcal{X}^2.
\]
We then recall the uniform boundedness \eqref{boundedness at T plus} and conclude that 
\[
  \liminf_{t\rightarrow T_+} \int_{|x|>|T_+-t|} |\nabla_{t,x} u(x,t)|^2 dx \leq \eps^2.
\]
It immediately follows that $\|(v_0,v_1)\|_{\dot{H}^1 \times L^2(\Rm^3)} \leq \eps$. Thus by small data scattering $v$ is a global scattering solution. Let $G^+$ be the nonlinear radiation profile of $v$. Clearly $\|G^+\|_{L^2(\Rm)} \lesssim_1 \eps$. We then recall the coincidence of $u$ and $v$ in the region $\{(x,t): |x|>|T_+-t|, t\in (-T_-,T_+)\}$ and obtain
\[
  \int_{|x|>|T_+-t|} |\nabla_{t,x} u(x,t)|^2 dx \leq 2\eps^2, \qquad \forall t\in (-T_-,T_+).
\]
In addition, the inequality above (at time $0$) enables us to choose a small constant $\delta'$ so that 
\[
 \int_{|x|> T_+ -\delta'} |\nabla_{t,x} u(x,0)|^2 dx \leq 3\eps^2
\]
By small data scattering, the exterior solution (we still call it $u$ by the coincidence of solutions in the overlapping region of their domains) to (CP1) in the exterior region $\Omega_{T_+ - \delta'}$ with initial data $(u(\cdot,0), u_t(\cdot,0))$ must be defined for all time $t$ and satisfy
\[
 \int_{|x|>T_++|t|-\delta'} |\nabla u(x,t)|^2 dx \leq 6\eps^2, \qquad \forall t\in \Rm;
\]
and 
\[
 \|\chi_{T_+-\delta'} u\|_{Y(\Rm)} \lesssim_1 \eps.
\]
Let $G^-$ be the nonlinear radiation profile of this exterior solution $u$ in the negative time direction. We have $\|G^-\|_{L^2([T_+-\delta', +\infty))} \lesssim_1 \eps$. Now we fix a time $t_1 \in (T_+ - \delta' , T_+)$ and consider the exterior solution $w$ to (CP1) in the region $\Omega_0$ with initial data (at time zero) $(u(\cdot,t_1), u_t(\cdot,t_1))$. Collecting the information of $v$ and $u$ (including the extension of $u$ defined above), we obtain that $w$ is defined for all $t$ and (please see figure \ref{figure exsol} for an illustration of $w$)
\[
 w(x,t) = \left\{\begin{array}{ll} v(x, t+t_1), & t\geq T_+ - t_1, \, |x|>|t|;\\ u(x, t+t_1), & t< T_+ - t_1,\, |x|>|t|. \end{array}\right.
\]
 \begin{figure}[h]
 \centering
 \includegraphics[scale=1.25]{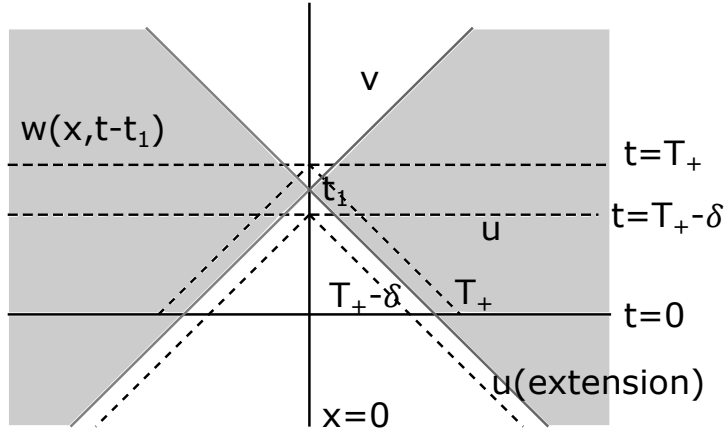}
 \caption{the exterior solution $w$} \label{figure exsol}
\end{figure}
We then recall the non-linear radiation profiles of $u$ and $v$, then conclude that $w$ is asymptotically equivalent to a linear free wave $w_L$, which is defined by its radiation profiles $G_L^\pm$:
\begin{align*}
 &G_L^+(s) = G^+(s-t_1), \; \forall s>0;& &G_L^-(s) = G^-(s+t_1),\; \forall s>0.&
\end{align*}
Since we have
\[
 \|G_L^\pm\|_{L^2(\Rm)}^2 = \|G_L^+\|_{L^2(\Rm^+)}^2 + \|G_L^-\|_{L^2(\Rm^+)}^2 = \|G^+\|_{L^2((-t_1,+\infty))}^2 + \|G^-\|_{L^2((t_1,+\infty))}^2 \lesssim_1 \eps^2,
\]
 Corollary \ref{nonlinear asymptotic with small norm profiles} and Proposition \ref{L2 small nonlinear} then give a global scattering solution $v^{t_1}$ so that $v^{t_1}$ is asymptoticly equivalent to $w$ and satisfies 
\[
 \sup_{t\in \Rm} \|\nabla_{x,t} v^{t_1}(\cdot, t)\|_{L^2(\Rm)}^2 + \sum_{k=-\infty}^\infty \|\chi_{2^k} v^{t_1}\|_{Y(\Rm)}^2 \lesssim_1 \eps^2. 
\]
We claim that if $\eps$ is sufficiently small, then $w = v^{t_1}$ in $\Omega_0$ thus the data of $u$ at $t_1$ are exactly the same as the data of $v^{t_1}$, which is very small. This is a contradiction thus finishes the proof of the first part. The remaining work is to show $w= v^{t_1}$ if $\eps$ is sufficiently small. Since they are asymptotically equivalent to the same linear free wave $w_L$, it suffices to show the characteristic number of $w$ respect to $v^{t_1}$ is zero. If the characteristic number were $\alpha \neq 0$, we could apply Lemma \ref{partition of u with u alpha} on $v^{t_1}, u^\alpha$ in the region $\Omega_{r_\alpha}$, as long as $\eps$ is sufficiently small and obtain 
\[
 \|u(\cdot,t_1) - u^\alpha - v^{t_1}(\cdot,0)\|_{\dot{H}^1(\{x:|x|>r_\alpha\})} \lesssim_{\gamma, c_1, A} \eps.
\]
This implies ($C(\gamma, c_1, A)$ is a constant solely determined by the parameters $\gamma, c_1, A$)
\[
 \sqrt{A_\mathcal{X}^2 +\eps^2}\geq \|u(\cdot,t_1)\|_{\dot{H}^1(\{x: |x|>r_\alpha\})} \geq A - C(\gamma, c_1, A) \eps.
\]
Finally we let $\eps\rightarrow 0^+$. This immediately gives a contradiction with our assumption $A_\mathcal{X} < A$.

\paragraph{Part 2} Now we have obtained $T_+ = +\infty$. We first choose a large radius $R'>0$, so that the initial data $(u_0,u_1)$ satisfy $\|(u_0,u_1)\|_{\dot{H}^1 \times L^2(\{x: |x|>R'\})} \ll1$. It immediately follows the small data scattering theory that the exterior solution $u$ to (CP1) in the region $\Omega_{R'}$ with initial data $(u_0,u_1)$ is globally defined for all $t\in \Rm$ and scatters in the exterior region with
\begin{equation} \label{exterior Y norm step 2}
 \|\chi_{R'} u\|_{Y(\Rm)} < +\infty. 
\end{equation}
 We still call it $u$, because of the coincidence of this exterior solution and the original solution in the overlapping region of their domains. Our assumption of uniformly bounded $\dot{H}^1$ norm guarantees that 
 \[
  |u(r,t)| \lesssim_1 r^{-1/2} \|u(\cdot,t)\|_{\dot{H}^1(\Rm^3)} \lesssim r^{-1/2}, \qquad \forall r,t>0.
 \]
 Combining this estimate with $\|u\|_{Y([0,T])} < +\infty$ and \eqref{exterior Y norm step 2}, we have 
\[
 \|\chi_{\{(x,t): |x| > |t-\beta|\}} u\|_{Y(\Rm)} < +\infty, \qquad \forall \beta \geq R'.
\]
We may view $u$ as an exterior solution (time-translated) defined in the region $\{(x,t): |x|>|t-\beta|\}$ for any time $\beta \geq R'$ and then apply Lemma \ref{scatter profile of nonlinear solution} to conclude that there exist nonlinear profiles $G^\pm$ satisfying  
\begin{align*}
 &G^+\in L^2([-\beta, +\infty)), \; \forall \beta \geq R';& &G^- \in L^2([R', +\infty));& 
\end{align*}
so that given any $\beta\geq R'$, we have
 \begin{align*}
  \lim_{t\rightarrow +\infty} \int_{t-\beta}^\infty \left(\left|G^+(r-t) - r u_t (r, t)\right|^2 + \left|G^+(r-t) + r u_r (r, t)\right|^2\right) dr & = 0; \\
  \lim_{t\rightarrow -\infty} \int_{-t+R'}^\infty \left(\left|G^-(r+t) - r u_t(r,t)\right|^2 +  \left|G^-(r+t) - r u_r(r,t)\right|^2\right)dr & = 0.
 \end{align*}
 We then use the assumption  
 \[
  A_1 \doteq  \limsup_{t\rightarrow +\infty} \|(u(\cdot,t), u_t(\cdot,t))\|_{\dot{H}^1 \times L^2} < A, 
 \]
 and obtain 
 \[
  \|G^+\|_{L^2([-\beta,+\infty))} \lesssim_1 A_1, \, \forall \beta \geq R' \qquad \Rightarrow \qquad  \|G^+\|_{L^2(\Rm)} \lesssim_1 A_1 < +\infty.
 \]
 Next we consider the linear free wave $w_L^{t_1}$ asymptotically equivalent to the time-translated solution $u(x, t_1 + t)$, given a very large time $t_1>R'$. We use the radiation profiles given above and obtain that the radiation profiles $G_{t_1}^\pm$ of $w_L^{t_1}$ satisfy (recall $G_{t_1}^+(-s) = - G_{t_1}^-(s)$)
 \begin{align*}
  &G_{t_1}^+ (s) = G^+(s-t_1),\; \forall s>0; & &G_{t_1}^-(s) = G^-(s+t_1),\; \forall s>0&
 \end{align*}
 We claim that the following limits hold for $w_L^{t_1}$:
 \begin{itemize}
  \item[(i)] $\displaystyle \lim_{t_1\rightarrow +\infty} \sum_{k=-\infty}^\infty \|\tilde{\chi}_{2^k} w_L^{t_1}\|_{Y(\Rm)}^2 = 0$, thus $\displaystyle \lim_{t_1\rightarrow +\infty} \|\chi_0 w_L^{t_1}\|_{Y(\Rm)} = 0$;
  \item[(ii)] $\displaystyle \lim_{R_1 \rightarrow +\infty} \left(\limsup_{t_1\rightarrow +\infty} \int_{||x|-t_1|>R_1} |\nabla_{t,x} w_L^{t_1} (x,0)|^2 dx\right) = 0$;
  \item[(iii)] $\displaystyle \lim_{t_1 \rightarrow +\infty} \|w_L^{t_1}\|_{Y(\Rm^+)} = 0$.
 \end{itemize}
 It suffices to show that these limits hold if $G^+$ is compactly supported and $G^- = 0$. The general case immediately follows because we may split $G^+$ into two parts: a compactly-supported part and another part with small $L^2$ norm, observe the fact $\|G_{t_1}^-\|_{L^2(\Rm^+)}\rightarrow 0$ as $t_1\rightarrow +\infty$, and then apply the Strichartz estimates and Proposition \ref{l2 linear}. Now let us assume that $G^+$ is supported in $[-R,R]$ and $G^- = 0$. When $t_1$ is sufficiently large we have 
 \begin{align*}
  &\hbox{Supp}\, G_{t_1}^- \subset [-R-t_1, R-t_1];& &\|G_{t_1}^-\|_{L^2(\Rm)} = \|G^+\|_{L^2(\Rm)}.&
 \end{align*}
 We then apply Remark \ref{channel by small support data} and Lemma \ref{wide channel by localized radiation}:
 \[
  \|\tilde{\chi}_{2^k} w_L^{t_1}\|_{Y(\Rm)} \lesssim_1 \left\{\begin{array}{ll} (R/t_1)^{1/2} (2^k/t_1)^{1/10} \|G^+\|_{L^2}, & 2^k < t_1\\
  (R/2^k)^{1/2}\|G^+\|_{L^2}, & 2^k \geq t_1 \end{array} \right.
 \]
 This immediately verifies the first limit of (i). The second limit is a direct consequence of the first one: 
  \begin{equation*}
  \|\chi_0 w_L^{t_1}\|_{Y(\Rm)} \leq \left(\sum_{k=-\infty}^\infty \|\chi_{2^k} w_L^{t_1}\|_{Y(\Rm)}^5\right)^{1/5} \leq \left(\sum_{k=-\infty}^\infty \|\chi_{2^k} w_L^{t_1}\|_{Y(\Rm)}^2\right)^{1/2} \rightarrow 0.
 \end{equation*} 
 We next recall the explicit formula 
 \begin{align*}
  &w_L^{t_1}(r,0) = \frac{1}{r} \int_{-r}^r G_{t_1}^- (s) ds;& &\partial_t w_L^{t_1} (r,0) =\frac{G_{t_1}^- (r) - G_{t_1}^-(-r)}{r}&
 \end{align*}
 A straight-forward calculation shows that when $R_1 \geq R$ and $t_1$ is large, we have
 \begin{align*}
  \int_{||x|-t_1|>R_1} |\nabla_{t,x} w_L^{t_1} (x,0)|^2 dx & \lesssim_1 \int_{t_1+R_1}^\infty r^{-2} \left|\int_{-R}^R G^+(s) ds\right|^2 dr \\
  & \lesssim_1 (t_1+R_1)^{-1}\left|\int_{-R}^R G^+(s) ds\right|^2.
 \end{align*}
 This verifies (ii).  Our assumption on $G^+$ and $G^-$ implies $w_L^{t_1} (x,t) = 0$ if $|x|<t$, as long as $t_1$ is sufficiently large. A combination of this fact with (i) immediately gives (iii). Next we utilize (i) and consider the primary asymptotically equivalent solution $w^{t_1}$ of $w_L^{t_1}$, if $t_1$ is sufficiently large. We recall the properties of primary asymptotically equivalent solution given in Subsection \ref{sec: primary}:
\begin{itemize} 
 \item $w^{t_1}$ is defined in $\Omega_0$ and asymptotically equivalent to $w_L^{t_1}$ (or $u(x,t+t_1)$); 
 \item The initial data $(w_0^{t_1}, w_1^{t_1})$ of $w^{t_1}$ satisfy
 \begin{equation} \label{approximation of initial data part 2}
  \left\|(w_0^{t_1}, w_1^{t_1})-(w_{L}^{t_1}(\cdot,0), \partial_t w_{L}^{t_1}(\cdot,0))\right\|_{\dot{H}^1 \times L^2(\Rm^3)} \lesssim_\gamma \|\chi_0 w_L^{t_1}\|_{Y(\Rm)}^5 \rightarrow 0. 
 \end{equation}
 \item The solution $w^{t_1}$ also satisfies 
 \begin{equation} \label{l2 small of w t1}
  \lim_{t_1\rightarrow +\infty} \sum_{k=-\infty}^\infty \|\tilde{\chi}_{2^k} w^{t_1}\|_{Y(\Rm)}^2 = 0.
 \end{equation} 
\end{itemize} 
In addition, we may also combine \eqref{approximation of initial data part 2} with (ii), (iii) above to obtain 
 \begin{equation} \label{concentration of energy initial w}
 \lim_{R_1 \rightarrow +\infty} \left(\limsup_{t_1\rightarrow +\infty} \int_{||x|-t_1|>R_1} \left(|\nabla w_0^{t_1} (x)|^2 + |w_1^{t_1}(x)|^2\right) dx\right) = 0;
 \end{equation}
 and
 \begin{equation} \label{scattering of w positive time}
  \lim_{t_1 \rightarrow +\infty} \|\mathbf{S}_L (w_0^{t_1}, w_1^{t_1})\|_{Y(\Rm^+)} = 0.
 \end{equation}
 Now we let $t_1$ be a sufficiently large number and consider the solution $u(x,t+t_1)$ and $w^{t_1}$, which are both asymptotically equivalent to the same linear free wave $w_L^{t_1}$. Let $\alpha(t_1)$ be the characteristic number of $u(x,t+t_1)$ with respect to $w^{t_1}$. If $\alpha(t_1) = 0$, then we immediately have $(u(\cdot,t_1), u_t(\cdot,t_1)) = (w_0^{t_1}, w_1^{t_1})$. This finishes the proof if $t_1$ is sufficiently large because \eqref{scattering of w positive time} implies that $u$ must scatter in the positive time direction. Next we assume $\alpha(t_1) \neq 0$ for all sufficiently large time $t_1$, recall \eqref{l2 small of w t1}, apply Lemma \ref{partition of u with u alpha} on the weakly non-radiative solutions $u^{\alpha(t_1)}$, the exterior solution $w^{t_1}$ in the region $\Omega_{r_{\alpha(t_1)}}$ and obtain
\begin{equation} \label{asymptotic decomposition}
 \lim_{t_1\rightarrow +\infty} \|u(\cdot,t_1) - w_0^{t_1} - u^{\alpha(t_1)}\|_{\dot{H}^1(\{x: |x|>r_{\alpha(t_1)}\})} = 0
\end{equation}
In addition, Theorem \ref{main non-radiative} implies that there exists a constant $\tau_0 = \tau_0(\gamma)$ so that 
\begin{equation} \label{tail of u alpha}
  \|u^\alpha\|_{\dot{H}^1(\{x: |x|>\tau |\alpha|^2\})} \simeq_1 \tau^{-1/2}, \qquad \forall \tau \geq \tau_0.
\end{equation}
Therefore if we choose a constant $c_2 \gg \max\{A^{-2}, \tau_0\}$, then we must have $r_\alpha \leq c_2 |\alpha|^2$. We claim that $t_1^{-1} |\alpha(t_1)|^2 \rightarrow 0^+$ as $t_1 \rightarrow +\infty$. If this were false, then we would find a sequence of times $t_j \rightarrow +\infty$ and a positive constant $\kappa$, so that $|\alpha(t_j)|^2 \geq \kappa t_j$. We choose a large constant $\tau > \max\{\tau_0, c_2, 2/\kappa\}$ and consider the $\dot{H}^1$ norms of the functions $u(\cdot,t_j), w_0^{t_j}, u^\alpha(t_j)$ in the exterior regions $\{x: |x|>\tau |\alpha(t_j)|^2\}$. By our assumption on $t_j$ we have $\tau |\alpha(t_j)|^2 \geq \tau \kappa t_j > 2t_j$. We then use the asymptotic behaviour of $u$ and obtain 
\[
 \lim_{j \rightarrow +\infty} \|u(\cdot,t_j)\|_{\dot{H}^1 (\{x: |x|>\tau |\alpha(t_j)|^2\})} \leq \lim_{j \rightarrow +\infty} \|u(\cdot,t_j)\|_{\dot{H}^1 (\{x: |x|>2t_j\}} = 0.
\]
Similarly by \eqref{concentration of energy initial w} we also have 
\[
 \lim_{j \rightarrow +\infty} \|w_0^{t_j}\|_{\dot{H}^1 (\{x: |x|>\tau |\alpha(t_j)|^2\})}  \leq \lim_{j \rightarrow +\infty} \|w_0^{t_j}\|_{\dot{H}^1 (\{x: |x|>2t_j\})} = 0.
\]
Our choice of $\tau$ guarantees that $\tau |\alpha(t_j)|^2 > r_{\alpha(t_j)}$. Thus we may combine these limits with \eqref{asymptotic decomposition} and obtain
\[
 \lim_{j \rightarrow +\infty} \|u^{\alpha(t_j)}\|_{\dot{H}^1 (\{x: |x|>\tau |\alpha(t_j)|^2\})} = 0.
\]
This contradicts with \eqref{tail of u alpha}. Therefore the limit $t_1^{-1} |\alpha(t_1)|^2 \rightarrow 0^+$ must hold as $t_1 \rightarrow +\infty$. By \eqref{tail of u alpha} we have 
\[
 \lim_{t_1 \rightarrow +\infty} \|u^{\alpha(t_1)}\|_{\dot{H}^1(\{x: |x|>0.9 t_1\})} \lesssim_1 \lim_{t_1\rightarrow +\infty} |\alpha(t_1)| (t_1)^{-1/2} = 0. 
\]
We recall the assumption $\|u^\alpha\|_{\dot{H}^1(\{x: |x|>r_\alpha\})} = A$, thus have 
\[
 \lim_{t_1 \rightarrow +\infty} \|\nabla u^{\alpha(t_1)}\|_{L^2(\{x: r_{\alpha(t_1)}<|x|<0.9 t_1\})} = A. 
\]
We next combine \eqref{concentration of energy initial w} with \eqref{asymptotic decomposition} to conclude 
\begin{equation*} 
 \lim_{t_1 \rightarrow +\infty} \|\nabla u(\cdot,t_1) - \nabla u^{\alpha(t_1)}\|_{L^2 (\{x: r_{\alpha(t_1)}<|x|<0.9 t_1\})} = 0.
\end{equation*}
Finally we combine the two limits above and obtain 
\[
 \lim_{t_1 \rightarrow +\infty} \|u(\cdot,t_1)\|_{\dot{H}^1(\{x: r_{\alpha(t_1)}<|x|<0.9 t_1\})} = A.
\]
This contradicts with the uniform boundedness assumption on the norm of $u$, thus finishes the proof. 
\end{proof}

\section{Defocusing Equations}
We call a non-linear wave equation 
\[
 \partial_t^2 u - \Delta u = F(x,u)
\]
defocusing if the nonlinear term $F(x,u)$ satisfies
\[
 F(x,u) u \leq 0, \qquad \forall (x,u) \in \Rm^3 \times \Rm. 
\]
If the initial data $(u_0, u_1) \in \dot{H}^1 \times L^2(\Rm^3)$, then the energy defined below is a conserved quantity 
\[
 E(u,u_t) = \int_{\Rm^3} \left(\frac{1}{2}|\nabla u(x,t)|^2 + \frac{1}{2}|u_t(x,t)|^2 + V(x,u(x)) \right) dx = \hbox{Const}. 
\]
Here the potential is defined by 
\[
 V(x,u) = -\int_0^u F(x,v) dv \geq 0.
\]
Please note that our assumption on $F$ implies 
\[
 |V(x,u)| \leq \frac{\gamma}{6} |u|^6\qquad \Rightarrow \qquad \int_{\Rm^3} V(x,u(x)) dx \lesssim_\gamma \|u\|_{\dot{H}^1(\Rm^3)}^6.  
\]
We first consider the global behaviour of the radial non-radiative solutions to a defocusing equation. We recall that $\{(u^\alpha, R_\alpha)\}$ is the one-parameter family of non-radiative solutions defined earlier in this work. 

\begin{lemma} \label{lower bound of defocusing}
 Assume that the nonlinear term is defocusing. Then given any $\alpha \neq 0$ and $A>0$, there exists a radius $R \gtrsim_1 \alpha^2 A^{-2}$, so that $R>R_\alpha$ and 
 \[
  \|u^\alpha\|_{\dot{H}^1(\{x: |x|>R\})} = A. 
 \]
\end{lemma}
\begin{proof}
Without loss of generality we assume $\alpha >0$. We recall that $w(r) = r u^\alpha (r)$ satisfies the elliptic equation $-w_{rr} = r F(r,w/r)$ with $(w,w_r)\rightarrow (\alpha, 0)$ as $r\rightarrow +\infty$. Thus we have
\begin{equation} \label{wr in term of w}
 w_r = \int_r^\infty \tau F(\tau, w(\tau)/\tau) d\tau. 
\end{equation} 
These give the asymptotic behaviour of $w$ near infinity. 
\begin{align*}
 &w(r) \geq \alpha,& &\partial_r w(r) \leq 0,& &r\gg 1.&
\end{align*}
We claim that these inequalities hold for all $r > R_\alpha$. (i.e. the solution $u^\alpha$ is still meaningful) This follows a continuity argument. If there existed a radius $r_1 > R_\alpha$ so that $w(r_1) < \alpha$, we would be able to find a number $\beta < \alpha$ but $\beta > \max\{w(r_1),0\}$. We then consider the set $\{r\geq r_1: w(r) = \beta\}$. By the intermediate value theorem, this set is non-empty. We choose 
\[
 r_2 = \max\{r\geq r_1: w(r) = \beta\}. 
\]
The continuity and asymptotic behaviour of $w$ guarantees that $r_2$ is well-defined. By continuity we must have $w(r) > \beta > 0$ for all $r>r_2$. We recall \eqref{wr in term of w} and obtain 
\[
 w_r (r) \leq 0, \; \forall r\geq r_2; \qquad \Rightarrow \qquad w(r_2) \geq \alpha.
\]
This is a contradiction. Thus we always have $w(r) \geq \alpha$. The inequality $w_r (r) \leq 0$ immediately follows \eqref{wr in term of w}. Finally we have 
\[ 
  u_r^\alpha = \frac{w_r}{r} - \frac{w}{r^2} \leq \frac{-\alpha}{r^2}. 
\]
Thus 
\[
 \|u^\alpha\|_{\dot{H}^1(\{x: |x|>R\})} \gtrsim_1 \alpha R^{-1/2}, \qquad r>R_\alpha. 
\]
This finishes the proof. 
\end{proof}

\begin{proof}[Proof of Corollary \ref{main corollary defocusing}]
 Any solution $u$ to (CP1) satisfies 
 \[
  \sup_{t\in (-T_-,T_+)} \left\|(u(\cdot,t), u_t(\cdot,t))\right\|_{\dot{H}^1 \times L^2(\Rm^3)}\leq (2E)^{1/2}. 
 \]
Lemma \ref{lower bound of defocusing} implies that given any $\alpha \neq 0$, there exists a radius $r_\alpha > R_\alpha$ so that $r_\alpha \gtrsim_1 |\alpha|^2/E$ and the non-radiative solution $u^\alpha$ satisfies 
\[
 \|u^\alpha\|_{\dot{H}^1(\{x: |x|>r_\alpha\})} = 2E^{1/2}.
\]
Finally we apply Theorem \ref{main global} to finish the proof. 
\end{proof}

\section*{Acknowledgement}
The authors are financially supported by National Natural Science Foundation of China Project 12071339.


\begin{thebibliography}{99}
 \bibitem{newradiation} R. C\^{o}te, and C. Laurent. {``Concentration close to the cone for linear waves.''} \textit{arXiv preprint} 2109.08434. 
 \bibitem{tkm1} T. Duyckaerts, C.E. Kenig, and F. Merle. {``Universality of blow-up profile for small radial type II blow-up solutions of the energy-critical wave equation.''} \textit{The Journal of the European Mathematical Society} 13, Issue 3(2011): 533-599.
  \bibitem{se} T. Duyckaerts, C.E. Kenig, and F. Merle. {``Classification of radial solutions of the focusing, energy-critical wave equation.''} \textit{Cambridge Journal of Mathematics} 1(2013): 75-144.
 \bibitem{dkm2} T. Duyckaerts, C.E. Kenig, and F. Merle. {``Scattering for radial, bounded solutions of focusing supercritical wave equations.''} \textit{International Mathematics Research Notices} 2014:  224-258.
 \bibitem{dkm3} T. Duyckaerts, C.E. Kenig, and F. Merle. {``Scattering profile for global solutions of the energy-critical wave equation.''} \textit{Journal of European Mathematical Society} 21 (2019): 2117-2162.
 \bibitem{oddhigh} T. Duyckaerts, C. E. Kenig, and F. Merle. {``Soliton resolution for the critical wave equation with radial data in odd space dimensions.''}  arXiv preprint 1912.07664.
\bibitem{radiation1} F. G. Friedlander. {``On the radiation field of pulse solutions of the wave equation.''}  \textit{Proceeding of the Royal Society Series A} 269 (1962): 53-65.
\bibitem{radiation2} F. G. Friedlander. {``Radiation fields and hyperbolic scattering theory.''} \textit{Mathematical Proceedings of Cambridge Philosophical  Society} 88(1980): 483-515.
\bibitem{strichartz} J. Ginibre, and G. Velo. {``Generalized Strichartz inequality for the wave equation.''} \textit{Journal of Functional Analysis} 133(1995): 50-68.
\bibitem{loc1} L. Kapitanski. {``Weak and yet weaker solutions of semilinear wave equations''} \textit{Communications in Partial Differential Equations} 19(1994): 1629-1676.
 \bibitem{channel5d} C. E. Kenig, A. Lawrie, B. Liu and W. Schlag. {``Relaxation of wave maps exterior to a ball to harmonic maps for all data''} \textit{Geometric and Functional Analysis} 24(2014): 610-647.
\bibitem{kenig} C. E. Kenig, and F. Merle. {``Global Well-posedness, scattering and blow-up for the energy critical focusing non-linear wave equation.''} \textit{Acta Mathematica} 201(2008): 147-212.
 \bibitem{shenradiation} L. Li, R. Shen and L. Wei.  {``Explicit formula of radiation fields of free waves with applications on channel of energy''}, to appear in \textit{Analysis and PDE}. 
 \bibitem{nonradialCE} L. Li, R. Shen, C. Wang and L. Wei. {``Asymptotic behaviour of non-radiative solution to the wave equations.''} \textit{arXiv} 2201.02286.
 \bibitem{ls} H. Lindblad, and C. Sogge. {``On existence and scattering with minimal regularity for semi-linear wave equations''} \textit{Journal of Functional Analysis} 130(1995): 357-426.
 \bibitem{enscatter1} K. Nakanishi. {``Unique global existence and asymptotic behaviour of solutions for wave equations with non-coercive critical nonlinearity.''} \textit{Communications in Partial Differential Equations} 24(1999): 185-221.
 \bibitem{enscatter2} K. Nakanishi. {``Scattering theory for nonlinear Klein-Gordon equations with Sobolev critical power.''} \textit{International Mathematics Research Notices} 1999, no.1: 31-60.
 \bibitem{ss2} J. Shatah, and M. Struwe. {``Well-posedness in the energy space for semilinear wave equations with critical growth''} \textit{International Mathematics Research Notices} 7(1994): 303-309.
 \bibitem{shen2} R. Shen. {``On the energy subcritical, nonlinear wave equation in $\Rm^3$ with radial data''}  \textit{Analysis and PDE} 6(2013): 1929-1987.
 \bibitem{struwe} M. Struwe. {``Globally regular solutions to the $u^5$ Klein-Gordon equation.''} \textit{Annali della Scuola Normale Superiore di Pisa - Classe di Scienze} 15(1988): 495-513.
\end{thebibliography}
\end{document}